\documentclass[11pt]{article}
\usepackage{dsfont}
\usepackage{enumitem}
\RequirePackage{amsthm,amsmath,amsfonts,amssymb}

\newtheorem{theorem}{Theorem}[section]

\newtheorem{lemma}[theorem]{Lemma}
\newtheorem{remark}[theorem]{Remark}
\newtheorem{proposition}[theorem]{Proposition}
\newtheorem{corollary}[theorem]{Corollary}
\newtheorem{definition}[theorem]{Definition}
\newcommand{\dif}{\mathop{}\!\mathrm{d}}
\newcommand{\R}{\mathbb{R}}
\newcommand{\pr}{\mathbb{P}}

\newcommand{\Poi}{\mathrm{Poi}}
\newcommand{\Bin}{\mathrm{Bin}}
\newcommand{\E}{{\mathbb{E}}}
\newcommand{\N}{{\mathbb{N}}}

\newcommand{\erfc}{\operatorname{erfc}}
\newcommand{\ind}{\mathds{1}}
\newcommand{\diam}{\mathrm{diam}}

\newcommand{\gap}{\mathrm{gap}}
\newcommand{\Card}{\mathrm{Card}}
\newcommand{\n}{\mathbf{n}}
\newcommand{\T}{\mathcal{T}}

\usepackage{amsmath} 
\usepackage{amssymb}
\usepackage{hyperref} 
\usepackage{graphicx}
\begin{document}

\title{Singular Birth-Death model}

\date{}

\author{Ruibo Kou \thanks{Mathematics Institute, University of Warwick, Coventry, CV4 7AL, UK, \texttt{Ruibo.Kou@warwick.ac.uk}}}

\maketitle
\begin{abstract}
We introduce and study a one-dimensional singular birth-death model constructed from It\^o's Brownian excursion measure. In this model, particles are born at the origin at infinite rate, evolve as Brownian excursions away from the origin, and are killed upon returning to the origin. Although the number of particles alive at any fixed deterministic time is almost surely finite, we prove that the set of exceptional times at which infinitely many particles are alive is almost surely dense and has Hausdorff dimension $1/2$. We further investigate the spatial structure of the particle system at such exceptional times. In particular, we identify the critical near-origin profile of the infinite particle cloud, which is governed by a scale of order $\log\log(1/|x|)^{-1}$ as $x\downarrow0$. Finally, we establish laws of large numbers and central limit theorems for two natural observables of the system: the $\alpha$-power sum of the positions of all alive particles and the top particle position.
\end{abstract}
\noindent
\section{Introduction and main results}\label{Intro}
\subsection{Introduction}
The classical birth and death processes introduced by Feller in \cite{feller2,feller1} have a  long history as model of population dynamics. It describe changes in population size through births and deaths.
The addition of spatial motion to birth-death models, allowing particles to move as well as undergo birth and death, is a comparatively recent development. A general framework for combining Markovian particle motions with branching was developed by Ikeda, Nagasawa and Watanabe \cite{ClassicalBranching} in the 1960s. The model and its various extensions has wide applications such as the mathematical study of cancer evolution \cite{cancerapplication} and the convergence of diffusion Monte Carlo algorithm \cite{Fan}. The interplay between birth and death rate has particular interest when both birth and death become large. In this paper we introduce a singular birth-death model arising as the infinite rate limit of a one dimensional particle system in which births and deaths occurs at the origin.

To motivate the construction, consider a Poisson point process of birth times with rate $\mu>0$. At each birth time, a new particle is created at the origin and evolves independently as a Brownian motion started from zero. Each particle is killed when its accumulated local time at the origin exceeds an independent exponential clock with rate $\lambda>0$. Therefore each particle is created at rate $\mu$ in time and killed at rate $\lambda$ with respect to its local time at zero. At this stage, both birth and death rates are finite and the system consists of independent Brownian particles with singular birth and death position.

We are interested in the limit in which both the birth and death rates tend to infinity. More precisely, we let
\begin{equation*}
    \mu,\lambda\to\infty,\quad\frac{\mu}{\lambda}=\rho>0.
\end{equation*}
Consequently, particles are killed increasingly rapidly after returning to the origin, while the increasing birth rate compensates for this loss, yielding a nontrivial limit. The resulting limit is a particle system in which particles are born at the origin at infinite rate and evolve along Brownian excursions away from the origin.

One reason to study such a limit is that taking a limit can reveal a simpler mathematical object, as is familiar from the study of scaling limits. The same idea motivates the high rate limit considered here. When the birth and killing rates tend to infinity with their ratio fixed, the limiting system admits an explicit representation as a Poisson point process whose intensity measure is characterized by Brownian excursion measure.

A heuristic motivation for the model comes from the infrared problem in quantum field theory. In models describing the interaction of matter with massless bosonic fields, a particle may be accompanied by a cloud of arbitrarily soft bosons, whose number diverges as an infrared cutoff is removed. This phenomenon has been studied in models involving massless scalar bosons, beginning with the work of Nelson \cite{Nelson} and Fr\"ohlich \cite{Frohlich}, and is closely related to the soft-photon problem in quantum electrodynamics. (See \cite{Kulish} as an example.) In particular, Hirokawa \cite{hirokawa} established rigorous results on the divergence of the soft-boson number and the absence of ground states in the Nelson model. Inspired by the picture of an infinitely massive particle acting as a fixed source, we consider a simplified classical system in which independent Brownian particles are produced at the origin and absorbed through a mechanism governed by their local time at the origin. In the singular limit, both the birth rate and the rate of absorption with respect to local time become infinite. We investigate the exceptional times at which infinitely many particles coexist, even though their number is almost surely finite at every deterministic time. Our model is not derived from quantum field theory, and this analogy serves only as an informal motivation for studying the emergence of an infinite particle population.

The Brownian net and the Brownian fan provide further motivation for this infinite rate limit.
The Brownian net, introduced by Sun and Swart in \cite{BN}, is an interacting particle system of one dimensional branching-coalescing Brownian motions arsing from the scaling limit of branching-coalescing random walks on a discrete space-time lattice. Under this scaling limit, the branching rate tends to infinity. The Brownian net can also be constructed from the Brownian web by introducing branching through a Poisson marking of special points, as developed by Newman, Ravishankar and Schertzer in \cite{BW}. In this construction, paths may branch from an existing trajectory with infinite branching rate as the intensity measure of the Poisson marking is infinite, and subsequently coalesce with them. Our construction isolates a similar branching and coalescing mechanism by considering a single distinguished 'parent' particle fixed at the origin, from which particles are born and to which they return. We retain only the offspring produced by this parent particle and do not introduce interactions between distinct offspring. In particular, when an offspring returns to the origin in the limiting system, it coalesces with the parent particle.

A more direct connection is the Brownian fan introduced by Hairer and Weare in \cite{Fan}. The Brownian fan is constructed through a recursive branching mechanism in which each particle produces offspring along its trajectory according to the Poisson point process of the type  described above. (See Section \ref{Brownian fan} or \cite[Definition 2.9]{Fan} for a precise definition.) Each offspring can in turn produce further generations, resulting a system with an infinite rate branching-killing mechanism. Our singular birth-death model corresponds to the first generation dynamics of the Brownian fan: in the Brownian fan if we fix the parent trajectory at the origin and retain only its direct offspring, the resulting system is our singular birth-death model. This connection is particularly relevant to the study of exceptional times at which infinitely many particles are alive.

\begin{figure}[htpb]
\centering
\includegraphics[scale=0.8]{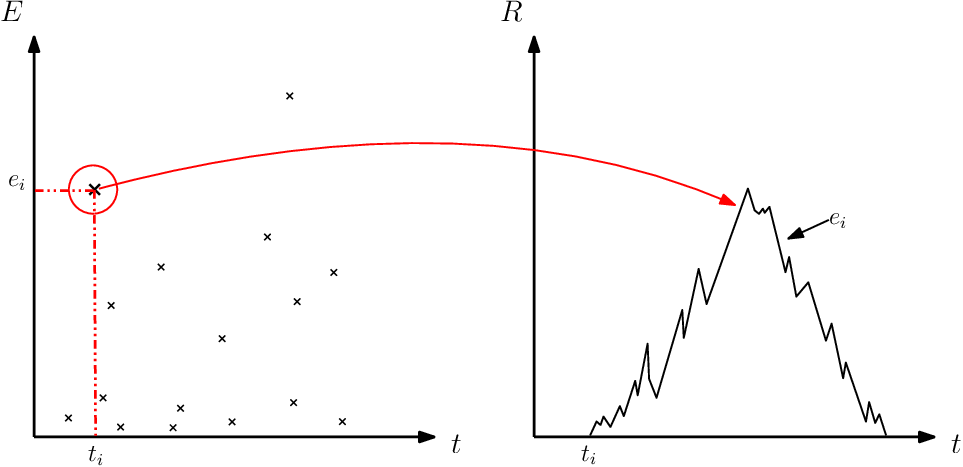}
\caption{A map from points in $\R^+\times\mathcal{E}$ to paths on $\R^+\times \R$}
\label{PPP}
\end{figure}

The singular birth-death model admits a natural description in terms of It\^o's Brownian excursion measure $\n$. It is initiated by Itô in \cite{Ito_excursion}, which provides a powerful framework for decomposing recurrent Markov processes into excursions away from a distinguished recurrent state. Our construction uses the same excursion measure in a different way. Rather than placing excursions on the local time axis and concatenating them to form a single path, we place excursions on the time axis and interpret each excursion as the trajectory of a distinct particle. More precisely, the limiting system is described by a Poisson point process 
\begin{equation*}
    \Theta=\sum_{i\geq 1}\delta_{(t_i,e_i)}
\end{equation*}
on the product space of the time axis $\R^+$ and the excursion space with intensity measure
\begin{equation*}
    \rho\dif t\otimes \n.
\end{equation*}
Each atom $(t_i,e_i)$ represents a particle born at time $t_i$ whose trajectory is given by
\begin{equation*}
    x_i(t)=e_i(t-t_i),
\qquad
t_i\leq t<t_i+\xi(e_i).
\end{equation*}
The particle is removed at time $t_i+\xi(e_i)$, where $\xi(e_i)$ denotes the life time of the excursion. (See Fig \ref{PPP} for an illustration of this construction.)

 Our first object is to investigate the exceptional time phenomenon quantitatively. Let $X_t$ denote the point process for the particles alive at time $t$ and let $N_t= X_t(\R)$ be their number. We first show that although $N_t$ is almost surely finite at every deterministic time, the random set
 \begin{equation}\label{set of excptional times}
     \T=\{\tau\in\R^+:N_\tau=\infty\}
 \end{equation}
 is almost surely dense in $\R^+$. Furthermore we show that the Hausdorff dimension 
 \begin{equation*}
     \dim \T=\frac{1}{2}\quad\text{almost surely.}
 \end{equation*}

Hairer and Weare \cite[Section 2]{Fan} establish the existence of exceptional times for the Brownian fan. As a corollary of our dimension result, we show that the set of exceptional times at which infinitely many particles alive in the Brownian fan has Hausdorff dimension $1/2$ almost surely. Moreover, for our singular birth-death model, we prove that the origin is the unique accumulation point of the particle positions at every exceptional time and that the corresponding set of exceptional time-space points has Hausdorff dimension one. These results provide a quantitative description of the exceptional times arising from the first generation excursion mechanism associated with the Brownian fan.

Our second object is to investigate the spatial distribution of the infinitely many particles present at exceptional times. Let $X_t$ be the point process of particles alive at time $t$. Since the origin is their only possible accumulation point, it is natural to study the particle configuration
\begin{equation*}
    \langle X_t,f\rangle=\int_\R f(x)X_t(\dif x)
\end{equation*}
for continuous functions vanishing at the origin. This functional captures the accumulation of particles near the origin at exceptional times. We establish a necessary and sufficient criterion for the existence of exceptional times at which this quantity is infinite and determine the Hausdorff dimension of the corresponding set of times, which is equal to $1/2$ almost surely. In particular, our results identify a critical near-origin scale of order
\begin{equation*}
    \left(\log\log(1 /|x|)\right)^{-1},
\end{equation*}
providing a detailed spatial description of  the infinite particle cloud.

Finally, we study the long-time behavior of two natural observables: the sum of the $\alpha$-powers of the absolute positions of all particles alive at time $t$, and the position of the top particle. We establish laws of large numbers and central limit theorems for both observables. Together, these results characterize the exceptional time behavior and the long time evolution of the singular birth-death model.

The rest of the paper is organized as follows. In the rest of Section \ref{Intro} we introduce necessary notations and state our main results rigorously. In Section \ref{Construction of the model} and \ref{convergence of n^lambda} we construct the singular birth-death model as the limit of particle systems with finite birth-death rate. This is achieved by establishing the convergence
\begin{equation*}
    \lambda\n^\lambda\to\n,\quad\lambda\to\infty
\end{equation*}
 under the $M_0$ topology, (see \cite{M0convergence1,M0convergence2} or Definition \ref{M_0 topology} for the definition of $M_0$ convergence) where $\n^\lambda$ is the law of Brownian motion absorbed at zero with rate $\lambda$. In Section \ref{proof of theorem 1 and 2} we prove our Hausdorff dimension results for the sets of exceptional times by constructing a Cantor-like subset. Section \ref{proof of theorem 3} and \ref{proof of theorem 4} are devoted to the proofs of law of large numbers and central limit theorem for the sum of $\alpha$-powers and positions of the top particles respectively.

\subsection{Main results}
We starts by introducing some notations. For any $e\in C(\R^+)$ Let
\begin{equation*}\label{Excursion space}
    \xi(e):=\inf\{s\in\R^+:e(t)=0\text{ for all }t\geq s\}\end{equation*}
be the life time of path $e$ and
\begin{equation}\label{Excursion space} 
    \mathcal{E}:=\{e\in C(\R^+):0<\xi(e)<\infty\}.
\end{equation}
be the space of excursions. The It\^o's Brownian excursion measure is a measure $\n$ on $\mathcal{E}$ such that for any $0<t_1<\cdots<t_k$ we have
\begin{equation*}
\begin{split}
    &\n\left\{e\in \mathcal{E}:e(t_j)\in\dif x_j,\ j=1,\cdots,k,\ t_1<\cdots<t_k,\ \xi(f)\geq t_k\right\}\\&=\n_{t_1}(\dif x_1)P_{t_2-t_1}^0(x_1,\dif x_2)\cdots P^0_{t_k-t_{k-1}}(x_{k-1},\dif x_k),
    \end{split}
\end{equation*}
where 
\begin{equation}\label{transition density of killed Brwonian motion}
    P_t^0(x,\dif y)=\frac{1}{\sqrt{2\pi t}}\left(\exp\left(-\frac{(x-y)^2}{2t}\right)-\exp\left(-\frac{(x+y)^2}{2t}\right)\right)\ind_{\{xy>0\}}\dif y
\end{equation}
is the transition density of the Brownian motion killed at zero and
\begin{equation}\label{entrance law}
    \n_t(\dif x)=\frac{|x|}{\sqrt{2\pi t^3}}\exp\left(\frac{-x^2}{2t}\right)\dif x\,\quad \text{for }x\in\R,
\end{equation}
is called the entrance law of $\n$. (See \cite{Ito_excursion} for the derivation of the entrance law.) 
Let $\Theta$ be a Poisson point process on $\R^+\times \mathcal {E}$ with intensity measure 
\begin{equation*}
    \rho\dif t\otimes\n.
\end{equation*}
The intensity measure is sigma-finite, so $\Theta$ is almost surely a countable simple point process. We may therefore write
\begin{equation*}
    \Theta=\sum_{i=1}^\infty \delta_{(t_i,e_i)}.
\end{equation*}
\begin{definition}[Singular Birth-Death Model]\label{def of sbd}
    The singular birth-death model with parameter $\rho$ is a point process 
    \begin{equation}\label{def of X_t}
         X_t=\sum_{i=1}^\infty \delta_{e_i(t-t_i)}\ind_{\{t_i\leq t\leq t_i+\xi(e_i)\}}=\sum_{i=1}^{N_t}\delta_{x_i(t)},
    \end{equation}
    where $N_t=X_t(\R)$ is the number of particles alive at time $t$ and
$x_i(t)=e_i(t-t_i)$ denotes the position at time $t$ of the $i$-th particle alive at that time. 
\end{definition}
Thus, $X_t$ records the positions of all particles alive at time $t$. Our first result describes the size of $\T$ defined by (\ref{set of excptional times}) and the spatial configuration of the system at exceptional times.
\begin{theorem}\label{exceptional times}
Almost surely, the following statements hold.
    \begin{enumerate}[label=(\roman*)]
    \item The set $\T$ is dense in $\R^+$.
    \item The Hausdorff dimension of $\T$ is
    \begin{equation*}
        \dim(\T)=\frac{1}{2}.
    \end{equation*}
    \item At each exceptional time $\tau\in \T$, zero is the unique cluster point of the positions of particles.
    \item The set of time-space exceptional points 
    \begin{equation*}
        \mathcal{S}:=\{(\tau,x_i(\tau)):\tau\in \T,\ i\in\N\}
    \end{equation*}
    has Hausdorff dimension $1$.
    \end{enumerate}
\end{theorem}
We next investigate the spatial distribution of particles near the origin at exceptional times. For a measure $\nu$ on $\R$ and a nonnegative measurable function $f\colon\mathbb R\to[0,\infty]$, write
\begin{equation*}
   \langle\nu,f\rangle:=\int_\R f(x)\nu(\dif x).
\end{equation*}

\begin{theorem}\label{loglog theorem}
    Let $f\colon\mathbb R\to[0,\infty)$ be an even continuous function that is non-decreasing on $[0,\infty)$.
    \begin{enumerate}[label=(\roman*)]
    \item The process $\left(\langle X_t,f\rangle\right)_{t\in\R^+}$ is a discontinuous process if and only if 
    \begin{equation}\label{requirement of test function}
        \int_0^\epsilon \frac{f(x)}{x\log(1/x)}\dif x=\infty\quad\text{for some }0<\epsilon<1.
    \end{equation}
    Moreover, almost surely there exists $\tau\in \T$ such that $\langle X_\tau,f\rangle=\infty$ if and only if $f$ satisfies (\ref{requirement of test function}).
    \item Suppose $f$ satisfies (\ref{requirement of test function}), the set
    \begin{equation*}
        \T(f):=\{\tau\in \T:\langle X_\tau,f\rangle=\infty\}
    \end{equation*}
    has Hausdorff dimension $1/2$ almost surely.
    \end{enumerate}   
\end{theorem}
\begin{remark}

    By Theorem \ref{exceptional times} $(iii)$, zero is the only accumulation point of the particle configuration at an exceptional time. Consequently, the criterion in Theorem \ref{loglog theorem} depends only on the behavior of $f$ near zero. Moreover, since $\n$ is symmetric, the configurations of particles on the positive and negative half spaces have the same distribution. It is therefore natural to restrict attention to even functions.

    For an even continuous function $g_p$ such that $g_p$ is non decreasing on $[0,\infty)$. If 
    \begin{equation*}
        g_p\ind_{\{|x|<\exp(-\exp(1))\}}=(\log\log(1/|x|))^{-p}.
    \end{equation*}   
     Then  it is easy to check that $g_p$ satisfies (\ref{requirement of test function}) if and only if $p\leq1$. Hence we get the smallest space that $X_t$ lives in is
    \begin{equation*}
        \bigcap_{p>1}\{\nu:\langle\nu,g_p(x)\rangle<\infty\}.
    \end{equation*}
\end{remark}
The next two results are about the long term behavior for $X_t$ with two natural observables of the system. The first is the sum of the $\alpha$-powers of the absolute positions of all particles alive at time $t$. Let $f_\alpha(x)=|x|^\alpha$ for $\alpha\in(0,1)$.
We set
\begin{equation}\label{def Z_t^alpha}
    Z_t^\alpha:=\langle X_t,f_\alpha\rangle=\sum_{i=1}^{N_t}|x_i(t)|^\alpha.
\end{equation}
Note that for any $\alpha>0$, $f_\alpha$ fails to satisfy (\ref{requirement of test function}) hence $(Z_t^\alpha)_{t\in\R^+}$ is a well defined continuous process.
\begin{theorem}\label{Z_t^alpha}
    Let $Z_t^\alpha$ be defined by (\ref{def Z_t^alpha}).
    \begin{enumerate}[label=(\roman*)]
    \item Almost surely, $Z_t^\alpha$ is locally $\gamma-$H\"older continuous  for all $\gamma\in(0,\alpha/2)$.
    \item We have law of large numbers for $Z_t^\alpha$:
    \begin{equation*}
        \frac{Z_t^\alpha}{t^{(1+\alpha)/2}}\to D_\alpha,\ a.s.,
    \end{equation*}
    where $D_\alpha=\rho2^{(3+\alpha)/2}\int_0^\infty x^\alpha\erfc(x)\dif x<\infty$.
    \item We have central limit theorem for $Z_t^\alpha$:
    \begin{equation*}
        \frac{Z^\alpha_t-D_\alpha t^{(1+\alpha)/2}}{\sqrt{\widetilde D_\alpha t^{(1+2\alpha)/2}}}\to N(0,1)\text{ in distribution},
    \end{equation*}
    where $N(0,1)$ is a standard normal random variable and 
    \begin{equation}\label{def of D_alpha}
        \widetilde D_\alpha=\rho 2^{(3+2\alpha)/2} \int_0^\infty x^{2\alpha}\erfc(x)\dif x
    \end{equation}
    is a finite constant.
    \end{enumerate}
\end{theorem}
Our final result concerns the position of the the top particle. To formulate it, for any $t\in\R^+$ define
\begin{equation}\label{def of M_t}
    M_t:=\max_{1\leq i\leq N_t}\{x_i(t)\}\vee 0
\end{equation}
The convention $M_t=0$ is used when there is no particle on the positive half-line. Let
\begin{equation*}
    A(y):=\int_y^\infty \erfc(x)\dif x.
\end{equation*}
Since $A$ is continuous and strictly decreasing it admits an inverse $A^{-1}$.
Let 
\begin{equation*}
    u_t:=A^{-1}\left(\frac{1}{\sqrt{2t}}\right).
\end{equation*}
Now we can state the theorem about $M_t$:
\begin{theorem}\label{theorem of M_t}
    Let $M_t$ be defined by (\ref{def of M_t}).
    \begin{enumerate}[label=(\roman*)]
    \item Almost surely, $M_t$ is locally $\gamma-$H\"older continue for all $\gamma\in(0,1/2)$.
    \item We have law of large numbers for $M_t$:
    \begin{equation*}
        \frac{M_t}{\sqrt{t\log t}}\to 1,\ a.s.
    \end{equation*}
    \item We have central limit theorem for $M_t$:
    \begin{equation*}
        \frac{M_t-\sqrt{2t}u_t}{v_t}\to G\text{ in distribution},
    \end{equation*}
    where $v_t=\sqrt{t}/\sqrt{2}u_t$ and $G$ is a shifted Gumbel random variable whose distribution function is given by 
    \begin{equation*}
        \pr[G\leq x]=\exp\left(-\exp\left(-(x-\log \rho)\right)\right).
    \end{equation*}
    \end{enumerate}
\end{theorem}
\begin{remark}
By symmetry, the corresponding results hold for the bottom particle. More precisely, if $\widetilde{M_t}:=\min_{1\leq i\leq N_t}\{x_i(t)\}\wedge 0$ then $(-\widetilde M_t)_{t\geq 0}$ satisfies the same conclusions as $(M_t)_{t\geq 0}$.

The term $\sqrt{2t}u_t\sim \sqrt{t\log t}$ for large $t$ is asymptotically equal to $\E[M_t]$. The term $v_t\sim\sqrt{t/\log t}$ is asymptotically equal to the standard deviation of $M_t$.
\end{remark}

We conclude the introduction by fixing some notations. In the rest of this paper, for any two functions $f$ and $g$ from $\R$ to $\R$ and $x_0\in\overline{\R}$ the statement
\begin{equation*}
    f(x)\approx g(x),\ \text{as }x\to x_0
\end{equation*}
is defined to be
\begin{equation*}
    f(x)=O(g(x))\ \text{and }g(x)=O(f(x))\ \text{as }x\to x_0.
\end{equation*}
We also use
\begin{equation*}
    a_n\approx b_n\ \text{as }n\to\infty
\end{equation*}
for two real valued sequences with similar meaning.
Moreover, for any set $A\subseteq\R^n$ we use 
\begin{equation*}
    \dim(A),\quad \diam(A)\quad \Card(A)\quad \text{and}\quad |A|
\end{equation*}
to represent the Hausdorff dimension, diameter, cardinality and Lebesgue measure of $A$ respectively. We also use
\begin{equation*}
    \Poi(\lambda)\quad \mathrm{Exp}(\lambda)\quad\Bin(n,p_0)
\end{equation*}
To represent Poisson, exponential and binomial random variables with parameter $\lambda$, $(n,p_0)$ respectively whereas we use $\exp(x)$ to represent exponential function. Finally $C$ denotes a finite positive constant whose value may change from line to line and any subscripts indicate its dependence on parameters.

\section{Construction of the model}\label{Construction of the model}\ 

In this section, we justify the interpretation of $\Theta$ as a birth-death particle system with infinite birth and killing rates. We first construct a particle system with finite birth rate and finite killing rate, where killing occurs at a constant rate with respect to the local time accumulated at the origin. We then let both rates tend to infinity while keeping their ratio fixed.

Let $(B_t)_{t\geq0}$ be a standard Brownian motion and let $L_t^0$ denote its local time at zero. For $\lambda>0$, let $\mathrm{Exp}(\lambda)$ be an exponential random variable that is independent of the Brownian motion, and define
\begin{equation*}
    B_t^\lambda=\begin{cases}
        B_t&\quad\text{if }L_t^0\leq\text{Exp}(\lambda)\\
        0&\quad\text{if }L_t^0>\text{Exp}(\lambda) 
    \end{cases}.
\end{equation*}
Thus $B_t^\lambda$ is a particle absorbed at zero with rate $\lambda$. The distribution of $(B^\lambda_\cdot)$ is denoted by $\n^\lambda$, which a measure supported on $\mathcal{E}$. For any pair of positive real numbers $(\mu,\lambda)$, let $\Theta^{\mu,\lambda}$ be a Poisson point process on $\R^+\times \mathcal{E}$ with intensity measure
\begin{equation*}
    \mu\dif t\otimes\n^\lambda.
\end{equation*}
As in the Definition \ref{def of sbd}, each atom $(t_i,e_i)\in\Theta^{\mu,\lambda}$ represents a particle born at time $t_i$ with $\{e(s-t_i):{s\in[t_i,t_i+\xi(e)]}\}$ being its trajectory. 

\begin{definition}[Finite rate singular birth-death model]
    The singular birth-death model with finite birth-death rate $(\mu,\lambda)$ is a point process
    \begin{equation*}
        X^{\mu,\lambda}_t=\sum_{i\geq 1}\delta_{e_i(t-t_i)}\ind_{\{t_i\leq t\leq t_i+\xi(e_i)\}}=\sum_{i=1}^{N^{\mu,\lambda}_t}\delta_{x^{\mu,\lambda}_i(t)},
    \end{equation*}
where $N_t^{\mu,\lambda}=X_t^{\mu,\lambda}(\R)$ is the number of particles alive at time $t$ and $x^{\mu,\lambda}_i(t)$ denotes the potion at time $t$ of $i$-the particle alive at that time.
\end{definition}
Since $\n^\lambda$ is a probability measure, the projection of $\Theta^{\mu,\lambda}$ onto the time coordinate is a Poisson point process on $\R^+$ with intensity $\mu\dif s$. Thus particles are born at rate $\mu$ and each particle is killed at rate $\lambda\dif L^0_t$. Fix $\rho>0$ and set
\begin{equation*}
    \rho=\frac{\mu}{\lambda}.
\end{equation*}
We now let $\lambda\to\infty$. In this regime, both the birth and death rates tend to infinity.  

The key point for showing our finite rate singular birth-death model converges to the (infinite rate) singular birth-death model is to show that $\lambda\n^\lambda\to\n$ under a suitable topology on the space of measures supported on $\mathcal{E}$. The main issue is that $\n(\mathcal{E})=\infty$ while $\n(\{e\in\mathcal{E}:\|e\|_\infty>\epsilon\})<\infty$ for all $\epsilon>0$. Therefore we need to avoid $0$ in $\mathcal{E}$, which leads as to the $M_0$ topology introduced by Hult and Lindskog in \cite{M0convergence1}.

Let $\mathbf S$ be a separable Banach space, and let $\mathbf S_0=\mathbf S\setminus\{0\}$ be equipped with the relative topology and the $\sigma$-algebra $\mathcal{F}_0=\{A:A\subseteq \mathbf S_0,A\in\mathcal F\}$, where $\mathcal{F}$ is the Borel $\sigma$-algebra on $\mathbf S$. Let $\mathcal C_b(\mathbf S)$ denote the set of bounded and real-valued continuous functionals on $\mathbf S$. Define
\begin{equation*}
    \mathcal{C}_0:=\{F\in \mathcal C_b(\mathbf S):\text{there exits $\epsilon$ such that $F(x)=0$ for all $\|x\|<\epsilon$}\}.
\end{equation*}
We write $B_r$ for the ball centered at zero with radius $r$. Let $M_0$ be the set of Borel measures on $\mathbf S$ such that $\mathbf S\setminus B_r$ has finite measure for all $r>0$.
\begin{definition}[Convergence under $M_0$ topology]\label{M_0 topology}
    For a sequence of measures $\mu_n$ and $\mu$ in $M_0$, we say $\mu_n\to \mu$ under $M_0$ topology if 
    \begin{equation*}
        \int_\mathbf S F\dif \mu_n\to\int_\mathbf S F\dif \mu
    \end{equation*}
    for all $F\in\mathcal{C}_0$.
\end{definition}

Let $C_0(\R^+)$ be the set of continuous functions from $\R^+$ to $\R$ that vanish at infinity.  Equipped with the supremum norm, $C_0(\R^+)$ is a separable Banach space.  Both $\n$ and $\lambda\n^\lambda$ can therefore be viewed as measures supported on $C_0(\R^+)$.  (Since $(\mathcal E,\|\cdot\|)$ is not complete, we work in a larger space $C_0(\R^+)$.) In section \ref{convergence of n^lambda}, by taking $\mathbf S=C_0(\R^+)$, we shall prove that $\lambda \n^\lambda\to \n$ under $M_0$ topology and as a consequence, $\Theta^{\rho\lambda,\lambda}\to \Theta$ in distribution. To be more precise, we have the following theorem, which provides another way to construct Brownian excursion measure $\n$ by taking the limit of Brownian motion killed with finite rate.

\begin{theorem}\label{M0 convergence of lambda n lambda}
    Let $C_0(\R^+)$ be the set of continuous functions $e$ such that $\lim_{t\to\infty}e(t)=0$ and equip it with the supremum norm $\|\cdot\|$. Then for any bounded real value functional $F_\epsilon$ such that $F_\epsilon(e)=0$ for all $\|e\|<\epsilon$, we have
    \begin{equation*}
        \int_{C_0(\R^+)}F_\epsilon(e)\lambda \n^\lambda(\dif e)\to\int_{C_0(\R^+)}F_\epsilon(e)\n(\dif e).
    \end{equation*}
    In another word, $\lambda \n^\lambda \to \n$ under $M_0$ topology. As a consequence we have
    \begin{equation*}
        \Theta^{\rho\lambda,\lambda}\to\Theta,\quad\text{in distribution under $M_0$ topology}.
    \end{equation*}
\end{theorem}

In this section, we provides the intuition for the convergence for $\lambda \n^\lambda$. The one dimensional marginal of $\n^\lambda$ is given by
\begin{equation*}
    \begin{split}
        \n_t^\lambda(\dif x)&:=\pr_0[B_t^\lambda\in\dif x]=\pr_0[B_t\in\dif x,\ L^0_t\leq \text{Exp}(\lambda)]\\
        &=\int_0^\infty\int_0^u\pr_0[B_t\in\dif x,\ L^0_t\in\dif l]\lambda e^{-\lambda u}\dif l\dif u\\
        &=\int_0^\infty\int_l^\infty\pr_0[B_t\in\dif x,\ L^0_t\in\dif l]\lambda e^{-\lambda u}\dif u\dif l\\
        &=\int_0^\infty\pr_0[B_t\in\dif x,\ L^0_t\in \dif l]e^{-\lambda l}\dif l.
    \end{split}
\end{equation*}
The joint density of $(B_t,L_t^0)$ is well studied in \cite[Exercise 8.7]{localTime}, which is given by
\begin{equation*}
\begin{split}
    q_t(x,l)\dif x\dif l&:=\pr_0[B_t\in\dif x,L^0_t\in\dif l]\\&=\frac{|x|+l}{\sqrt{2\pi t^3}}\exp\left(-\frac{(|x|+l)^2}{2t}\right)\dif x\dif l,\ \text{for }l\geq0\text{ and } x\in\R.
\end{split}
\end{equation*}
Hence we get 
\begin{equation}\label{first compute of n_t^lambda}
    \n^\lambda_t(\dif x)=\int_0^\infty q_t(x,l)e^{-\lambda l}\dif l\dif x.
\end{equation}
Let $\tau_0=\inf\{t\in\R^+:B_t=0\}$ then using the strong Markov property, we can derive the transition density $P_t^\lambda(x,\dif y)$ of $(B_t^\lambda)_{t\geq0}$, which is given by
\begin{equation}\label{transition density of killed Brownian motion with rate lambda}
     P_t^\lambda(x,\dif y)=
    \begin{cases}
    \n_t^\lambda(\dif y),\quad&\text{for}\quad x=0,\\
    P_t^0(x,\dif y)+\int_0^t\pr_x[\tau_0\in\dif s]\n_{t-s}^\lambda(\dif y)\dif s,\quad&\text{for}\quad x\neq0,
    \end{cases}
\end{equation}
where $P_t^0(x,\dif y)$ is the transition density of Brownian motion killed at zero given by (\ref{transition density of killed Brwonian motion}). Therefore, we get the finite dimensional distribution of $(B^\lambda_\cdot)$:
\begin{equation*}
    \begin{split}
        &\n^\lambda(\{e\in \mathcal{E}:e(t_j)\in\dif x_j,\ j=1,\cdots,k,\ t_1<\cdots<t_k,\ \xi(e)\geq t_k\})\\
        &=\n^\lambda_{t_1}(\dif x_1)P_{s_2}^\lambda(x_1,\dif x_2)\cdots P^\lambda_{s_k}(x_{k-1},\dif x_k),
    \end{split}
\end{equation*}
where $s_j=t_j-t_{j-1}$ for $2\leq j\leq k$. By taking $\lambda\to \infty$ we get from (\ref{first compute of n_t^lambda}) that for any $t>0,x\in\R$, $\n_t^\lambda(\dif x)\to0$ as $\lambda\to\infty$. Hence according to (\ref{transition density of killed Brownian motion with rate lambda}) we have
\begin{equation*}
    P^\lambda_t(x,\dif y)\to P_t^0(x,\dif y).
\end{equation*}
Indeed, after a particle returns to zero, its remaining lifetime tends to zero as $\lambda\to\infty$. Moreover by (\ref{first compute of n_t^lambda}),
\begin{equation*}
        \lambda \n_t^\lambda(\dif x)=\int_0^\infty q_t(x,l)\lambda \exp\left(-\lambda l\right)\dif l \dif x=\E[q_t(x,\mathrm{Exp}(\lambda))]\dif x.
\end{equation*}
Since $\mathrm{Exp}(\lambda)$ converges to zero in distribution as $\lambda\to\infty$, we get 
\begin{equation*}
    \lambda \n_t^\lambda(\dif x)\to q_t(x,0)\dif x=\frac{|x|}{\sqrt{2\pi t^3}}\exp\left(\frac{-x^2}{2t}\right)\dif x,
\end{equation*}
which is the entrance law of the It\^o's Brownian excursion measure given by (\ref{entrance law}). Hence we get
\begin{equation}\label{uniequness of n}
\begin{split}
    &\lambda \n^{\lambda}(\{e\in \mathcal{E}:e(t_j)\in\dif x_j,\ j=1,\cdots,k,\ t_1<\cdots<t_k,\ \xi(e)\geq t_k\})\\ &\to \n_{t_1}(\dif x_1)P_{s_2}^0(x_1,\dif x_2)\cdots P^0_{s_k}(x_{k-1},\dif x_k)\\
    &=\n(\{e\in \mathcal{E}:e(t_j)\in\dif x_j,\ j=1,\cdots,k,\ t_1<\cdots<t_k,\ \xi(e)\geq t_k\})
    \end{split}.
\end{equation}
Therefore, the finite dimensional distribution of $\lambda \n^\lambda$ converges to $\n$, we get a unique candidate for the limit point of $\lambda \n^\lambda$.

We next recall the Brownian scaling property of $\n$, which is the key reason behind the $1/2$ Hausdorff dimension of Theorem \ref{exceptional times} and \ref{loglog theorem}. 
\begin{proposition}\label{Brownian excursion scaling}
    Let $h>0$ and $\Phi_h:\mathcal{E}\to\mathcal{E}$ be a map such that $(\Phi_he)(s)=h^{-1}e(h^2s)$ for all $s\geq0$ and $e\in\mathcal{E}$. Then $\n\circ\Phi^{-1}_h=h^{-1}\n$.
\end{proposition}
\begin{proof}
    Let $\n'=\n\circ\Phi_h^{-1}$. We firstly compute the single marginal of $\n'$. Note that
    \begin{align*}
        \n'_t(\dif x)&:=\n'\left(\{e\in\mathcal{E}: e(t)\in\dif x,\ \xi(e)\geq t\}\right)\\
        &=\n\left(\{e\in\mathcal{E}:(\Phi_he)(t)\in\dif x,\xi(\Phi_h(e))\geq t\}\right)\\
        &=\n\left(\{e\in\mathcal{E}:h^{-1}e(h^{2}t)\in\dif x,h^{-2}\xi(e)\geq t\}\right).
    \end{align*}
    Using the entrance law of $\n$ given by (\ref{entrance law}) we obtain
    \begin{equation*}
        \n'_t(\dif x)=\frac{|hx|}{\sqrt{2\pi(h^2t)^3}}h\exp\left(-\frac{x^2}{2t}\right)\dif x=h^{-1}\n_t(\dif x).
    \end{equation*}
    Therefore, we use Markov property to get
    \begin{align*}
        &\n'\left(\{e\in \mathcal{E}:e(t_j)\in\dif x_j,\ j=1,\cdots,k,\ t_1<\cdots<t_k,\ \xi(f)\geq t_k\}\right)\\
        &=\n'_{t_1}(\dif x_1)P_{s_2}^0(hx_1,\dif hx_2)\cdots P^0_{s_k}(hx_{k-1},\dif hx_k)\\
        &=h^{-1}\n_{t_1}(\dif x_1)P_{s_2}^0(x_1,\dif x_2)\cdots P^0_{s_k}(x_{k-1},\dif x_k)\\
        &=h^{-1}\n\left(\{e\in \mathcal{E}:e(t_j)\in\dif x_j,\ j=1,\cdots,k,\ t_1<\cdots<t_k,\ \xi(f)\geq t_k\}\right),
    \end{align*}
    where in the second equality we use the fact that $P^0_{h^2t}(hx,\dif hy)=P^{0}_t(x,\dif y)$ as the distribution of killed Brownian motion is invariant under the Brownian scaling.
\end{proof}

We conclude this section by determining the distribution of the particle configuration $X_t$ given by  (\ref{def of X_t}) at a fixed deterministic time. Throughout this paper for any simple measure $\mathbf{X}$ we write $x\in \mathbf{X}$ to mean that $x$ is an atom of the measure.
\begin{proposition}\label{Poi X_t}
For each deterministic time $t\in\R^+$, $X_t$ is a Poisson point process with intensity measure $\Lambda_t(\dif x)=\rho\erfc(\frac{|x|}{\sqrt{2t}})\dif x$.
\end{proposition}
\begin{proof}
    For any $I\subseteq \R$, let $A(I)=\{(s,e)\in\R^+\times \mathcal{E}: e(t-s)\in I,s+\xi(e)>t\}$. Then $X_t(I)=\Theta(A(I))$ is a Poisson random variable. And if $I_1\cap I_2=\emptyset$, we have $A(I_1)\cap A(I_2)=\emptyset$. Therefore, $X_t(I_1)=\Theta(A(I_1))$ is independent of $X_t(I_2)=\Theta(A(I_2))$ as $\Theta$ is a Poisson point process. It remain to determine the intensity measure of $X_t$. Note that
    \begin{equation*}
        \begin{split}
            \rho \dif t\otimes \n(A(I))&=\rho \int_0^t \n(\{e\in \mathcal{E}:e(t-s)\in I,\xi(e)>t-s\})\dif s\\
            &=\rho \int_0^t\int_I\n_{t-s}(\dif x)\dif s\\
            &=\rho \int_I\int_0^t\frac{|x|\exp({-x^2/2s})}{\sqrt{2\pi s^3}}\dif s\dif x\\
            &=\rho \int_I\erfc\left(\frac{|x|}{\sqrt{2t}}\right)\dif x,
        \end{split}
    \end{equation*}
    where in the second equation we use definition of entrance law of Brownian motion given by (\ref{entrance law}).
    Hence we get the desired intensity measure for $X_t$.
\end{proof}
\begin{remark}\label{distribution of X_t+}
    Let $X_t^+$ denote the restriction of $X_t$ to $[0,\infty)$, which is the point process of particles whose position is above zero. Then by the similar computation we get $X_t^+$ is a Poisson point process with intensity measure $\rho \erfc({x}/{\sqrt{2t}})\dif x$ for $x\geq 0$. 
\end{remark}

It is easy to check that $\int_\R \Lambda(\dif x)<\infty$. Therefore, according to Proposition \ref{Poi X_t}, $\E[N_t]$ is finite for any $t\in\R^+$. So at any deterministic time, there are finitely many alive particles.

\section{Exceptional times and points}\label{proof of theorem 1 and 2}
\subsection{Lower bound for the dimension of exceptional times}

Throughout this subsection, $f$ is an even continuous function that is non decreasing on $[0,\infty)$ satisfies condition (\ref{requirement of test function}) in Theorem \ref{loglog theorem}. Recall that 
\begin{equation*}
    \T=\{\tau\in\R^+:N_\tau=\infty\},\quad \T(f)=\{\tau\in \T:\langle X_\tau,f\rangle=\infty\}
\end{equation*}
are sets defined in Theorem \ref{exceptional times} and \ref{loglog theorem}. The main goal of this subsection is to prove the following lemma: 
\begin{proposition}\label{main lemma for lower bound}
    For every closed interval $[a,b]\subseteq\R^+$, almost surely, $\dim(\T(f)\cap[a,b])\geq1/2$.
\end{proposition}
This lemma is crucial because it implies both the existence of dense sets of exceptional times and the lower bound for the Hausdorff dimensions of $\T$ and $\T(f)$.
\begin{corollary}\label{exists of T, T(f), and lowerbound}
    Almost surely, $\T$ and $\T(f)$ are dense subsets of $\R^+$ with $\dim(\T)\geq\dim(\T(f))\geq1/2$.
\end{corollary}
\begin{proof}
Since $\T(f)\subseteq \T$, it is enough to prove the corresponding statements for $\T(f)$. By Lemma \ref{main lemma for lower bound}, almost surely, $\T(f)\cap [a,b]$ is non empty for every pair of rational numbers $a,b$ with $0\leq a< b$. Hence $\T(f)$, and therefore also $\T$ is dense in $\R^+$. The dimension inequalities follow directly from Lemma \ref{main lemma for lower bound} and the inclusion $\T(f)\subseteq \T$.
\end{proof}

The proof of Lemma \ref{main lemma for lower bound} is based on the construction of a Cantor-like random set contained in $\T(f)$. We shall then apply the following lemma by Orey and Taylor \cite[Lemma 2.2]{Cantor-like_set}, which they used to prove the lower bound for the Hausdorff dimension of the set of fast times of Brownian motion.
\begin{lemma}\label{Cantor dim}
    Suppose that $K$ is a Cantor-like set with a representation
    \begin{equation*}    K=\bigcap_{n=1}^\infty\bigcup_{i=1}^{M_n}J_{n,i},
    \end{equation*}
    where the $J_{n,i}\ (1\leq i\leq M_n)$ are disjoint closed subintervals of $[0,1]$. Let $c,\delta>0$ be two positive numbers. Then
    \begin{equation*}
        \dim K\geq s
    \end{equation*}
    if, for every interval $J\subseteq[0,1]$ with $\diam(J)<\delta$, there is a finite integer $n(J)$ such that
    \begin{equation*}
        M_n(J)\leq c\left(\diam(J)\right)^sM_n\quad\text{for}\quad n\geq n(J),
    \end{equation*}
    where $M_n(J)$ denotes the number of intervals $J_{n,i}$ contained in $J$.
\end{lemma}

We now introduce some notations needed for the construction. For any $h>0$, define 
\begin{equation*}
    G_h:=\left\{e\in\mathcal{E}:e(u)\geq h,\ \text{for all }u\in[h^2,2h^2],\ \max_{0\leq s\leq\xi(e)}|e(u)|\leq 2h\right\}.
\end{equation*}
Let $c_0=\n(G_1)$ then by Proposition \ref{Brownian excursion scaling} we get
\begin{equation}\label{scale n(G_h)}
    \n(G_h)=\n\left(\left\{h^{-1}e(h^2u)\geq1,\ \text{for all }u\in[1,2],\ \max_{0\leq u\leq h^{-2}\xi(e)}|e(u)|\leq 2\right\}\right)=h^{-1}c_0.
\end{equation}
Let 
\begin{equation}\label{def of p_0}
    p_0=\pr[\Poi(c_0)\geq 1],
\end{equation}
where $\Poi(c_0)$ is a Poisson random variable with parameter $c_0$. For $n\in\N$ and $q>4$ we set 
\begin{equation*}
    h_n=\exp(-q^n).
\end{equation*}
 The assumption $q>4$ is only technical, at the end of the argument we shall let $q$ tend to infinity. Besides, for $n\geq 2$ let 
 \begin{equation*}
     R_n=\left\lceil \frac{h_{n-1}^2}{2h_n}\right\rceil,\quad m_n=\left\lceil\frac{R_n p_0}{2}\right\rceil\quad\text{and}\quad M_n=\prod_{k=1}^n m_n, 
 \end{equation*}
 where $\lceil\cdot\rceil$ is the ceiling function and $m_1=R_1=1$. We define a sequence of Poisson point process on the time axis $\R^+$ by
\begin{equation*}
    \Theta^{(n)}(I):=\Theta(I\times G_{h_n}),\quad\text{for any}\quad I\subseteq \R^+.
\end{equation*}
According to (\ref{scale n(G_h)}), it is easy to check that the intensity measure of $\Theta^{(n)}$ is
 \begin{equation*}
     \n(G_{h_n})\dif t=h_n^{-1}c_0\dif t.
 \end{equation*}
  This sequence of Poisson point process can be interpreted as follows: for any $s\in\Theta^{(n)}$, there exists a particle born at $s$ and its corresponding path remains away from zero during a time interval of length at least $h_n^2$. Since the sets $G_{h_n}$ are disjoint for different $n$, the Poisson point processes $\left(\Theta^{(n)}\right)_{n\geq 1}$ are independent. Finally, we let
\begin{equation*}
    \mathcal{F}_n=\sigma\left(\{\Theta^{(1)},\cdots,\Theta^{(n)}\}\right)
\end{equation*}
be the sigma algebra generated by $\Theta^{(1)},\cdots,\Theta^{(n)}$.

\begin{figure}[htpb]
\centering
\includegraphics[scale=0.45]{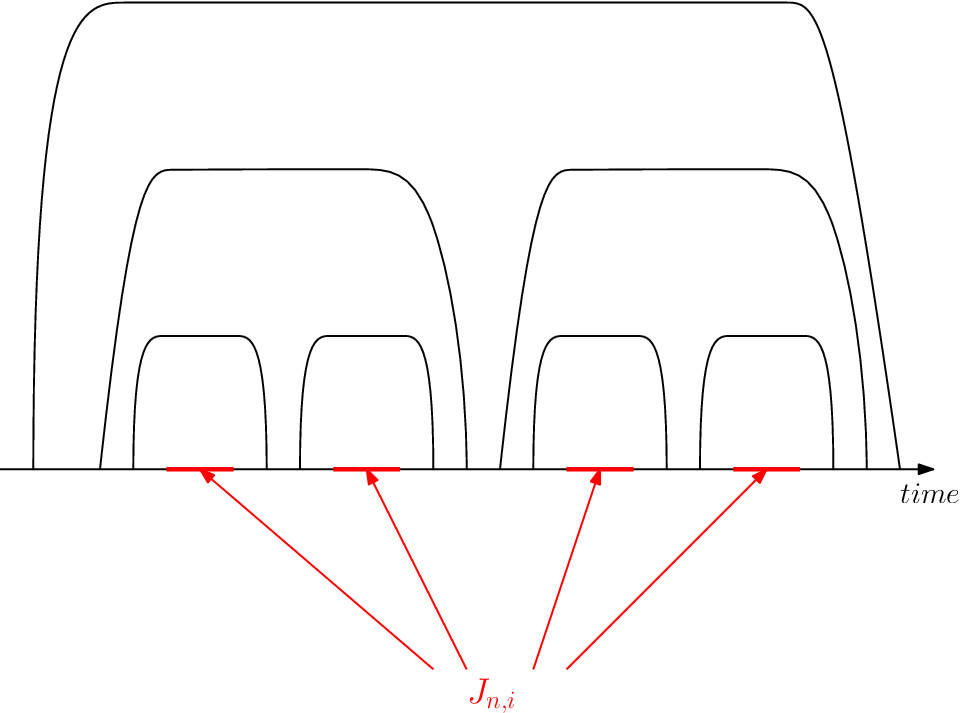}
\caption{The Cantor-like set on the event $E_\infty$ with $J_{n,i}$ marked with red color and path of partilces marked as black curve.}
\label{Cantor-like}
\end{figure}

\begin{lemma}\label{construction of E_infty}
    There exists a decreasing sequence of events $E_1\supseteq E_2\supseteq\cdots$ and a collection of closed random intervals $\left((J_{k,i})_{i=1}^{M_k}\right)_{k=1}^\infty$ such that the following properties hold for every $n\in\N$.
    \begin{enumerate}[label=(\roman*)]
        \item $E_n\in\mathcal{F}_n$.
        \item For any $1\leq i\leq M_k$ and $k\leq n$, almost surely $|J_{k,i}|=h_k^2$ and $J_{k,i}$ is measurable with respect to $\mathcal{F}_k$.
        \item On the event $E_n$, for any $k\leq n$ and $t\in J_{k,i}$, there exists $(t_j,e_j)\in\Theta$ such that $e_j(t-t_j)\in[h_k,2h_k]$.
        \item On the event $E_n$, for every $k\leq n$ the intervals are sub-intervals of $[0,1]$ and $\gap(J_{k,i},J_{k,j})\geq h_k-2h_{k}^2$ whenever $i\neq j$, where $\gap([a,b],[c,d])=(c-b)\vee (a-d)$ is the gap of two closed intervals. In particular on the event $E_n$, the intervals $(J_{k,i})_{i=1}^{M_k}$ are disjoint for all $k\leq n$.
    \end{enumerate}
    Moreover, for any $n\in\N$ we have
    \begin{equation*}
    \pr[E_n]=
        \begin{cases}
            \pr\left[\Poi\left(c_0\frac{1-2h^2_1}{h_1)}\right)\geq1\right],\quad&\text{for }n=1,\\ \\
            \pr\left[\Poi\left(c_0\frac{1-2h^2_1}{h_1)}\right)\geq1\right]\prod_{k=2}^n\pr\left[\Bin(R_k,p_0)\geq m_k\right]^{M_{k-1}},\quad&\text{for }n\geq 1,
        \end{cases}
    \end{equation*}
    where $\Bin(R_k,p_0)$ is a Binomial distributed random variables with parameter $(R_k,p_0)$.
\end{lemma}
\begin{remark}
    This Lemma provides a Cantor-like set with representation
    \begin{equation*}
        K^{q}:=\bigcap_{n=1}^\infty\bigcup_{j=1}^{M_n}J_{n,i}.
    \end{equation*}
    We write $K^q$ as each interval $J_{n,i}$ depends on the sequence $h_n=\exp(-q^n)$, so $K^q$ is a random set depends on $q$. The strategy for proving Lemma \ref{main lemma for lower bound} is as follows: first we will show that on the event 
    \begin{equation}\label{def of E_infty}
        E_{\infty}:=\bigcap_{n=1}^\infty E_n,
    \end{equation} we have $K^q\subseteq \T(f)$. (See figure \ref{Cantor-like} as an illustration.) Second by applying Lemma \ref{Cantor dim}, we shall prove $\dim(K^q)\geq s(q)$ on the event $E_\infty$, where $s(q)$ is a function of $q$ such that $s(q)\to1/2$ from below as $q\to\infty$. Then show that $\pr[E_\infty]>0$, which implies that $\dim(T(f))\geq s(q)$ with positive probability. Finally we shall use Kolmogorov's 0-1 law to make the probability equal to 1 and take $q$ to infinity.
\end{remark}

\begin{proof}
    We prove this Lemma by induction. First we construct $E_1$ and $J_{1,1}$. Let $\sigma_1=\inf\{s\in[0,1-2h_1^2]:s\in\Theta^{(1)}\}$ and $E_1=\{\sigma_1<\infty\}$. Define $\tau_1$ to be
    \begin{equation*}
        \tau_1:=
        \begin{cases}
            \sigma_1+h_1^2\quad&\text{on}\quad E_1\\
            1\quad&\text{on}\quad E_1^c
        \end{cases}
    \end{equation*}
    and $J_{1,1}=[\tau_1,\tau_1+h^2_1]$. Clearly $|J_{1,1}|=h_1^2$. Since both $E_1$ and $J_{1,1}$ are determined by $\Theta^{(1)}$ properties $(i)$ and $(ii)$ hold. Property $(iii)$ follows from the definition of $\Theta^{(1)}$ and $(iv)$ is trivial as there is only one interval $J_{1,1}$. Since $\Theta^{(1)}$ is a Poisson point process with intensity $c_0h_1^{-1}\dif t$ we have
    \begin{equation*}
        \pr\left[E_1\right]=\pr\left[\Theta^{(1)}([0,1-2h_1^2])\geq 1\right]=\pr\left[\Poi\left(c_0\frac{1-2h^2_1}{h_1}\right)\geq 1\right].
    \end{equation*}
    Suppose that the events $E_1,\cdots,E_{n-1}$ and closed intervals $\left((J_{k,i})_{i=1}^{M_k}\right)_{k=1}^{n-1}$ have been constructed and satisfy the required properties up to $n-1$. We shall construct $E_n$ and intervals $(J_{n,i})_{i=1}^{M_n}$. On the event $E_{n-1}$, for each $i=1,\cdots,M_{n-1}$, we write $J_{n-1,i}=[\tau_{n-1,i},\tau_{n-1,i}+h_{n-1}^2]$ and cut it into close intervals of the form 
    \begin{equation*}
        W_j^{n,i}=[\tau_{n-1,i}+2(j-1) h_n,\tau_{n-1,i}+(2(j-1)+1)h_n].
    \end{equation*}
     for $1\leq j\leq R_n$. Each interval $W_j^{n,i}$ has length $h_n$, and for each fixed $i$ the gap between any two consecutive such intervals is at least $h_n$. (See figure \ref{W_j^{n,i}} as an illustration.) By the induction hypothesis we have $\gap(J_{n-1,i},J_{n-1,i+1})\geq h_{n-1}-2h^2_{n-1}$. Hence 
     \begin{equation*}
         \gap(W_{R_n}^{n,i},W_{1}^{n,i+1})\geq h_{n-1}-2h^2_{n-1}-h_n\geq h_n,
     \end{equation*}
     where the last inequality is followed by $q>4$. (See also figure \ref{W_j^{n,i}} as an illustration. In fact here is why we require $q>4$.) Therefore, any two intervals in the collection $\left((W^{n,i}_j)_{j=1}^{R_n}\right)_{i=1}^{M_{n-1}}$ is separated by gaps of at least $h_n$. For complement we define $W_j^{n,i}=[i+jh_n,i+(j+1)h_n)$ for all $j=1,\cdots,R_n$ and $i=1,\cdots,M_{n-1}$ on the event $E^c_{n-1}$ so that they are disjoint intervals with $|W_j^{n,i}|=h_n$, but the proof only need $W_j^{n,i}$ to be well defined on the event $E_n$. For each $i$ and $j$ set
    \begin{equation*}
        Y_j^{n,i}=\ind_{\{\Theta^{(n)}(W^{n,i}_j)\geq1\}}\quad\text{and}\quad X^{n,i}=\sum_{j=1}^{R_n}Y_j^{n,i}.
    \end{equation*}
    We construct $E_n$ by
    \begin{equation*}
        E_n=E_{n-1}\cap\left\{\bigcap_{i=1}^{M_{n-1}}\{X^{n,i}\geq m_n\}\right\}.
    \end{equation*}
    Note that the collection of intervals $\left((W_j^{n,i})_{j=1}^{R_n}\right)_{i=1}^{M_{n-1}}$ is determined by $\left((J_{k,i})_{i=1}^{M_k}\right)_{k=1}^{n-1}$ which is measurable with respect to $\mathcal{F}_{n-1}$ according to the induction hypothesis. And $\Theta^{(n)}$ is measurable with respect to $\mathcal{F}_n$. Hence we get $Y_{j}^{n,i}$ is measurable with respect to $\mathcal{F}_n$. Therefore $E_n\in\mathcal{F}_n$. This proves property $(i)$. 
    
    Now we need to construct $(J_{n,i})_{i=1}^{M_n}$. On the event $E_n$ for each $i=1,\cdots M_{n-1}$ there are at least $m_n$ intervals of the form $W^{n,i}_j$ such that $\Theta^{(n)}(W^{n,i}_j)\geq 1$. We denote those intervals by $W_{j_1}^{n,i},W^{n,i}_{j_2},\cdots W_{j_{m_n}}^{n,i}$. For every such intervals $W_{j_k}^{n,i}$, there exists $s\in W_{j_k}^{n,i}$ and $e\in G_{h_n}$ such that $(s,e)\in \Theta$. Let $J_j^{n,i}=[s+h_n^2,s+2h_n^2]$. On the event $E_n^c$, let $J_j^{n,i}=[1,1+h_n^2]$ for all $i,j$. Relabeling the intervals gives $J_{n,1},\cdots,J_{n,M_n}$. By construction, each $J_{n,j}$ depends on $\left((Y_j^{n,i})_{j=1}^{R_n}\right)_{i=1}^{M_{n-1}}$ which is a collection of $\mathcal{F}_{n}$ measurable random variables. Hence we get $J_{n,j}$ is measurable with respect to $\mathcal{F}_n$ with $|J_{n,j}|=h_n^2$, which proves property $(ii)$. Property $(iii)$ follows directly from the definition of $J_{n,j}$. Indeed if $t\in J_{j}^{n,i}$ then $t=t_i+u$ for some $(t_i,e)\in\Theta$, $e\in G_{h_n}$ and $u\in[h_n^2,2h_n^2]$. Since $e\in G_{h_n}$ we get $e(t-t_i)=e(u)\in [h_n,2h_n]$.
  
    It remains to proof property $(iv)$. Note that on the event $E_n$ 
    the intervals $W_j^{n,i}$ are separated by gaps of at least $h_n$. Since each interval $J^{n,i}_j$ is of the form $[s+h_n^2,s+2h_n^2]$ with $s\in W_j^{n,i}$, passing from $W_j^{n,i}$ to $J_j^{n,i}$ can reduce the gap by at most $2h_n^2$. Therefore, any two distinct intervals of the form $J_{n,i}$ are separated by at least $h_n-2h_n^2$. This proves property $(iv)$.

    Finally we compute the probability of $E_n$. By the tower property we have
    \begin{equation}\label{probabilit of E_n}
        \pr[E_n]=\E\left[\E\left[\ind_{\{E_{n-1}\}}\ind_{\{\cap_{i=1}^{M_{n-1}}\{X^{n,i}\geq m_n\}\}}\middle|\mathcal{F}_{n-1}\right]\right].
    \end{equation}
    For each $i$,
    \begin{equation*}
        \E\left[\ind_{\{X^{n,i}\geq m_n\}}\middle|\mathcal{F}_{n-1}\right]=\pr\left[\sum_{j=1}^{R_n}Y_j^{n,i}\geq m_n\middle|\mathcal{F}_{n-1}\right].
    \end{equation*}
    Note that $W_j^{n,i}$ are disjoint sets with $|W_j^{{n,i}}|=h_n$ for all $i,j$. Hence condition on $\mathcal{F}_{n-1}$, $\left(\Theta^{(n)}(W_j^{n,i})\right)_{j=1}^{m_n}$ is a collection of independent Poisson random variables with parameter 
    \begin{equation*}
        |W_j^{n,i}|\n(G_{h_n})=c_0h_nh_n^{-1}=c_0,
    \end{equation*}
    where the equality follows by (\ref{scale n(G_h)}). Hence condition on $\mathcal{F}_{n-1}$, $X^{n,i}=\sum_{j=1}^{R_n}Y_j^{n,i}$ is a binomial distributed random variable with parameter $(R_n,p_0)$, where $p_0$ is given by (\ref{def of p_0}). Thus
    \begin{equation*}
        \pr\left[\sum_{j=1}^{R_n}Y_j^{n,i}\geq m_n\middle|\mathcal{F}_{n-1}\right]=\pr\left[\Bin(R_n,p_0)\geq m_n\right].
    \end{equation*}
    Again, the intervals $W^{n,i}_j$ are disjoint for different $i=1,\cdots, M_n$. We get $X^{n,i}$ are independent for different $i$ after conditioning on $\mathcal{F}_{n-1}$. Hence 
    \begin{equation*}
        \E\left[\ind_{\{\cap_{i=1}^{M_{n-1}}{\{X^{n,i}\geq m_n}\}}\middle|\mathcal{F}_{n-1}\right]=\pr\left[\Bin(R_n,p_0)\geq m_n\right]^{M_{n-1}}.
    \end{equation*}
    By (\ref{probabilit of E_n}) and the induction hypothesis, we get
    \begin{equation*}
        \pr\left[E_n\right]=\pr\left[\Poi(c_0\frac{1-2h^2_1}{h_1})\geq1\right]\prod_{k=2}^n\pr\left[\Bin(R_k,p_0)\geq m_k\right]^{M_{k-1}}
    \end{equation*}
    as desired.
\end{proof}
\begin{figure}[htpb]
\centering
\includegraphics[scale=0.45]{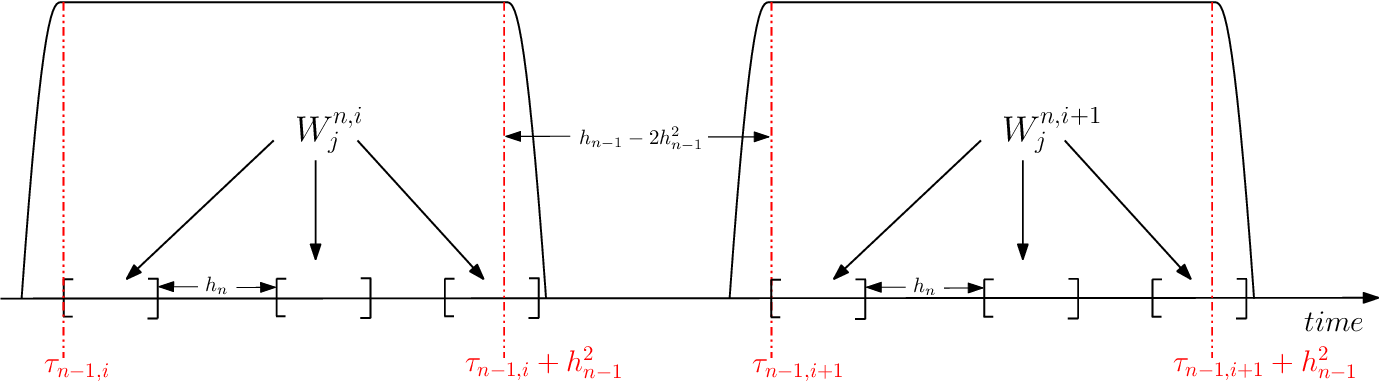}
\caption{The strategy for choosing $W_j^{n,i}$, where the black curve is a path of particle and the end points of $J_{n-1,i},J_{n-1,i+1}$ are marked with red color.}
\label{W_j^{n,i}}
\end{figure}

\begin{lemma}\label{step1}
    On the event $E_\infty$, we have $K^q\subseteq T(f)$ for every $f$ satisfying (\ref{requirement of test function}).
\end{lemma}
\begin{proof}
    According to $(iii)$ of Lemma \ref{construction of E_infty}, on the event $E_\infty$ for any $n\in\N$ and $t\in J_{n,i}$ with $1\leq i\leq M_n$, there exists an alive particle at time $t$ whose position $x_n(t)\in [h_n,2h_n]$. Since $h_n=\exp(-q^n)$ we have $2h_{n+1}<h_{n}$ for all $n\in\N$. Hence for any $t\in K^q$, there exists a sequence of distinct particles alive at time $t$ with positions $x_n(t)\geq h_n$ for all $n\in\N$. For any $f$ satisfying (\ref{requirement of test function}), using the change of variables we obtain
    \begin{align*}
        \int_0^\epsilon \frac{f(x)}{x\log(1/x)}\dif x=\int_{\log(\epsilon^{-1})}^\infty \frac{f(e^{-y})}{y}\dif y=\infty.
    \end{align*}
    Choose $N$ sufficiently large such that $q^n\geq \log(\epsilon^{-1})$ for all $n\geq N$ we get
    \begin{align*}
        \infty=\int_{q^N}^\infty \frac{f(e^{-y})}{y}\dif y=\sum_{n=N}^\infty \int_{q^n}^{q^{n+1}}\frac{f(e^{-y})}{y}\dif y\leq\sum_{n=N}^\infty f(h_n)\log q,
    \end{align*}
    where the inequality follows from the fact that $f$ is non-deceasing on $[0,\infty)$. Hence we get
    \begin{equation*}
        \sum_{n=1}^\infty f(h_n)=\infty.
    \end{equation*} 
    Now fix $t\in K^q$. As observed above, there are distinct particles alive at time $t$ with positions $x_n(t)\geq h_n$, hence
    \begin{equation*}
        \langle X_t,f\rangle\geq \sum_{n=1}^\infty f(x_n(t))\geq\sum_{n=1}^\infty f(h_n)=\infty,
    \end{equation*}
    which implies that $t\in \T(f)$. Hence $K^q\subseteq \T(f)$ on the event $E_\infty$.
\end{proof}

\begin{lemma}
    Let 
     \begin{equation*}
        s(q)=\frac{1}{2}\left(\frac{q-2}{q-1}-\frac{1}{q}\right).
    \end{equation*}
    On the event $E_\infty$, $\dim (K^q)\geq s(q)$ for all $q>4$.
\end{lemma}
\begin{proof}
    In order to apply Lemma \ref{Cantor dim}, we need to show that for any interval $J$ with $|J|<\delta$, there exists $n(J)$ such that 
    \begin{equation}\label{condition for lowerbound}
        M_n(J)\leq 3 |J|^s M_n,\quad\text{for }n\geq n(J)\text{ and }s<s(q).
    \end{equation}
    To do this, we firstly need an estimation of $M_n$. Recall that 
    \begin{equation*}
        M_n=\prod_{k=1}^n m_k,\quad m_k=\left\lceil \frac{R_kp_0}{2}\right\rceil,\quad R_k=\left\lceil\frac{h^2_{k-1}}{2h_k}\right\rceil\approx h_{k-1}^{2-q}.
    \end{equation*}
    We get for $n\geq 2$:
    \begin{align*}
        \log M_n&=\sum_{k=1}^n\log m_k=\sum_{k=2}^n \log \left\lceil\frac{R_kp_0}{2}\right\rceil\approx\sum_{k=2}^n\log \frac{p_0}{2}+(2-q)\log h_{k-1}\\
        &=(n-1)\log{\frac{p_0}{2}}+(q-2)\sum_{k=2}^n q^{k-1}=-C(n-1)+A(q)(q^n-q),
    \end{align*}
    where
    \begin{equation*}
        A(q)=\frac{q-2}{q-1}.
    \end{equation*}
    Hence 
    \begin{equation}\label{estimation of M_n}
        M_n\approx\exp\left(-Cn+A(q)q^n\right).
    \end{equation}
    Now we return to check (\ref{condition for lowerbound}). Let $h_k$ be such that 
    \begin{equation*}
        h_k^2\leq |J|<h_{k-1}^2.
    \end{equation*}
    By Lemma \ref{construction of E_infty} $(iv)$, on the event $E_\infty$ we have $\gap(J_{k,i},J_{k,j})\geq h_k-2h_k^2$ for $i\neq j$. Hence $J$ contains at most $|J|/(h_k-2h_k^2)+1$ intervals of the form $J_{k,i}$. According to our construction, each $J_{k,i}$ contains $M_n/M_k$ intervals of the form $J_{n,j}$ ($j=1,\cdots,M_n$). We can choose $k$ large enough (or choose $\delta$ small enough equivalently) so that $h_k-2h_k^2\geq h_k/2$. Hence 
    \begin{equation}\label{upperbound of M_n(J)}
        M_n(J)\leq \left(\frac{|J|}{h_k/2}+1\right)\frac{M_n}{M_k}=\left(\frac{2|J|}{h_k}+1\right)\frac{M_n}{M_k}
    \end{equation}
    We consider two cases. Case 2 is the crucial one, as it is precisely this case that limits the attainable lower bound on the Hausdorff dimension to $s(q)$.
    
    \textbf{Case 1}: $|J|\leq h_k$. Then by (\ref{upperbound of M_n(J)}) we get $M_n(J)\leq 3 M_n M^{-1}_{k}$. For any $q>4$, by (\ref{estimation of M_n}) we have 
    \begin{equation*}
        h_k^{A(q)-q^{-1}}M_k\approx\exp\left(q^{k-1}-Ck\right)\to\infty\quad \text{as }k\to\infty.
    \end{equation*}
    Hence $M_k^{-1}\leq h_k^{A(q)-q^{-1}}$ for $k$ being large enough. Thus
    \begin{equation*}
        M_{k}^{-1}\leq (h_k^2)^{(A(q)-q^{-1})/2}\leq |J|^{s(q)}.
    \end{equation*}
    Therefore, 
    \begin{equation*}
        M_n(J)\leq3 M_n|J|^{s(q)},\quad\text{for }n\geq k.
    \end{equation*}
    By Lemma \ref{Cantor dim}, we know that in this case
    \begin{equation*}
        \dim(K^q)\geq s(q)\quad\text{for any }q>4.
    \end{equation*}
    \textbf{Case 2}: $|J|> h_k$. Note that
    \begin{equation*}
        h_k^{A(q)-q^{-2}}M_k\approx\exp\left(q^{k-2}-Ck\right)\to\infty\quad \text{as }k\to\infty.
    \end{equation*}
    Hence $M_k^{-1}\leq h_k^{A(q)-q^{-2}}$ for large enough $k$. Therefore in this case we have
    \begin{equation}\label{upperbound of M_n(J) in case 2}
        M_n(J)\leq 3h_k^{-1}|J|M_k^{-1}M_n\leq 3M_n|J|h^{A(q)-1-q^{-2}}_k.
    \end{equation}
     For any $s\leq 1-(q/2)\left(1-A(q)+q^{-2}\right)=s(q)$, it is easy to check
    \begin{equation*}
        h_{k-1}^{2(1-s)}\leq h_{k-1}^{q\left(1-A(q)+q^{-2}\right)}=h_k^{1-A(q)+q^{-2}}.
    \end{equation*}
    Since $|J|\leq h_{k-1}^2$, we get $|J|^{1-s}\leq h_{k-1}^{2(1-s)}\leq h_k^{1-A(q)+q^{-2}}$, which is equivalent to 
    \begin{equation*}
        |J|\leq|J|^sh_k^{1-A(q)+q^{-2}}.
    \end{equation*}
    Hence by (\ref{upperbound of M_n(J) in case 2}) we get $M_n(J)\leq 3 M_n|J|^s$. Applying Lemma \ref{Cantor dim}, we have $\dim(K^q)\geq s$ for any $s\leq s(q)$. Consequently we get $\dim (K^q)\geq s(q)$ 
\end{proof}

\begin{lemma}\label{step3}
    Let $E_\infty$ be defined as in (\ref{def of E_infty}). Then we have $\pr[E_\infty]>0$.
\end{lemma}
\begin{proof}
    Note that $E_\infty=\cap_{n=1}^\infty E_n$ with $(E_n)_{n\geq1}$ being an decreasing sequence of events. Then we have $\pr[E_\infty]=\lim_{n\to\infty}\pr[E_n]$. By Lemma \ref{construction of E_infty}, we have 
    \begin{equation*}
        \pr\left[E_n\right]=\pr\left[\Poi\left(c_0\frac{1-2h^2_1}{h_1}\right)\geq1\right]\prod_{k=2}^n\pr\left[\Bin(R_k,p_0)\geq m_k\right]^{M_{k-1}}.
    \end{equation*}
    Hence it is enough to show 
    \begin{equation*}
        -\sum_{n=2}^\infty M_{n-1}\log\left(\pr[\Bin(R_n,p_0)\geq m_n]\right)<\infty.
    \end{equation*}
    We will show that this series is easily summable.
    Let $X\sim\Bin(R_n,p_0)$ then $\E[X]=R_np_0\leq 2m_n$. By using Markov inequality and the moment generating function of $X$, for any $a>0$ we get the following Chernoff bound,
    \begin{align*}
        \pr\left[X\leq \frac{\E[X]}{2}\right]&\leq \exp\left({\frac{a\E[X]}{2}}\right)\left(1-p_0+p_0\exp(-a)\right)^{R_n}\\
        &\leq\exp\left(\frac{a\E[X]}{2}+p_0(\exp(-a)-1)R_n\right)\\
        &=\exp\left(\E[X](\frac{a}{2}+\exp(-a)-1)\right).
    \end{align*}
    By taking $a=\log(2)$ we get for some $C>0$,
    \begin{equation*}
        \pr\left[X\leq\frac{\E[X]}{2}\right]\leq \pr[\Bin(R_n,p_0)\leq m_n]\leq\exp(-CR_n).
    \end{equation*}
    Hence from the fact that $\log(1-x)\geq-2x$ for $x$ being small enough, we get
    \begin{equation}\label{take log}
        \sum_{n=2}^\infty M_{n-1}\log\left(\pr[\Bin(R_n,p_0)\geq m_n]\right)\geq-2\sum_{n=2}^\infty M_{n-1}\exp\left(-C R_n\right).
    \end{equation}
    Note that by (\ref{estimation of M_n}), we have
    \begin{equation*}
        M_n\approx\exp\left(-Cn+A(q)q^n\right)\leq\exp\left(A(q)q^n\right).
    \end{equation*}
    Thus for large enough $n$ we get
    \begin{equation*}
         M_{n-1}\exp\left(-CR_n\right)\leq\exp\left(-CR_n+A(q)q^n\right)\approx\exp\left(-C\exp\left({(q-2)q^{n-1}}\right)+A(q)q^n\right),
    \end{equation*}
    which is summable. Hence by (\ref{take log}) we get
    \begin{equation*}
        \sum_{n=2}^\infty M_{n-1}\log\left(\pr[\Bin(R_n,p_0)\geq m_n]\right)>-\infty,
    \end{equation*}
    which is what we desired to prove.
\end{proof}
Now we are ready to finish this section by proving Lemma \ref{main lemma for lower bound}. 
\begin{proof}[proof of Lemma \ref{main lemma for lower bound}]
    Without loss of generality, we assume $a=0$ and $b=1$ because the proof can be easily generalized to any closed interval $[a,b]\subseteq\R^+$.
    From Lemma \ref{step1} to \ref{step3} we know that for any $q>4$, we have $\dim\left(\T(f)\cap[0,1]\right)\geq s(q)$ with positive probability. It is suffices to show that this is a tail event. For any $k\in\N$, let 
    \begin{equation}\label{def of tilde G_k}
    \tilde G_k=
    \begin{cases}
                \{e\in\mathcal{E}:\max_{0\leq u\leq\xi(e)}|e(u)|\in[h_{k},h_{k-1}]\}&\quad\text{for }k\geq2,\\
              \{e\in\mathcal{E}:\max_{0\leq u\leq\xi(e)}|e(u)|\geq h_1\}&\quad\text{for }k=1.
    \end{cases}
    \end{equation}
    and define $\Theta^{(k)}$ to be the Poisson point process on the $\R^+\times\mathcal{E}$ such that
    \begin{equation*}
        \tilde\Theta^{(k)}(W)=\Theta\left(W\cap([0,1]\times \tilde G_k)\right)
    \end{equation*}
    for any $W\subseteq \R^+\times\mathcal{E}$. Let $\mathcal{F}^k$ be the sigma algebra generated by $\tilde \Theta^{(k)}$, then $(\mathcal{F}^k)_{k\geq 1}$ is a sequence of independent sigma algebra as $(\tilde{G}_k)_{k\geq 1}$ is pairwise disjoint sequence. Note that 
    \begin{equation}\label{upper bound of n(tilde G_k)}
        \n(\tilde{G}_k)\leq \n\left(\{e\in\mathcal{E}:\max_{0\leq u\leq \xi(e)}|e(u)|\geq h_k\}\right)=\rho h^{-1}_k
    \end{equation}
    (see \cite[Proposition 2]{Rogers_excursion} for proof of above equality.) We get
    \begin{equation*}
        \E[\tilde\Theta^{(k)}(\R^+\times\mathcal{E})]=\E[\Theta([0,1]\times \tilde G_k)]\leq\rho h^{-1}_k<\infty.
    \end{equation*}
    Hence the number of particles whose maximum height is away from zero at least $h_k$ is almost surely finite for all $t\in[0,1]$ and $k\in\N$. Therefore the set $\T(f)\cap[0,1]$ does not depend on those finitely many particles, which implies that 
    \begin{equation*}
        \left\{\dim\left(\T(f)\cap[0,1]\right)\geq s(q)\right\}\in\bigcap_{n=1}^\infty\bigcup_{k=n}^\infty \mathcal{F}^k
    \end{equation*}
    is a tail event. By Kolmogorov's 0-1 law, we get $\dim\left(\T(f)\cap[0,1]\right)\geq s(q)$ almost surely. Note that $s(q)\to 1/2$ as $q\to \infty$, therefore we have $\dim\left(\T(f)\cap[0,1]\right)\geq 1/2$ almost surely.
\end{proof}
\begin{remark}\label{no cluster point except zero}
    In this proof we have shown that the number of particles born in $[0,1]$ with top positions being away from zero at least $\epsilon$ is almost surely finite. Thus at any exceptional time $\tau\in \T$, the collection of positions of alive particles contains no cluster point other than zero.
\end{remark}

\subsection{Upper bound for the dimension of exceptional times}\label{Brownian fan}
To obtain an upper bound for $\dim \T$, we shall prove that the Hausdorff dimension of the set of exceptional times at which the Brownian fan contains infinitely many particles is almost surely at most $1/2$. Since our singular birth-death model coincides with the first generation of the Brownian fan, this will imply that 
\begin{equation*}
    \dim \T\leq 1/2.
\end{equation*} 
We begin by recalling the definition of the Brownian fan introduced in \cite{Fan}. Define the extended excursion space $\bar{\mathcal E}$ by
\begin{equation*}
  \bar{\mathcal{E}}:=\{w\in C(\R^+):w \text{ is constant outside a unique time interval }
[\eta(w),\eta(w)+\xi(w)]
\}.
\end{equation*}
For any $t\in\R^+$ and $w\in \bar{\mathcal{E}}$ define a map $S_{t,w}:\mathcal{E}\to\bar{\mathcal{E}}$ by
\begin{equation*}
(S_{t,w}e)(s)=
\begin{cases}
    e(s-t)+w(t)\quad&\text{for }t\leq s<t+\xi(e),\\
    w(t)\quad&\text{else}.
\end{cases}
\end{equation*}
Therefore $S_{t,w}$ shifts the starting point of  excursion $e$ from $(0,0)$ to $(t,w(t))$. Fix $T\in\R^+$ and let $\mathbf{0}\in C([0,T]])$ be the zero function. We recursively define a sequence of Poisson point processes. Let $\mu^1_T$ be a Poisson point process on $\bar{\mathcal{E}}$ with intensity measure 
\begin{equation*}
    \rho\int_0^T \n\circ S_{t,\mathbf{0}}^{-1}\dif t.
\end{equation*}
For $n\geq 2$, let $\mu^n_T$ be a (conditional independent of $\mu^{l}_T$ for $l\leq n-1$) Poisson point process with intensity measure 
\begin{equation*}
    \rho\int_{\bar{\mathcal{E}}}\int_{\eta(w)}^{\eta(w)+\xi(w)}\n\circ S_{t,w}^{-1}\dif t \mu^{n-1}_T(\dif w).
\end{equation*}
\begin{definition}
    The Brownian fan is a point process on $\bar{\mathcal{E}}$ given by
    \begin{equation*}
        \mu_T^\infty=\sum_{n=1}^\infty\mu_T^n.
    \end{equation*}
\end{definition}
We claim that the first generation $\mu^1_T$ coincides with our singular birth-death model. Indeed, each atom $(t_i,e_i)\in\Theta$ with birth time $t_i\leq T$ can be view as a path $e_i$ shifted 
so that it starts at $(t_i,0)$. Therefore if we define the map $\Phi:[0,T]\times\mathcal{E}\to\bar{\mathcal{E}}$ given by $\Phi(t,e)=S_{t,\mathbf{0}}e$, then our singular birth-death model can be viewed as the pushforward point process $\Theta\circ\Phi^{-1}$, the corresponding pushforward point process is a Poisson point process on $\bar{\mathcal{E}}$ with intensity measure $(\rho\dif t\otimes \n)\circ\Phi^{-1}$. For any $A \subseteq \bar{\mathcal{E}}$, we have
\begin{equation*}
    \begin{split}
        (\rho\dif t\otimes\n) \circ\Phi^{-1}(A)&=\rho\dif t\otimes \n\left(\{(t,e)\in[0,T]\times \mathcal{E}:\Phi(t,e)\in A\}\right)\\
        &=\rho\int_0^T\int_{\mathcal{E}}\ind_{\{S_{t,\mathbf{0}}e\in A\}}\n(\dif e)\dif t\\
        &=\rho\int_{0}^T\n\left(\{e\in\mathcal{E}:S_{t,\mathbf{0} }e\in A\}\right)\dif t\\
        &=\rho\int_0^T \n\circ S_{t,\mathbf{0}}^{-1}(A)\dif t,
    \end{split}
\end{equation*}
which is precisely the intensity measure of $\mu^1_T$. Define the set of excursions alive at time $t$ by
\begin{equation*}
    \mathcal{N}_t=\{w\in\bar{\mathcal{E}}:\eta(w)<t<\eta(w)+\xi(w)\}.
\end{equation*}
Then the number of particles born before time $T$ and alive at $t$ for the Brownian fan is given by
\begin{equation*}
    N^{\mathrm{Fan}(T)}_t=\mu^\infty_T(\mathcal{N}_t).
\end{equation*}
We recall the following theorem from \cite[Theorem 2.14]{Fan}, which provides an upper bound on the exponential moment of $N^{\mathrm{Fan}(T)}_t$.
\begin{theorem}\label{number of partilces in fan}
    There exist constants $C,K>0$ such that 
    \begin{equation*}
        \E\left[\exp\left(K^{-1}\exp(-Kt)N_t^{\mathrm{Fan}(T)}\right)\right]\leq C,
    \end{equation*}
    holds uniformly over all $t,T>0$.
\end{theorem}
For any $ \delta>0$, let 
\begin{equation*}
    A_{\delta,T}=\{w\in\bar{\mathcal{E}}:\eta(w)\leq T,\ \xi(w)\geq \delta\}
\end{equation*}
be the set of excursions that start before time $T$  and have lifetime at least $\delta$. Then $\mu^\infty_T(A_{\delta,T})$ is the total number of particles in the Brownian fan born before time $T$ and have lifetime at least $\delta$.
\begin{proposition}\label{finite partilce with long life}
    For any $\delta>0$,
    \begin{equation*}
        \E[\mu^\infty_{T}(A_{\delta,T})]<\infty.
    \end{equation*}
    In particular, almost surely the total number of particles born before time $T$ with lifetime at least $\delta$ is finite.
\end{proposition}
\begin{proof}
    Let $F:\bar{\mathcal{E}}\to\R$ be the functional given by
    \begin{equation*}
        F(w)=\ind_{\{\xi(w)\geq \delta,\eta(w)\leq T\}}.
    \end{equation*}
    Then for any $n\geq 2$ conditioning on the previous generation $\mu^{n-1}_{T}=\sum_{i\geq 1}\delta_{w_i}$ we get
    \begin{equation}\label{compute expectation by condition}
        \begin{split}
            \E\left[\mu^n_T(A_\delta)\middle|\mu^{n-1}\right]&=\E\left[\int_{\bar{\mathcal{E}}}F(v)\mu^{n}_T(\dif v)\middle|\mu_T^{n-1}\right]
            \\&=\E\left[\rho\int_{\bar{\mathcal{E}}}F(v)\sum_{i\geq 1}\int_{\eta(w_i)}^{\eta(w_i)+\xi(w_i)}\n\circ S^{-1}_{t,w_i}(\dif v)\dif t\middle|\mu^{n-1}_T\right]\\
            &=\E\left[\rho\sum_{i\geq 1}\int_{\eta(w_i)}^{\eta(w_i)+\xi(w_i)}\int_{{\mathcal{E}}}F(S_{t,w_i}v)\n(\dif v)\dif t\middle|\mu^{n-1}_T\right]\\
            &=\E\left[\rho\sum_{i\geq1}\int_{\eta(w_i)}^{\eta(w_i)+\xi(w_i)}\int_{\mathcal{E}}\ind_{\{\xi(e)\geq \delta\}}\ind_{\{t\leq T\}}\n(\dif e)\dif t\middle|\mu_T^{n-1}\right]\\
            &=\E\left[\rho\n(\{e\in\mathcal{E}:\xi(e)\geq \delta\})\sum_{i\geq1}\int_{\eta(w_i)}^{\eta(w_i)+\xi(w_i)}\ind_{\{t\leq T\}}\dif t\middle| \mu^{n-1}_T\right]\\
            &=\rho\n(\{e\in\mathcal{E}:\xi(e)\geq \delta\})\E\left[\int_0^T\sum_{i\geq 1}\ind_{\{\eta(w_i)<t<\eta(w_i)+\xi(w_i)\}}\dif t\middle|\mu^{n-1}_T\right]\\
            &=\rho\n(\{e\in\mathcal{E}:\xi(e)\geq \delta\})\int_0^T\mu^{n-1}_T(\mathcal{N}_t)\dif t.
        \end{split}
    \end{equation}
    Note that
    \begin{equation*}
        \rho \n(\{e\in\mathcal{E}:\xi(e)\geq\delta\})=\rho \int_\R \n_{\delta}(\dif x)=\rho\sqrt{\frac{2}{\pi\delta}},
    \end{equation*}
    where $\n_t(\dif x)$ is the entrance law of Brownian excursion defined in (\ref{entrance law}). Hence for $n\geq 2$,
    \begin{equation*}
        \E[\mu^n_T(A_{\delta,T})]=\rho \sqrt\frac{2}{\pi\delta}\int_0^T\E[\mu_T^{n-1}(\mathcal{N}_t)]\dif t.
    \end{equation*}
    For the first generation,
    \begin{equation*}
        \E\left[\mu^1_T(A_{\delta,T})\right]=\int_0^T \rho \n(\{e\in\mathcal{E}:\xi(e)\geq\delta\})\dif t=T\rho\sqrt{\frac{2}{\pi\delta}}.
    \end{equation*}
    Consequently,
    \begin{equation}\label{particles with long life}
        \begin{split}
            \E\left[\mu^\infty_T(A_{\delta,T})\right]&=\E\left[\sum_{n=1}^\infty\mu^n_T(A_{\delta,T})\right]\\
            &=\rho\sqrt{\frac{2}{\pi\delta}}\left(T+\int_0^T\E\left[\sum_{n=2}^\infty\mu_T^{n-1}(\mathcal{N}_t)\right]\dif t\right)\\
            &=\rho\sqrt{\frac{2}{\pi\delta}}\left(T+\int_0^T\E[N^{\mathrm{Fan}(T)}_t]\dif t\right).
        \end{split}
    \end{equation}
    By Theorem \ref{number of partilces in fan} and the inequality $1+x\leq\exp(x)$ we get
    \begin{equation*}
        \E\left[K^{-1}\exp(-Kt)N^{\mathrm{Fan}(T)}_t\right]\leq C-1.
    \end{equation*}
    Therefore
    \begin{equation}\label{upper bound of N_t^Fan}
        \E\left[N^{\mathrm{Fan}(T)}_t\right]\leq (C-1)K\exp(Kt).
    \end{equation}
    It follows by (\ref{particles with long life}) that 
    \begin{equation*}
        \E\left[\mu_T^\infty(A_\delta)\right]\leq\rho\sqrt{\frac{2}{\pi\delta}}\left(T+\int_0^T(C-1)K\exp(Kt)\dif t\right)\leq C_{\delta,T}<\infty
    \end{equation*}
    as desired.
\end{proof}
For a fix $T>0$, define the set of exceptional times of the Brownian fan by
\begin{equation*}
    \T^{\mathrm{Fan}(T)}:=\{\tau\in\R^+:N^{\mathrm{Fan}(T)}_\tau=\infty\}.
\end{equation*}
Then $\T\cap[0,T]\subseteq\T^{\mathrm{Fan}}$ since all particles in the first generation of the Brownian fan are born during the time interval $[0,T]$.

\begin{remark}\label{independent of short time} In the proof of Proposition \ref{finite partilce with long life} we showed that the number of particles with life time greater than $\delta$ is finite almost surely. Therefore removing all particles with lifetime greater than $\delta$ leaves both $\T$ and $\T^\mathrm{Fan}$ unchanged.
\end{remark}

\begin{lemma}\label{upper bound}
    For any $T\in\R^+$, almost surely $\dim(\T^{\mathrm{Fan}(T)}\cap [0,T])\leq 1/2$.
\end{lemma}
\begin{proof}
    For any $\delta>0$ let $J(\delta)$ be the collection of intervals given by
    \begin{equation*}
        J(\delta)=\left\{[\eta(w),\eta(w)+\xi(w)]:w\in\mu^\infty_T,\ \eta(w)<T,\ \ \xi(w)<\delta\right\}.
    \end{equation*}
    Then all intervals in $J(\delta)$ has diameter less than $\delta$. We claim that $J(\delta)$ is a cover of $\T^{\mathrm{Fan}(T)}\cap[0,T]$. Indeed, for any $\tau\in\T^{\mathrm{Fan}(T)}\cap[0,T]$, there are infinitely many particles alive at time $\tau$. By Proposition \ref{finite partilce with long life} only finitely many of these particles can have lifetime at least $\delta$. Hence there exists $w\in\mu^\infty_T$ with $\xi(w)<\delta$ such that $\tau\in[\eta(w),\eta(w)+\xi(w)]$. Consequently, $\T^{\mathrm{Fan}(T)}\cap[0,T]$ is covered by intervals in $J(\delta)$. To finish the proof it is enough to show that for any $\alpha>1/2$ and $\delta>0$
    \begin{equation*}
        H^\alpha_\delta:=\sum_{I\in J(\delta)}(\diam(I))^\alpha<\infty,\quad\text{a.s.}
    \end{equation*}
    Let $F:\bar{\mathcal{E}}\to\R$ be given by
    \begin{equation*}
        F(v)=\ind_{\{\xi(v)<\delta,\eta(v)<T\}}(\xi(v))^\alpha.
    \end{equation*}
    By a computation analogous to (\ref{compute expectation by condition}) we get for $n\geq 2$,
    \begin{equation*}
        \E\left[\int_{\bar{\mathcal{E}}}F(v)\mu_T^n(\dif v)\right]=\int_\mathcal{E}\rho(\xi(e))^\alpha \ind_{\{\xi(e)<\delta\}}\n(\dif e)\E\left[\int_0^T\mu_T^{n-1}(\mathcal N_t)\dif t\right].
    \end{equation*}
    And for $n=1$,
    \begin{equation*}
        \E\left[\int_{\bar{\mathcal{E}}}F(v)\mu_T^n(\dif v)\right]=\int_\mathcal{E}T\rho(\xi(e))^\alpha \ind_{\{\xi(e)<\delta\}}\n(\dif e).
    \end{equation*}
    For any $\delta>0$,
    \begin{equation*}
        \rho\int_\mathcal{E}(\xi(e))^\alpha\ind_{\{\xi(e)<\delta\}}\n(\dif e)=\rho \int_0^\delta (\xi(e))^\alpha\n(\{\xi(e)\in\dif t\}).
    \end{equation*}
    Since $\n(\xi(e)\in\dif t)$ is proportionally to $t^{-3/2}\dif t$, it follows that
    \begin{equation*}
        \int_0^\delta t^\alpha \n(\xi(e)\in\dif t)=\int_0^\delta Ct^{\alpha-3/2}\dif t.
    \end{equation*}
    For any $\alpha>1/2$ this integral is finite and equals $C_\alpha \delta^{\alpha-1/2}$. Therefore, 
    \begin{equation*}
        \begin{split}
            \E\left[H_\delta^\alpha\right]&=\E\left[\int_{\bar{\mathcal E}}F(v)\mu^\infty_T(\dif v)\right]\\
            &=\rho\int_{\mathcal{E}}(\xi(e))^\alpha\ind_{\{\xi(e)<\delta\}}\left(T+\int_0^T\E[N^\mathrm{Fan}_t]\dif t\right)\\
            &\overset{(\ref{upper bound of N_t^Fan})}{\leq} C_{\alpha,T}\rho\int_\mathcal E (\xi(e))^\alpha\ind_{\{\xi(e)<\delta\}}\n(\dif e)\\
            &\leq C_{\alpha,T}\delta^{\alpha-1/2}<\infty.
        \end{split}
    \end{equation*}
    Hence for any $\alpha>1/2$, $\dim(\T^{\mathrm{Fan}(T)}\cap[0,T])\leq 1/2$ almost surely.
\end{proof}
We can now prove $(ii)$ and $(iii)$ of Theorem \ref{exceptional times}.
\begin{proof}[proof of Theorem \ref{exceptional times}]
For (ii), by Lemma \ref{upper bound}, for any $T>0$ almost surely, 
\begin{equation*}
    \dim(\T^{\mathrm{Fan}}\cap[0,T])\leq\frac{1}{2}.
\end{equation*}
Since $\T\cap[0,T]\subseteq \T^{\mathrm{Fan}(T)}$ we obtain
\begin{equation*}
     \dim(\T\cap[0,T])\leq\frac{1}{2}\quad \text{a.s.}
\end{equation*}
    Hence 
    \begin{equation*}
        \dim(\T)=\dim\big(\cup_{n=1}^\infty (\T\cap[0,n])\big)=\limsup_{n\to\infty}\dim\big(\T\cap [0,n]\big)\leq \frac{1}{2}.
    \end{equation*}
    On the other hand, Corollary \ref{exists of T, T(f), and lowerbound} gives $\dim(\T)\geq 1/2$ almost surely. Therefore,
    \begin{equation*}
        \dim(\T)=\frac{1}{2}\ a.s.
    \end{equation*}
    Part $(iii)$ is a direct consequence of Remark \ref{no cluster point except zero}. 
\end{proof}

\subsection{Hausdorff dimension of the set of exceptional times of the Brownian fan}
As a corollary of Theorem \ref{exceptional times} $(ii)$ we can also determine the Hausdorff dimension of the exceptional times of the Brownian fan.
\begin{corollary}
    Let $\T^{\mathrm{Fan}(T)}$ be the set of exceptional times for the Brownian fan defined as before. Then
    \begin{equation*}
        \dim(\T^{\mathrm{Fan}(T)})=\frac{1}{2}
    \end{equation*}
    almost surely.
\end{corollary}
\begin{proof}
    The lower bound follows by Proposition \ref{main lemma for lower bound} and the inclusion $\T\cap[0,T]\subseteq \T^{\mathrm{Fan}(T)}$:
    \begin{equation*}
        \frac{1}{2}\leq \dim(\T\cap[0,T])\leq\dim\T^{\mathrm{Fan}(T)}.
    \end{equation*}
    For the upper bound, choose $S>T$ and consider a Brownian fan $\mu_S^\infty$. By restricting the initial path $\mathbf0$ from $C[0,S]$ to $C[0,T]$, we obtain a natural coupling of $\mu^\infty_S$ and $\mu_T^\infty$. Therefore,
    \begin{equation*}
        \T^{\mathrm {Fan}(T)}\cap[0,S]\subseteq \T^{\mathrm {Fan}(S)}\cap[0,S].
    \end{equation*}
    By Lemma \ref{upper bound} almost surely:
    \begin{equation*}
        \dim\left(\T^{\mathrm {Fan}(T)}\cap[0,S]\right)\leq\dim\left(\T^{\mathrm {Fan}(S)}\cap[0,S]\right)\leq 1/2.
    \end{equation*}
    Hence
     \begin{equation*}
        \dim(\T^{\mathrm{Fan}(T)})=\dim\big(\cup_{n=1}^\infty (\T^{\mathrm{Fan}(T)}\cap[0,n])\big)=\limsup_{n\to\infty}\dim\big(\T^{\mathrm{Fan}(T)}\cap [0,n]\big)\leq \frac{1}{2}.
    \end{equation*}
    Combining the upper and lower bounds yields $\dim(\T^{\mathrm{Fan}(T)})=1/2$ almost surely.
\end{proof}

\subsection{Dimension of the exceptional points, Theorem \ref{exceptional times} (iv)} The upper bound of $\dim \mathcal S$ is essentially follows by the H\"older continuity of the Brownian excursion. We shall use the following Proposition from \cite[Proposition 2.3]{dimS}.  
\begin{proposition}\label{Holder map upper bound}
    Let $A\subseteq \R^n$ and $g:A\to\R^m$ be a $\gamma$-H\"older continuous function. Then 
    \begin{equation*}
        \dim g(A)\leq \frac{1}{\gamma}\dim A.
    \end{equation*}
\end{proposition}
\begin{lemma}\label{upper bound of S}
    Almost surely, $\dim S\leq 1$.
\end{lemma}
\begin{proof}
    Recall that $\Theta=\sum_{i=1}^\infty \delta_{(t_i,e_i)}$. For any $(t_i,e_i)\in\Theta$, define
    \begin{equation*}
        \T_i:=[t_i,t_i+\xi(e_i)]\cap \T\subseteq \T. 
    \end{equation*}
    Recall that the trajectory of the particle $x_i$ is given by $x_i(t)=\left(e_i(t-t_i)\right)_{t\in[t_i,t_i+\xi(e_i)]}$. Hence 
    \begin{equation}\label{cover of S}
        \mathcal{S}\subseteq \cup_{i=1}^\infty \{(t,x_i(t)):t\in\T _i\}:=\cup_{i=1}^\infty \mathcal{S}_i.
    \end{equation}
    Almost surely, Brownian excursion is $\gamma$-H\"older continuous for all $\gamma<1/2$. Hence the map $t\mapsto (t,x_i(t))$ is also H\"older $\gamma$-continuous map for every $\gamma<1/2$. Therefore by Proposition \ref{Holder map upper bound} we get 
    \begin{equation*}
        \dim \mathcal{S}_i\leq \frac{\dim \T_i}{\gamma}\leq \frac{\dim \T}{\gamma}.
    \end{equation*}
   By taking $\gamma \to 1/2$ and applying Theorem \ref{exceptional times} $(ii)$ we obtain $\dim S_i\leq 1$ almost surely. Therefore by the inclusion given in (\ref{cover of S}) we get $\dim S\leq \dim(\cup_{i=1}^\infty \mathcal{S}_i)\leq 1$ almost surely.
\end{proof}
We now turn to the lower bound of $\dim \mathcal{S}$. Let 
\begin{equation*}
    \mathcal{E}^1:=\left\{e\in\mathcal {E}:\xi(e)\geq1\right\},\quad\Theta^1=\Theta|_{\R^+\times \mathcal{E}^1},\quad\text{and }\Theta^c=\Theta_{\R^+\times\mathcal{E}\setminus\mathcal{E}^1 }
\end{equation*}
By Remark \ref{independent of short time}, $\Theta^1$ is independent of $\T$ and $\Theta^c$ generate the same set of exceptional times as $\Theta$.

The intensity measure of $\Theta^1$ can be written as
    \begin{equation*}
        \rho\n\left(\left\{\xi(e)\geq 1\right\}\right)\dif t\otimes \frac{\n\left(\dif e,\xi(e)\geq 1\right)}{\n(\xi(e)\geq 1)}=:\rho\n\left(\left\{\xi(e)\geq 1\right\}\right)\dif t\otimes Q.
    \end{equation*}
    Since $Q$ is a probability measure supported on $\mathcal{E}^1$ we have the following representation of $\Theta^1$: first sample a Poisson point process of birth times on $\R^+$ with intensity $\rho \n({\xi(e)\geq 1})\dif t$ and then attach independent marks with common distribution $Q$ to its points. Let 
    \begin{equation*}
        \sigma=\inf\{t:(t,e)\in\Theta^1\}
    \end{equation*}
    and denote by $e^1$ the mark corresponding to the first point so that
\begin{equation*}
    (\sigma,e^1)\in\Theta^1.
\end{equation*}
Therefore, $\left(e^1({t-\sigma})\right)_{t\in[\sigma,\sigma+\xi(e)]}$ is the first particle with life time at least one. Moreover, $e^1$ has distribution $Q$ and is independent of $\sigma$ and $\T$. Define 
 \begin{equation}\label{definition of S^1}
     \mathcal{S}^1:=\{(t,e^1(t-\sigma):t\in\T_\sigma\},\quad\text{where }\T_\sigma=[\sigma+1/2,\sigma+1]\cap\T.
 \end{equation}
    Clearly, $\mathcal{S}^1$ is a subset of $\mathcal{S}$. We shall prove that $\dim\mathcal{S}^1\geq 1$ almost surely using energy method.

     By the Markov property the joint law of $e^1(s)$ and $e^1(s+h)$ under $Q$ is given by
    \begin{equation}\label{transition density of Q}
        Q\left(\left\{e^1(s)\in\dif x,e^1(s+h)\in\dif y\right\}\right)=\frac{\n_s(\dif x)P^0_h(x,y)\pr_y[\tau_0>1-s-h]}{\n(\xi(e)\geq 1)}\dif y,
    \end{equation}
  where $\tau_0$ is the first hitting time of zero by Brownian motion. The following estimate will be used to control the energy.
  \begin{lemma}\label{finite expectation under Q}
      For any $\alpha\in(0,1)$ and $h,s>0$ satisfying $s+h<1$,
      \begin{equation*}
          \E_Q\left[\left(h^2+(e^1(s+h)-e^1(s))^2\right)^{-\alpha/2}\right]\leq C_{\alpha}s^{-1/2} h^{-\alpha/2}.
      \end{equation*}
  \end{lemma}
 \begin{proof}
     By (\ref{transition density of Q}) we have
     \begin{equation*}
         \begin{split}
             \E_Q\left[\left(h^2+(e^1(s+h)-e^1(s))^2\right)^{-\alpha/2}\right]=\frac{\int_0^\infty \n_s(\dif x)\int_0^\infty \frac{P_h^0(x,y)\pr[\tau_0>1-s-h]}{\left(h^2+(y-x)^2\right)^{\alpha/2}}\dif y}{\n(\{\xi(e)\geq 1\})}
         \end{split}
     \end{equation*}
     Let
     \begin{equation*}
         P_h(z)=\frac{1}{\sqrt{2 \pi h}}\exp\left(-\frac{z^2}{2h}\right),
     \end{equation*}
     then from the fact that $P^0_h(x,y)\leq P_h(x-y)$ and $\pr[\tau_0>1-s-h]\leq 1$ we obtain
     \begin{equation*}
         \begin{split}
             \E_Q\left[\left(h^2+((e^1(s+h)-e^1(s))^2\right)^{-\alpha/2}\right]&\leq\frac{1}{\n(\{\xi(e)\geq 1\})}\int_0^\infty \n_s(\dif x)\int_0^\infty \frac{P_h(y-x)}{\left(h^2+(y-x)^2\right)^{\alpha/2}}\dif y\\
             &\leq \frac{1}{\n(\{\xi(e)\geq 1\})}\int_0^\infty \n_s(\dif x)\int_\R\frac{P_h(z)}{(h^2+z^2)^{\alpha/2}}\dif z\\
             &=\frac{\n(\{\xi(e)\geq s\})}{\n(\{\xi(e)\geq 1\})}\int_\R \frac{P_h(z)}{(h^2+z^2)^{\alpha /2}}\dif z.
         \end{split}
     \end{equation*}
     Making the change of variable $z=\sqrt{h}u$ we have
     \begin{equation*}
         \int_\R \frac{P_h(z)}{(h^2+z^2)^{\alpha/2}}\dif z=\int_\R \frac{\exp\left(-u^2/2\right)}{\sqrt{2\pi} h^{\alpha/2}(h+u^2)^{\alpha/2}}\dif u\leq \int_\R \frac{\exp\left(-u^2/2\right)}{\sqrt{2\pi} h^{\alpha/2}u^\alpha}\dif u=C_\alpha h^{-\alpha/2}.
     \end{equation*}
     Hence 
     \begin{equation*}
         \E_Q\left[\left(h^2+(e^1(s+h)-e^1(s))^2\right)^{-\alpha/2}\right]\leq \frac{\n(\{\xi(e)\geq s\})}{\n(\{\xi(e)\geq 1\})}C_\alpha h^{-\alpha/2}=C_{\alpha}s^{-1/2}h^{-\alpha/2}
     \end{equation*}
 as desired.  \end{proof}
To apply energy method, we next construct a measure supported on $\T_\sigma$. For any $0<\epsilon<1$ let
\begin{equation*}
    A(\epsilon)=\{t\in [\sigma+1/2,\sigma+1]:\text{ there exists } (s,e)\in\Theta^c \text{ with } t-s\in(0,\epsilon\wedge \xi(e))\}
\end{equation*}
be the set of times $t$ such that there exists a particle born in time interval $(t-\epsilon,t)$ and stay alive at time $t$.  Then $A(\epsilon)$ is an approximation of $\T_\sigma$ in the sense that $\cap_{\epsilon>0}A(\epsilon)=\T_\sigma$. Moreover each $A(\epsilon)$ is open set because all particles have continuous paths. Consequently, almost surely $\T_\sigma$ is a Borel set with $\dim (\T_\sigma)=1/2$. Hence for any $\beta<1/2$ there exists a random measure $\mu_\beta$ (See \cite[Theorem 4.3 (b)]{dimS} for the existence of such measure) supported on $\T_\sigma$ such that 
\begin{equation*}
    \int_{\T_\sigma}\int_{\T_\sigma} \frac{\mu_\beta(\dif t)\mu_\beta(\dif s)}{|t-s|^\beta}<\infty,\quad\text{a.s.}
\end{equation*}
We are now ready to establish the lower bound.
\begin{lemma}\label{lower bound of dim S}
    Almost surely, $\dim \mathcal{S}\geq 1$.
\end{lemma}
\begin{proof}
Recall the definition of $\mathcal{S}^1$ given by (\ref{definition of S^1}). We define a pushforward measure $\nu_\beta$ supported on $\mathcal{S}^1$ by
\begin{equation*}
    \nu_\beta(A)=\mu_\beta(\{t:(t,e^1)\in A\}).
\end{equation*}
Note that $\T_\sigma$ depends only on $\sigma$ and $\T$ whereas $e^1(\cdot)$ is independent of $(\sigma,\T)$ and has distribution $Q$. Conditioning on $\mu_\beta$ and $\sigma$ and applying Lemma \ref{finite expectation under Q} we obtain the $\gamma-$energy of $\nu_\beta$ as follows:
\begin{align*}
    \E\left[\int_{\mathcal{S}^1}\int_{\mathcal{S}^1}\frac{\nu_\beta(\dif x)\nu_\beta(\dif y)}{|x-y|^\gamma}\middle|\mu_\beta,\sigma\right]&=\E\left[\int_{\T_\sigma}\int_{\T_\sigma} \frac{\mu_\beta(\dif t)\mu_\beta(\dif s)}{\left(|t-s|^2+|e^1(t-\sigma)-e^1(s-\sigma)|^2\right)^{\gamma/2}}\middle| \mu_\beta,\sigma\right]\\&=\int_{\sigma+\frac{1}{2}}^{\sigma+1}\int_{\sigma+\frac{1}{2}}^{\sigma+1} \frac{\mu_\beta(\dif t)\mu_\beta(\dif s)}{\E_Q\left[\left(|t-s|^2+|e^1(t-\sigma)-e^1(s-\sigma)|^2\right)^{\gamma/2}\right]}\\&\leq \int_{\sigma+\frac{1}{2}}^{\sigma+1}\int_{\sigma+\frac{1}{2}}^{\sigma+1} C_\gamma\left(t\wedge s-\sigma\right)^{-1/2}\frac{\mu_\beta(\dif t)\mu_\beta(\dif s)}{|t-s|^{\gamma/2}}\\&\leq\int_{\sigma+\frac{1}{2}}^{\sigma+1}\int_{\sigma+\frac{1}{2}}^{\sigma+1}\sqrt{2} C_\gamma\frac{\mu_\beta(\dif t)\mu_\beta(\dif s)}{|t-s|^{\gamma/2}}<\infty,\ {\text{a.s.}}
\end{align*}
provided $\gamma<2\beta$.
Therefore, we get the above energy is almost surely finite under the unconditional probability measure. Hence $\dim \mathcal{S}^1\geq \gamma$ for any $\gamma<1$ as $\beta$ can be chosen arbitrarily close to $1/2$. Finally since $\mathcal{S}^1\subseteq \mathcal {S}$, we conclude that $\dim\mathcal{S}\geq 1$ almost surely. 
\end{proof}
Combining Lemma \ref{upper bound of S}
 and Lemma \ref{lower bound of dim S}, we finish the proof of Theorem{\ref{exceptional times}} $(iv)$.

 \subsection{Nonexistence of exceptional time}
By $(ii)$ of Theorem \ref{exceptional times} and Corollary \ref{exists of T, T(f), and lowerbound} we know that for any $f$ satisfies (\ref{requirement of test function}), almost surely $\dim\left(\T(f)\right)=1/2$ because $\T(f)\subseteq \T$, which provides the upper bound of $\dim\left(\T(f)\right)$. Hence to finish the proof of Theorem \ref{loglog theorem}, it suffices to show that the process $\left(\langle X_t,f\rangle\right)_{t\in\R^+}$ is a continuous process for any even continuous function $f$ such that $f$ is non-decreasing on $[0,\infty)$ with 
\begin{equation}\label{requirement of test function finite version}
    \int_0^\epsilon\frac{f(x)}{x\log(1/x)}\dif x<\infty,\quad\text{for some }\epsilon>0.
\end{equation}
Recall the definition of $\tilde G_k$ given by (\ref{def of tilde G_k}), for any $T>0$ we define a sequence of process by
\begin{equation*}
    Z^k(t):=\sum_{\substack{(t_i,e_i)\in\Theta,\ e_i\in\tilde G_k\\ 0\leq t-t_i\leq \xi(e_i)}}f(|e_i(t-t_i)|),\quad 0\leq t\leq T.
\end{equation*}
Hence $\langle X_t,f\rangle=\sum_{k=1}^\infty Z^k(t)$. According to Remark \ref{no cluster point except zero}, we know that the set $\{(t,e)\in\Theta:e\in\tilde G_k,\ t\in[0,T]\}$ is almost surely finite. Thus $\left(Z^k(t)\right)_{k\geq 1}$ is a sequence of continuous processes. Therefore it is enough to show that
\begin{equation*}
    \sum_{k=2}^\infty\sup_{t\in[0,T]}Z^k(t)<\infty\quad\text{a.s.}
\end{equation*}
so that $Z^k(t)$ converge uniformly, which implies that $\langle X_t,f\rangle$ is continuous on $[0,T]$. For simplicity we let $T=1$, the proof can be easily generalized to any $T>0$. For any $t\in[0,1]$ and any closed interval $I\subseteq[0,1] $, let
\begin{align*}
        N_k(t)&=\Card\left(\left\{(s,e)\in\Theta:e\in\tilde G_k,\ s<t<s+\xi(e)\right\}\right),\\ N_k(I)&=\Card\left(\left\{(s,e)\in\Theta:e\in\tilde G_k,s\in I\right\}\right).
\end{align*}
Then both $N_k(t)$ and $N_k(I)$ are Poisson random variables.
We firstly need an estimation of $\E[N_k(t)]$ and $\E[N_k(I)]$.
\begin{lemma}\label{upperbound of Poi parameter}
    For any $k\geq 2$, $t\in[0,1]$ and $I=[a,b]\subseteq[0,1]$ with $|I|=h_{k}^2$, we have
    \begin{equation*}
        \E[N_k(t)]\leq {C_\rho}h_{k-1},\quad\E[N_k(I)]\leq \rho h_{k},
    \end{equation*}
    where $C$ is a positive constant.
\end{lemma}
\begin{proof}
    For $N_k(t)$, we have
    \begin{align*}
        \E[N_k(t)]&=\rho\int_0^t\n\left(\left\{e\in\tilde G_k,\xi(e)\geq t-s\right\}\right)\dif s\\
        &=\rho\int_0^t \n\left(\left\{e\in\tilde G_k,\xi(e)\geq s\right\}\right)\dif s\\
        &=\rho\int_0^t\int_{\mathcal{E}}\ind_{\{e\in\tilde G_k\}}\ind_{\{\xi(e)\geq s\}}\n(\dif e)\dif s\\
        &\leq \rho\int_{\mathcal{E}}\int_0^\infty \ind_{\{\xi(e)\geq s\}}\dif s\ind_{\{e\in\tilde G_k\}}\n(\dif e)\\
        &=\rho\int_{\mathcal{E}}\xi(e)\ind_{\{e\in\tilde G_k\}}\n(\dif e)\\
       &\overset{\tilde{e}(u)=h_{k-1}^{-1}e(h_{k-1}^2 u)}{=}\rho\int_{\mathcal{E}}h_{k-1}\xi(\tilde{e})\ind_{\{\tilde{e}:\max_{0\leq u\leq \xi(\tilde{e})}|\tilde e(u)|\in[h_{k}/h_{k-1},1]\}}\n(\dif \tilde e)\\
       &\leq h_{k-1}\rho\int_{\mathcal{E}}\xi(\tilde e)\ind_{\{\tilde{e}:\max_{0\leq u\leq \xi(\tilde{e})}|\tilde e(u)|\leq 1\}}\n(\dif \tilde e),
    \end{align*}
    where in the last equality we use the scaling property of It\^o's Brownian excursion measure given by Proposition \ref{Brownian excursion scaling}. Note that
    \begin{align*}
        \int_{\mathcal{E}}\xi(e)\ind_{\{e:\max_{0\leq u\leq \xi(e)}|e(u)|\leq 1\}}\n(\dif e)&=\int_\mathcal{E}\int_0^\infty \ind_{\{\xi(e)\geq u\}}\dif u\ind_{\{e:\max_{0\leq u\leq \xi(e)}|e(u)|\leq 1\}}\n(\dif e)\\
        &=\int_0^\infty \n\left(\{\xi(e)\geq u,\max_{0\leq u\leq \xi(e)}|e(u)|\leq 1\}\right)\dif u\\
        &\leq \int_0^\infty \n\left(\left\{\xi(e)\geq u,|e(u)|\leq 1\right\}\right)\dif u\\
        &=\int_0^\infty\int_{-1}^1 \n_t(\dif x)\dif t<\infty,
    \end{align*}
    where $\n_t(\dif x)$ is the entrance law given by (\ref{entrance law}). Hence $\rho\dif t\otimes \n\left(N_k(t)\right)\leq {C_\rho}h_{k-1}$ with
    \begin{equation*}
         C_\rho=\rho\int_0^\infty\int_{-1}^1 \n_t(\dif x)\dif t.
    \end{equation*}
    For $N_k(I)$, by (\ref{upper bound of n(tilde G_k)}) we can directly compute
    \begin{equation*}
        \E[N_k(I)]=\rho\n(\tilde G_k)|I|\leq \rho h_k^{-1}h_k^2=\rho h_k
    \end{equation*}
    as desired.
\end{proof}
We close this section by finishing the proof of Theorem \ref{loglog theorem}.

\begin{proof}[Proof of Theorem \ref{loglog theorem}]
    As we discussed in the beginning of this subsection, $(ii)$ and sufficiency of $(i)$ are direct consequences of Theorem \ref{exceptional times} $(ii)$ and Corollary \ref{exists of T, T(f), and lowerbound}. Hence we only need to prove the necessity of Theorem \ref{loglog theorem} $(i)$, which is enough to show 
    \begin{equation*}
        \sum_{k\geq 2}\sup_{t\in[0,1]}Z^k(t)<\infty\quad\text{a.s.}
    \end{equation*}
   From now on, we take $h_k=\exp(-\exp(k))$ and cut the interval $[0,1]$ into subintervals of the form $I_k^j=[s_k^{j-1},s_k^{j}]$, where 
    \begin{equation*}
        s_k^j=jh_k^2\quad\text{for }j=1,\cdots,\lceil h_k^{-2}\rceil.
    \end{equation*}
    Hence 
    \begin{equation}\label{upper bound of Z_^k(s)}
        \begin{split}
            \sup_{s\in I^j_k}Z^k(s)&\leq\left(N_k(s_k^{j-1})+N_k(I_k^j)\right)f(h_{k-1}),\\
            \sup_{s\in[0,1]}Z^k(s)&\leq \left(\max_{j=1,\cdots,\lceil h_k^{-2}\rceil}N_k(s_k^{j-1})+\max_{j=1,\cdots,\lceil h_k^{-2}\rceil}N_k(I_k^j)\right)f(h_{k-1}).
        \end{split}
    \end{equation}
    Note that $N_k(s_k^{j-1})$ is a Poisson random variable with parameter $\lambda_{k,j}\leq Ch_{k-1}$ by Lemma \ref{upperbound of Poi parameter}. Hence 
    \begin{equation*}
        \pr[N_k(s^j_k)\geq 6]=\sum_{l=6}^\infty \frac{\lambda_{k,j}^l \exp(-\lambda_{k,j})}{l!}\leq Ch_{k-1}^6.
    \end{equation*}
    Therefore, by the union bound, we get
    \begin{equation*}
        \pr\left[\max_{j=1,\cdots,\lceil h_k^{-2}\rceil}N_k(s_k^{j-1})\geq 6\right]\leq Ch_{k}^{-2}h_{k-1}^6=C_4\exp\left((2\exp(1)-6)\exp(k-1)\right),
    \end{equation*}
    which is summable. Applying Borel-Cantelli lemma, almost surely for large enough $k$ we get $\max_{j=1,\cdots,\lceil h_k^{-2}\rceil}N_k(s_k^{j-1})<6$. Similarly, by Lemma \ref{upperbound of Poi parameter}  we have
    \begin{equation*}
        \pr\left[N_k(I^j_k)\geq 6\right]\leq Ch_k^6\quad\text{for some }C>0.
    \end{equation*}
    Hence 
    \begin{equation*}
        \pr\left[\max_{j=1,\cdots,\lceil h_{k}^{-2}\rceil}N_k(I_k^j)\geq 6\right]\leq Ch_k^6h_k^{-2},
    \end{equation*}
    which is summable. Hence almost surely for large enough $k$ we get $\max_{j}N_k(I_k^j)<6$. Thus by (\ref{upper bound of Z_^k(s)}) there exists an almost surely finite number $N$ such that
    \begin{equation*}
        \sup_{s\in[0,1]}Z^k(s)\leq 10f(h_{k-1}),\quad\text{for all }k\geq N.
    \end{equation*}
    It remains to show $f(h_{k-1})$ is summable for large enough $k$. Note that for any even function $f$ that is non-decreasing on $[0,\infty)$ satisfies (\ref{requirement of test function finite version}), we get
    \begin{equation*}
        \int_0^\epsilon \frac{f(x)}{x\log(1/x)}\dif x=\int_{\log\epsilon^{-1}}^\infty \frac{f(\exp(-y))}{y}\dif y<\infty.
    \end{equation*}
    Choose $N$ large enough such that $\exp(k)\geq\log\epsilon^{-1}$ for all $k\geq N$. Then we have
    \begin{equation*}
        \infty>\int_{\log\epsilon^{-1}}^\infty \frac{f(\exp(-y))}{y}\dif y\geq \sum_{k=N}^\infty \int_{e^k}^{e^{k+1}}\frac{f(\exp(-y))}{y}\dif y\geq \sum_{k=N}^\infty f(h_{k+1}).
    \end{equation*}
    Hence $f(h_{k-1})$ is summable, which complete the proof.
\end{proof}

\section{Power sum of particles}\label{proof of theorem 3}
\subsection{H\"older continuity, Theorem \ref{Z_t^alpha} (i)}
The main idea for proving (i) is to apply Kolmogorov continuity theorem by controlling
\begin{equation*}
    \sup_{s\in[0,T-\Delta t]}\E[|Z^\alpha_{s+\Delta t}-Z^\alpha_s|^k].
\end{equation*}
For any $T\in\R^+$, even number $k$ and $\Delta t$ small enough. The first step is to derive the moment generating function of $Z^\alpha_{s+\Delta t}-Z^\alpha_{s}$. Recall that $X_s=\sum_{i=1}^{N_s}\delta_{x_i(s)}$, where $N_s$ is the total number of alive particles at time $s$. For any $t>0$ we shall divide $X_{s+t}$ into two parts. Let
\begin{equation*}
    X^{\text{old}}_{s+t}=\sum_{i=1}^{N_s}\delta_{x_i(s+t)}
\end{equation*}
be the point process of the alive particles at time $s$ evolving up to time $s+t$, which is the first part of our decomposition. To see the distribution of $X^{\text{old}}_{s+t}$, we can condition on $X_s$ and let $(x_i)_{i=1}^{N_s}$ be the positions of alive particles at time $s$. Let $y_i$ be distributed as $\bar{P}_t^0(x_i,\cdot)$ independently for different $i$,where $\bar{P}_t^0$ is the transition density of Brownian motion absorbed at zero given by 
\begin{equation*}
    \bar{P}^0_t(x,\dif y)=P_t^0(x,\dif y)+\erfc\left(\frac{|x|}{\sqrt{2t}}\right)\delta_0(\dif y).
\end{equation*}
Then $x_i(s+t)$ has the same distribution as $y_i$ for all $i=1,\cdots,N_s$. The second part  is newborn particles during time $[s,s+t]$ given by
\begin{equation*}
    X_{s+t}^{\text{new}}=\sum_{i=1}^\infty \delta_{e_i(s+t-t_i)}\ind_{\{s\leq t_i\leq s+t\}},
\end{equation*}
where $(t_i,e_i)\in\Theta$. Therefore we get an old-new decomposition of the point process for the alive particles:
\begin{equation}\label{decomposition of X_s+t}
    X_{s+t}=X^{\text{old}}_{s+t}+X^{\text{new}}_{s+t}.
\end{equation}
The advantages of this decomposition is that $X^{\text{old}}_{s+t}$ is independent of $X^{\text{new}}_{s+t}$ as the evolution of existing particles is independent of newborn particles according to the construction of our model. Moreover, note that by Proposition \ref{Poi X_t}, $X_s$ is a Poisson point process with intensity $\rho\erfc(|x|/\sqrt{2s})\dif x$. Hence by marking theorem for Poisson point process (see \cite[Theorem 5.6]{Marking_theorem} for details) the joint law of $(X_s,X_{s+t}^\text{old})$ is another Poisson point process on $\R\times \R$ with intensity measure 
\begin{equation}\label{intensity of old}
    \rho\erfc(|x|/\sqrt{2s})\bar{P}^0_t(x,\dif y)\dif x.
\end{equation}
Thus the moment generating function of $Z_{s+t}^\alpha-Z^\alpha_s$ can be treated as Laplace functional of two independent Poisson point processes.
\begin{lemma}
    For any $s,t>0$, the moment generating function of $Z^\alpha_{s+t}-Z^\alpha_s$ is given by
    \begin{equation*}
        M_{s,t}(\lambda):=\E[\exp\left({\lambda(Z^\alpha_{s+t}-Z^\alpha_{s})}\right]=\exp\left({I_{s,t}(\lambda)}\right),
    \end{equation*}
    where 
    \begin{equation}\label{formula of I_{s,t}}
    \begin{split}
        I_{s,t}(\lambda)&=\rho\int_\R\left(\exp\left({\lambda |x|^\alpha}\right)-1\right)\erfc\left(\frac{|x|}{\sqrt{2t}}\right)\dif x\\&+\rho\int_\R\int_\R\left(\exp\left({\lambda(|y|^\alpha-|x|^\alpha)}\right)-1\right)\erfc\left(\frac{|x|}{\sqrt{2s}}\right)\bar{P}_t^0(x,\dif y)\dif x
        \end{split}
    \end{equation}
\end{lemma}
\begin{proof}
    According to the definition of $(Z_u^\alpha)_{u\geq0}$ we have
    \begin{equation*}
        Z_{t+s}^\alpha-Z^\alpha_s=\langle X_{s+t}^{\text{new}},x^\alpha\rangle+\langle X_{s+t}^{\text{old}},x^\alpha\rangle-\langle X_s,x^\alpha\rangle.
    \end{equation*}
    Since $(X_{s+t}^{\text{old}},X_s)$ is a Poisson point process $\Phi$ on $\R\times\R$ whose intensity measure is given by (\ref{intensity of old}), we get the Laplace functional of $\Phi$ for the function $\varphi(x,y)=\lambda(|y|^\alpha-|x|^\alpha)$ is given by 
    \begin{equation*}
    \begin{split}
                \E\left[\exp\left({\int_{\R^2}\varphi(x,y)\dif\Phi}\right)\right]&=\E\left[\exp\left({\lambda(\langle X_{t+s}^{\text{old}},|x|^\alpha\rangle-\langle X_s,|x|^\alpha\rangle)}\right)\right]\\
                &=\E\left[\exp\left(\rho\int_\R\int_\R(\exp\left({\lambda(|y|^\alpha-|x|^\alpha)}\right)-1)\erfc\left(\frac{|x|}{\sqrt{2s}}\right)\bar{P}_t^0(x,\dif y)\dif x\right)\right].
    \end{split}
    \end{equation*}
    Similarly, by Proposition \ref{Poi X_t} we get $X_{s+t}^\text{new}$ is a Poisson point process with intensity measure $\rho\erfc(x/\sqrt{2t})$. Thus we obtain
    \begin{equation*}
        \E[\exp\left(\lambda\langle X_t^\text{new},|x|^\alpha\rangle\right)]=\exp\left(\rho\int_\R(\exp\left({\lambda |x|^\alpha}\right)-1)\erfc\left(\frac{|x|}{\sqrt{2t}}\right)\dif x\right).
    \end{equation*}
    Since $X_{s+t}^\text{new}$ is independent of $(X_{s+t}^{\text{old}},X_s)$, combining above two equalities we get
    \begin{equation*}
        \begin{split}
            M_{s,t}(\lambda)&=\E\left[\exp\left({\lambda(\langle X_{s+t}^{\text{old}},|x|^\alpha\rangle-\langle X_s,|x|^\alpha\rangle)}\right)\exp\left(\lambda\langle X_{s+t}^\text{new},|x|^\alpha\rangle\right)\right]
            \\&=\E\left[\exp\left({\lambda(\langle X_{s+t}^{\text{old}},|x|^\alpha\rangle-\langle X_s,|x|^\alpha\rangle)}\right)\right]\E\left[\exp\left(\lambda\langle X_{s+t}^\text{new},|x|^\alpha\rangle\right)\right]\\
            &=\exp\left(I_{s,t}(\lambda)\right),
        \end{split}
    \end{equation*}
    where $I_{s,t}(\lambda)$ has the desired formula.
\end{proof}
Now to compute even moments of $Z_{s+\Delta t}^\alpha-Z_{s}^\alpha$ we need to take the derivative of $M_{s,\Delta t}(\lambda)$ with respect to $\lambda$ and evaluate it at $\lambda=0$. In the rest of this section, for any function $f$ of $\lambda$, we use 
\begin{equation*}
    f^{(k)}(\lambda)
\end{equation*}
to represent the $k-$th derivative of $f$ with respect to $\lambda$. Note that $M_{s,t}(0)=1$ for all $s,t\in\R^+$ by By Fa\`a di Bruno's formula (see e.g. \cite{Faa})
\begin{equation}\label{k-th moement of Z_s+t-Z_s}
M_{s,t}^{(k)}(0)
=
\sum_{\substack{m_1,\dots,m_k\ge 0\\
m_1+2m_2+\cdots+km_k=k}}
\frac{k!}{m_1!m_2!\cdots m_k!}
\prod_{j=1}^k
\left(\frac{I_{s,t}^{(j)}(0)}{j!}\right)^{m_j}.
\end{equation}
Hence it is enough to compute $|I^{(k)}_{s,t}(0)|$ for small $t$.
\begin{lemma}\label{estimation of I^k_{s,t}}
    For any $s>0$ and $k\in\N$, $|I^{(k)}_{s,t}(0)|\leq C_{\rho,k,\alpha}(1+\sqrt{s})t^{\alpha k/2}$ for all $0\leq t\leq 1$.
\end{lemma}
\begin{proof}
    By the formula of $I_{s,t}(\lambda)$ given by (\ref{formula of I_{s,t}}) it is easy to check that 
    \begin{equation}\label{decomposition of I^k_s,t}
        I^{(k)}_{s,t}(0)=I_1^{(k)}+I_2^{(k)}+I_3^{(k)},
    \end{equation}
    where
    \begin{equation}\label{I1,I2,I3}
        \begin{split}
            I_1^{(k)}&=\rho\int_\R |x|^{k\alpha}\erfc\left(\frac{|x|}{\sqrt{2t}}\right)\dif x,\\
            I_2^{(k)}&=\rho\int_\R\int_\R (|y|^\alpha-|x|^\alpha)^k\erfc\left(\frac{|x|}{\sqrt{2s}}\right)\\
            &\times\frac{1}{\sqrt{2\pi t}}\left(\exp\left(-\frac{(y-x)^2}{2t}\right)-\exp\left(-\frac{(y+x)^2}{2t}\right)\right)\ind_{\{xy>0\}}\dif y\dif x\\
            I^{(k)}_3&=\rho(-1)^k\int_\R |x|^{k\alpha}\erfc\left(\frac{|x|}{\sqrt{2t}}\right)\erfc\left(\frac{|x|}{\sqrt{2s}}\right)\dif x.
        \end{split}
    \end{equation}
    For $I_1^{(k)}$, by changing of variable we get
    \begin{equation*}
    \begin{split}
        I_1^{(k)}=2\rho\int_0^\infty x^{k\alpha}\erfc\left(\frac{x}{\sqrt{2t}}\right)\dif x=\rho t^{\frac{1+k\alpha}{2}}2^{\frac{3+k\alpha}{2}}\int_0^\infty y^{k\alpha}\erfc(y)\dif y\leq C_{\rho,k,\alpha}t^{\frac{1+k\alpha}{2}}.
    \end{split}
    \end{equation*}
    for some constant $C_{\rho,k,\alpha}$.
    For $I^{(k)}_2$ we have
    \begin{equation*}
    \begin{split}
        |I_2^{(k)}|&\leq2\rho\int_0^\infty \erfc\left(\frac{x}{\sqrt{2s}}\right)\dif x\int_0^\infty\frac{1}{\sqrt{2\pi t}}|y^\alpha-x^\alpha\bigr|^k\exp\left(-\frac{(y-x)^2}{2t}\right)\dif y\\
        &=\frac{2\rho}{\sqrt{2\pi}}t^{\frac{\alpha k+1}{2}}\int_0^\infty \erfc\left(\frac{u\sqrt{t}}{\sqrt{2s}}\right)\dif u \int_0^\infty |v^\alpha-u^\alpha|^k \exp\left(-\frac{(u-v)^2}{2}\right)\dif v\\
        &\leq\frac{2\rho}{\sqrt{2\pi}}t^{\frac{\alpha k+1}{2}}\int_0^\infty \erfc\left(\frac{u\sqrt{t}}{\sqrt{2s}}\right)\dif u \int_0^\infty |v-u|^{\alpha k}\exp\left(-\frac{(u-v)^2}{2}\right)\dif v\\
        &=\frac{2\rho\sqrt{s}}{\sqrt{\pi}}t^{\frac{\alpha k}{2}}\int_0^\infty \erfc(w_1)\dif w_1\int_{-\frac{\sqrt{2s}w_1}{\sqrt{t}}}^\infty|w_2|^{\alpha k}\exp\left(-\frac{w_2^2}{2}\right)\dif w_2\\
        &\leq2\rho t^{\frac{\alpha k}{2}}\sqrt{\frac{s}{\pi}}\int_0^\infty \erfc(w_1)\dif w_1\int_\R |w_2|^{\alpha k}\exp\left(-\frac{w_2^2}{2}\right)\dif w_2\\
        &=C_{\rho, k,\alpha}\sqrt{s}t^{\frac{\alpha k}{2}}
    \end{split}
    \end{equation*}
    for all $t\geq 0$.  For $I_3^{(k)}$ we have
    \begin{equation*}
        |I^{(k)}_3|\leq I_1^{(k)}\leq C_{\rho,k,\alpha} t^{\frac{1+\alpha k}{2}}\leq C_{\rho,k,\alpha}t^{\frac{\alpha k}{2}}
    \end{equation*}
    for all $t\leq 1$. Therefore, 
    \begin{equation*}
    |I_{s,t}^{(k)}|\leq |I_1^{(k)}|+|I_2^{(k)}|+|I_3^{(k)}|\leq C_{\rho,k,\alpha}(1+\sqrt{s})t^{\frac{\alpha k}{2}}
    \end{equation*}
    for $t\geq  0$.
\end{proof}

\begin{proof}[proof of Theorem \ref{Z_t^alpha} $(i)$]
According to (\ref{k-th moement of Z_s+t-Z_s}) and Lemma \ref{estimation of I^k_{s,t}}, for $\Delta t\geq 0$ and $T\geq 1$:
\begin{equation}\label{asymptotic of Z_s+t-Z_t}
\begin{split}
    \sup_{s\in[0,T-\Delta t]}\E[&|Z_{s+\Delta t}^\alpha-Z_{s}^\alpha|^k]=\sum_{\substack{m_1,\dots,m_k\ge 0\\
m_1+2m_2+\cdots+km_k=k}}
\frac{k!}{m_1!m_2!\cdots m_k!}
\prod_{j=1}^k
\left(\frac{I_{s,\Delta t}^{(j)}(0)}{j!}\right)^{m_j}\\
&\leq \sum_{\substack{m_1,\dots,m_k\ge 0\\
m_1+2m_2+\cdots+km_k=k}}\frac{k!}{m_1!m_2!\cdots m_k!}\prod_{j=1}^k\left(C_{\rho,j,\alpha}(1+\sqrt{s})(\Delta t)^{\alpha j/2}\right)^{m_j}\\
&\leq \sum_{\substack{m_1,\dots,m_k\ge 0\\
m_1+2m_2+\cdots+km_k=k}}\frac{k!}{m_1!m_2!\cdots m_k!}\left(C_{\rho,j,\alpha}\sqrt{T}\right)^{\sum_{j=1}^k m_j}\left(\Delta t\right)^{\sum_{j=1}^k\alpha j m_j/2} \\
&\leq \sum_{\substack{m_1,\dots,m_k\ge 0\\
m_1+2m_2+\cdots+km_k=k}}\frac{k!}{m_1!m_2!\cdots m_k!}\left(C_{\rho,j,\alpha}\sqrt{T}\right)^{\sum_{j=1}^k jm_j}\left(\Delta t\right)^{\sum_{j=1}^k\alpha j m_j/2}\\
&\leq C_{\rho, k,\alpha}T^{k/2}(\Delta t)^{\alpha k/2}.
\end{split}
\end{equation}
Hence by Kolmogorov continuity theorem we get a modified process $(\tilde Z_t^\alpha)_{t\in[0,T]}$ such that almost surely $(\tilde Z_t^\alpha)_{t\in[0,T]}$ is $\gamma-$H\"older continuous for all $\gamma\in(0,\frac{\alpha}{2}-k^{-1})$. Moreover, let $J$ be a countable dense subset of $[0,T]$, then almost surely $\tilde Z_{t_j}^\alpha= Z_{t_j}^\alpha$ for all $t_j\in J$ and $\tilde Z_t^\alpha=\lim_{t_j\to t}Z^\alpha_{t_j}$ for any $t$ with $t_j\in J$ and $t_j\to t$. Since $f_\alpha(x)=|x|^\alpha$ satisfies the finite-integral condition given by (\ref{requirement of test function finite version}), the original process $(Z_t^\alpha)_{t\in[0,T]}$ is a continuous process. Therefore almost surely for any $t$ with $t_j\in J$ and $t_j\to t$ we get
\begin{equation*}
    \tilde Z_t^\alpha=\lim_{t_j\to t}Z_{t_j}^\alpha=Z_t^\alpha.
\end{equation*}
Hence we get the$Z_t^\alpha$ is $\gamma$-H\"older continuous for all $\gamma\in(0,\frac{\alpha}{2}-k^{-1})$ and all $k$ large enough such that $\frac{\alpha}{2}-k^{-1}\geq 0$. By taking $k\to\infty$ we get the desired result.
\end{proof}

\subsection{Law of large numbers and central limit theorem, Theorem \ref{Z_t^alpha} (ii) \& (iii)} 
To derive the law of large numbers for $Z_t^\alpha$, we first compute moments of $Z_t^\alpha$. According to the decomposition of $X_{s+t}$ given by (\ref{decomposition of X_s+t}), if we choose $s=0$ then $X_t=X_t^{\text{new}}$ because there is no particle at the beginning. Hence we get the moment generating function for $Z_t^\alpha$ by taking $s\to0$ in (\ref{formula of I_{s,t}}). The second part of (\ref{formula of I_{s,t}}) vanishes as $s\to 0$, therefore
\begin{equation}\label{MGF of Z_t^alpha}
    \E[\exp(\lambda Z_t^\alpha)]=M_{0,t}(\lambda)=\exp\left(I_{0,t}(\lambda)\right).
\end{equation}
Note that by (\ref{decomposition of I^k_s,t}) and (\ref{I1,I2,I3}), when $s\to0$, $I_2^{(k)},I_3^{(k)}\to0$ as well. Hence for $k\geq 1$ 
\begin{equation*}
    I_{0,t}^{(k)}=I_1^{(k)}.
\end{equation*}
In particular,
\begin{equation}\label{expecation of Z_t^alpha}
    \E[Z_t^\alpha]=I^{(1)}_{0,t}(0)=\rho\int_\R |x|^\alpha \erfc\left(\frac{|x|}{\sqrt{2t}}\right)\dif x=D_\alpha t^{(1+\alpha)/2},
\end{equation}
where $D_\alpha$ is the constant given by $(ii)$ of Theorem \ref{Z_t^alpha}. Using the moment generating function we can get the law of large number for $(Z_n^\alpha)_{n\in\N}$.
\begin{lemma}\label{LLN for discrete time Z_n}
    Almost surely, $Z^\alpha_n/\E[Z_n^\alpha]\to 1$ as $n\to \infty$.
\end{lemma}
\begin{proof}
The moment generating function of $Y_n:=Z_n^\alpha-\E[Z_n^\alpha]$ can be derived by using (\ref{MGF of Z_t^alpha}):
\begin{equation*}
    \E[\exp(\lambda Y_n)]=\exp\left(H_n(\lambda)\right),\quad H_n(\lambda)=I_{0,n}(\lambda)-\lambda \E[Z_n^\alpha].
\end{equation*}
Applying Fa\`a di Bruno's formula again we get the 6-th moment of $Y_n$:
\begin{equation}\label{6-th moment of Y_n}
    \E[Y_n^6]=
\sum_{\substack{m_1,\dots,m_k\ge 0\\
m_1+2m_2+\cdots+6m_6=6}}
\frac{6!}{m_1!m_2!\cdots m_6!}
\prod_{j=1}^6
\left(\frac{H_n^{(j)}(0)}{j!}\right)^{m_j}.
\end{equation}
Note that $H'_n(0)=0$ as $\E[Y_n]=0$ and by (\ref{decomposition of I^k_s,t}) we get for $k\geq 2$
\begin{equation*}
    H_n^{(k)}(0)=I^{(k)}_{0,n}(0)=I^{(k)}_1=C_\alpha(k)n^{(1+k\alpha)/2},
\end{equation*}
where $C_\alpha(k)=\rho2^{(3+k\alpha)/2}\int_0^\infty \erfc(x)x^{k\alpha}\dif x$. Therefore the only non zero term in (\ref{6-th moment of Y_n}) are
\begin{equation*}
    (m_1,\cdots,m_n)\in\left\{(0,1,0,1,0,0),(0,3,0,0,0,0),(0,0,2,0,0,0),(0,0,0,0,0,1)\right\}.
\end{equation*}
Hence
\begin{align*}
    \E[Y_n^6]&=I_1^{(6)}+15I_1^{(4)}I_1^{(2)}+10\left(I_1^{(3)}\right)^2+15\left(I_1^{(2)}\right)^3\\
    &=C_{\alpha}(6)n^{(1+6\alpha)/2}+15C_{\alpha}(4)C_{\alpha}(2)n^{(2+6\alpha)/2}\\&+10\left(C_{\alpha}(3)\right)^2n^{(2+6\alpha)/2}+15\left(C_\alpha(2)\right)^3n^{(3+6\alpha)/2}.
\end{align*}
Now by Markov inequality we get
\begin{equation*}
    \pr\left[\left|\frac{Z_n^\alpha-\E[Z_n^\alpha]}{\E[Z_n^\alpha]}\right|>\epsilon\right]\leq\frac{\E[Y_n^6]}{\epsilon^6\E[Z_n^\alpha]^6}=\frac{\E[Y_n^6]}{\epsilon^6D_\alpha^6n^{3(1+\alpha)}}=O(n^{-3/2}),\quad\text{as }n\to\infty,
\end{equation*}
which is a summable. Therefore, by Borel-Cantelli lemma we get $Z_n^\alpha/\E[Z_n^\alpha]\to1$ almost surely.
\end{proof}
The next step is to show $\left((Z_t^\alpha-Z_n^\alpha)/D_\alpha n^{(1+\alpha)/2}\right)_{t\in[n,n+1)}\to0$ as $n\to \infty$. To do this we need an estimation of the H\"older constant for $Z_t^\alpha$. We recall the following version of Garsia-Rodenmich-Rumsey inequality:
\begin{lemma}\label{GRR}
    For any $m>1$, $\delta>1/m$ and $f\in C[0,1]$, there exists a constant $C_{\delta,m}$ that is independent of $f$ such that 
    \begin{equation*}
        \sup_{t,s\in[0,1],t\neq s}\left(\frac{|f(t)-f(s)|}{|t-s|^{\delta-1/m}}\right)^m\leq C_{\delta,m}\int_0^1\int_0^1\frac{|f(t)-f(s)|^m}{|t-s|^{\delta m+1}}\dif t\dif s.
    \end{equation*}
\end{lemma}
This lemma is slightly different from the original version of Garsia-Rodenmich-Rumsey inequality derived in \cite[Lemma 1.1]{GRR}. By taking $\Psi(x)=x^m$ and $p(u)=u^{\delta+1/m}$ in \cite[Lemma 1.1]{GRR} we get this version of inequality.

\begin{proof}[proof of Theorem \ref{Z_t^alpha} $(ii)$]
    Note that by Theorem \ref{Z_t^alpha} $(i)$, we know that the sample paths of $Z_t^\alpha$ is $\gamma$-H\"older continuous for all $\gamma\in(0,\alpha/2)$. Let
    \begin{equation*}
        K_{\gamma,n}=\sup_{t,s\in[0,1],t\neq s}\frac{|Z^\alpha_{n+t}-Z^\alpha_{n+s}|}{|t-s|^\gamma}
    \end{equation*}
    be the H\"older constant. By setting $f(t)=Z_{n+t}^\alpha$ and $\delta=\gamma+1/m$ in Lemma \ref{GRR} and applying (\ref{asymptotic of Z_s+t-Z_t}) we get
    \begin{equation*}
        \begin{split}
        \E[K_{\gamma,n}^m]&\leq C_{\delta,m}\int_0^1\int_0^1\frac{\E[|Z^\alpha_{n+t}-Z^\alpha_{n+s}|^m]}{|t-s|^{\delta m+1}}\dif t \dif s\\
        &\leq 2C_{\delta,m}\int_0^1\int_0^{1-r}\frac{\E[|Z_{n+r+u}^\alpha-Z_{n+r}^\alpha|^m]}{u^{\gamma m+2}}\dif u\dif r\\
        &\leq2C_{\delta,m}\int_0^1\int_0^{1-r}\frac{C_{\rho,m,\alpha} (n+r)^{m/2}u^{\alpha m/2}}{u^{\gamma m+2}}\dif u \dif r\\
        &=C_{\gamma,\rho,m,\alpha}\int_0^1\int_0^{1-r}(n+r)^{m/2} u^{m(\alpha/2-\gamma)-2}\dif u\dif r\\
        &\leq C_{\gamma,\rho,m,\alpha}(n+1)^{m/2}\int_0^1\int_0^{1-r} u^{m(\alpha/2 -\gamma)-2}\dif u\dif r.
        \end{split}
    \end{equation*}
    Since $\alpha/2>\gamma$, for large enough $m$ we get $m(\alpha/2-\gamma)-2>0$, which implies that the above integral is finite. Thus 
    \begin{equation}\label{expectation of holder constant}
        \E[K_{\gamma,n}^m]\leq C_{\gamma,\rho,m,\alpha}(n+1)^{m/2}.
    \end{equation}
    By the H\"older continuity of $Z_t^\alpha$ we get
    \begin{equation*}
        \frac{|Z_t^\alpha-Z_n^\alpha|}{D_\alpha n^{(1+\alpha)/2}}\leq\frac{K_{\gamma,n}|t-n|^\gamma}{D_\alpha n^{(1+\alpha)/2}}\leq \frac{K_{\gamma,n}}{D_\alpha n^{(1+\alpha)/2}},\quad\text{for any }t\in[n,n+1).
    \end{equation*}
    So by (\ref{expectation of holder constant}) and Markov inequality we get
    \begin{align*}
        \pr\left[\sup_{t\in[n,n+1)}\frac{|Z^\alpha_t-Z^\alpha_n|}{D_\alpha n^{(1+\alpha)/2}}>\epsilon\right]\leq\frac{\E[K_{\gamma,n}^m]}{\epsilon^mD_\alpha^m n^{(1+\alpha)m/2}}\leq \frac{C_{\gamma, \rho,m,\alpha}(n+1)^{m/2}}{\epsilon^m D_\alpha^m n^{(1+\alpha)m/2}}\approx n^{-\alpha m /2},
    \end{align*}
    which is summable in $n$ for large enough $m$. Therefore, by Borel-Cantelli lemma and Lemma \ref{LLN for discrete time Z_n} we get 
    \begin{equation*}
        \frac{Z^\alpha_t}{D_\alpha \lfloor t\rfloor^{(1+\alpha)/2}}\to 1,\quad\text{a.s.},
    \end{equation*}
    where $\lfloor\cdot\rfloor$ is the floor function. Hence
    \begin{align*}
        \left|\frac{Z^\alpha_t}{D_\alpha t^{(1+\alpha)/2}}-\frac{Z^\alpha_t}{D_\alpha \lfloor t\rfloor^{(1+\alpha)/2}}\right|&=\frac{Z^\alpha_t}{D_\alpha \lfloor t\rfloor^{(1+\alpha)/2}}\left|\frac{\lfloor t\rfloor^{(1+\alpha)/2}-t^{(1+\alpha)/2}}{t^{(1+\alpha)/2}}\right|\\
        &\leq \frac{Z^\alpha_t}{{D_\alpha \lfloor t\rfloor^{(1+\alpha)/2}}}\frac{1}{t^{(1+\alpha)/2}}\to0
    \end{align*}
    as $t\to \infty$. Thus we get the desired law of large number for $Z_t^\alpha$.
\end{proof}
In the last part of this section, we show the central limit theorem of $Z_t^\alpha$ using the moment generating function derived before.
\begin{proof}[proof of Theorem \ref{Z_t^alpha} $(iii)$]
Let
\begin{equation*}
    Y_t=\frac{Z^\alpha_t-\E[Z^\alpha_t]}{\sigma_t},\quad\sigma_t^2=\widetilde{D}_\alpha t^{(1+2\alpha)/2}.
\end{equation*}
It is enough to show that the moment generating function of $Y_t$ converges to $\exp(\lambda^2/2)$. The moment generating function can be derived by (\ref{MGF of Z_t^alpha}) and (\ref{expecation of Z_t^alpha}) since
\begin{align*}
    \E[\exp(\lambda Y_t)]&=\E\left[\exp\left(\frac{\lambda Z_t^\alpha}{\sigma_t}\right)\exp\left(-\frac{\lambda \E[Z_t^\alpha]}{\sigma_t}\right)\right]\\
    &=\exp\left(2^{3/2}\rho\sqrt{t}\int_0^\infty \left(\exp\left(\frac{\lambda(\sqrt{2t}y)^\alpha}{\sigma_t}\right)-1\right)\erfc(y)\dif y\right)\\
    &\times\exp\left(-\frac{\lambda}{\sigma_t}D_\alpha t^{(1+\alpha)/2}\right)=\exp\left(L_t(\lambda)\right),
\end{align*}
where
\begin{equation*}
    L_t(\lambda)=2^{3/2}\rho\sqrt{t}\int_0^\infty\left(\exp\left(\frac{\lambda(\sqrt{2}y)^\alpha}{\sqrt{\widetilde D_\alpha}t^{1/4}}\right)-1-\frac{\lambda}{\sigma_t}(\sqrt{2t}y)^\alpha\right)\erfc(y)\dif y.
\end{equation*}
Therefore it is enough to show $L_t(\lambda)\to\lambda^2/2$ as $t\to \infty$. Write
\begin{equation*}
    r_t=\frac{\lambda 2^{\alpha/2}}{\sqrt{\widetilde  D_\alpha}t^{1/4}}
\end{equation*}
Then
\begin{equation*}
\begin{split}
    L_t(\lambda)&=2^{3/2}\rho\sqrt{t}\int_0^\infty \left(\exp(r_ty^\alpha)-1-r_t y^\alpha\right)\erfc(y)\dif y\\
    &=2^{3/2}\rho\sqrt{t}\int_0^\infty \sum_{k=2}^\infty \frac{(r_ty^\alpha)^k}{k!}\erfc(y)\dif y\\
    &=\rho\sqrt{2t}r_t^2\int_0^\infty y^{2\alpha}\erfc(y)\dif y+\rho 2^{3/2}V_t,
    \end{split}
\end{equation*}
where
\begin{equation*}
    V_t=\sum_{k=3}^\infty \sqrt{t}r_t^k\int_0^\infty \frac{y^{\alpha k}}{k!}\erfc(y)\dif y.
\end{equation*}
Note that $r_t\to0$ as $t\to\infty$ so using $\erfc(x)\leq \exp(-x^2)$ we get for large enough $t$ such that $r_t<1$, 
\begin{equation*}
\begin{split}
    V_t&\leq\sum_{k=3}^\infty \sqrt{t}r_t^k\sum_{k=0}^\infty \int_0^\infty \frac{y^{\alpha k}}{k!}\erfc(y)\dif y\\
    &=\sum_{k=3}^\infty \sqrt{t}r_t^k\int_0^\infty \exp(y^\alpha)\erfc(y)\dif y\\
    &\leq\sum_{k=3}^\infty \sqrt{t}r_t^k\int_0^\infty \exp(y^\alpha-y^2)\dif y\\
    &\leq C\sum_{k=3}^\infty \sqrt{t}r_t^k=\frac{C\sqrt{t}r_t^3}{1-r_t}
\end{split}
\end{equation*}
Hence 
\begin{equation*}
    L_t=\rho\sqrt{2t}r_t^2\int_0^\infty y^{2\alpha}\erfc(y)\dif y+O\left(\sqrt{t}r_t^3\right)
\end{equation*}
By the definition of $\widetilde  D_\alpha$ given by (\ref{def of D_alpha}), the first term is equal to $\lambda^2/2$. The second vanishes as $r_t^3\sqrt{t}=O(t^{-1/4})$. Hence $L_t(\lambda)\to\lambda^2/2$ as desired.
\end{proof}

\section{Position of top particles}\label{proof of theorem 4}
\subsection{H\"older continuity, Theorem \ref{theorem of M_t} (i)}
Let $M_t$ be the position of top particle at time $t$ defined by (\ref{def of M_t}). According to Remark \ref{no cluster point except zero}, $(M_t)_{t\geq 0}$ is a continuous process.
Indeed, for each $\epsilon>0$, the number of particles whose positions are at least $\epsilon$ away from zero is almost surely finite, and each particle path is continuous.

\begin{remark}\label{intuition for M_t}
The intuition for law of large numbers and central limit theorem for $M_t$ is that its asymptotic behavior is governed by particles whose trajectory are excursions with maximum positions away from the origin. Since the maximum position tends to infinity as $t\to\infty$, particles associated with small excursions are expected to have a negligible contribution. Nevertheless, to make our proofs rigorous, we retain these particles and control their contribution using an argument similar to the previous section. The main part of the proof therefore focuses on particles with large excursions, which determine the asymptotic behavior of $M_t.$    
\end{remark}

Recall the old-new decomposition of $X_{s+t}$ given by (\ref{decomposition of X_s+t}). Let $M_{s+t}^{\text{old}}$ be the position of top particle for the point process $X_{t+s}^{\text{old}}$, that is,
\begin{equation*}
    M_{s+t}^{\text{old}}:=\max_{1\leq i\leq N_s}x_i(s+t)\vee 0.
\end{equation*}
Similarly define 
\begin{equation*}
    M_{s+t}^{\text{new}}=\max\{e_i(s+t-t_i):(t_i,e_i)\in\Theta,\ s\leq t_i\leq s+t\}\vee 0
\end{equation*}
to be the position of the top particle born during the interval $[s,s+t]$. Then for all $s,t\in\R^+$
\begin{equation*}
    M_{s+t}=M_{s+t}^{\text{old}}\vee M_{s+t}^{\text{new}}.
\end{equation*}
To establish the H\"older continuity of $M_t$, it suffices to obtain suitable moment estimates for the increments $M_{t+s}-M_t$. As discussed in Remark \ref{intuition for M_t}, we mainly focus on $M_{t+s}^{\mathrm{old}}-M_t$, since particles born during a short time interval are unlikely to make large excursions away from the origin. Note that
\begin{equation}\label{upper bound of M_s+t-M_s}
    \begin{split}
        |M_{s+t}-M_s|^k&=|M_{s+t}^{\text{old}}\vee M_{s+t}^{\text{new}}-M_s|^k\\
        &\leq\left(|M_{s+t}^{\text{old}}-M_s|+M_{s+t}^{\text{new}}\right)^k\\
        &\leq 2^{k-1}\left(|M_{s+t}^{\text{old}}-M_s|^k+(M_{s+t}^{\text{new}})^k\right).
    \end{split}
\end{equation}
Therefore, it is enough to bound moments of $M_{s+t}^{\text{new}}$ and $|M_{s+t}^\text{old}-M_s|$ respectively and apply Kolmogorov continuity theorem.By Remark \ref{distribution of X_t+}, the restriction of $X_{s+t}^{\text{new}}$ to $\R^+$ has the same distribution of $X_t^+$ which is a Poisson point process with intensity measure $\rho\erfc(x/\sqrt{2t})\dif x$ on $\R^+$. Therefore for any $r>0$,
\begin{equation*}
    \pr[M_{s+t}^{\text{new}}\leq r]=\pr[X_{t}^{+}[r,\infty)=0]=\exp\left(-\rho\int_{r}^\infty \erfc(\frac{x}{\sqrt{2t}})\dif x\right).
\end{equation*}
\begin{remark}\label{distribution of M_t}
    Note that if we take $s=0$ then $M^{\text{new}}_t=M_t$, therefore we have the distribution of $M_t$ given by
    \begin{equation*}
        \pr[M_{t}\leq r]=\exp\left(-\rho\sqrt{2t}\int_{r/\sqrt{2t}}^\infty \erfc(x)\dif x\right).
    \end{equation*}
\end{remark}
\begin{lemma}\label{upper bound of M new}
    For any $t>0$ and $k\geq 1$, we have
    \begin{equation*}
        \E[(M_{t}^{\mathrm{new}})^k]=\E[M_t^k]\leq C_{\rho,k}t^{k/2}\left(\log(e+t)\right)^{k/2}.
    \end{equation*}
\end{lemma}
\begin{proof}
    Since $M_t$ has the same distribution as $M_t^{\mathrm{new}}$ it suffices to prove the bound for $\E[M_t^k]$. The $k$-th moment of $M_t$ can be computed by its distribution function:
    \begin{equation*}
        \E[M_t^k]=\int_0^\infty kr^{k-1}\pr[M_t>r]\dif r.
    \end{equation*}
    We split this integral at $R_t=\sqrt{4t\log(e+t)}$. For the first part we have
    \begin{equation*}
        \int_0^{R_t} kr^{k-1}\pr[M_t>r]\dif r\leq\int_0^{R_t}k r^{k-1}\dif r= R^k_t\leq C_k t^{k/2}\left(\log(e+t)\right)^{k/2}.
    \end{equation*}
    By Remark \ref{distribution of M_t} for any $r$ we obtain, using $\erfc(x)\leq \exp(-x^2)$, that 
    \begin{equation*}
        \begin{split}
            \pr[M_t>r]&=1-\exp\left(-\rho\sqrt{2t}\int_{r/\sqrt{2t}}^\infty \erfc(x)\dif x\right)\\
            &\leq \rho\sqrt{2t}\int_{r/\sqrt{2t}}^\infty \erfc(x)\dif x\\
            &\leq C_{\rho}\sqrt{t}\int_{r/\sqrt{2t}}^\infty \exp(-x^2)\dif x\\
            &\leq C_{\rho}\sqrt{t}\int_{r/\sqrt{2t}}^\infty \frac{x}{r/\sqrt{2t}}\exp(-x^2)\dif x\\
            &=C_{\rho}\frac{t}{r}\exp\left(-\frac{r^2}{2t}\right).
        \end{split}
    \end{equation*}
    Therefore, for the second part we get
    \begin{equation*}
        \begin{split}
            \int_{R_t}^\infty kr^{k-1}\pr[M_t>r]\dif r&\leq C_{\rho,k}t\int_{R_t}^\infty r^{k-2}\exp\left(-\frac{r^2}{2t}\right)\dif r\\
            &=C_{\rho,k}t^{(k+1)/2}\int_{R_{t}/\sqrt{2t}}^\infty u^{k-2}\exp(-u^2)\dif u.
        \end{split}
    \end{equation*}
    Note that there exists a constant $C_k$ such that $u^{k-2}\leq C_k\exp(u^2/2)$ for all $u\in \R^+$. Hence the above integral is bounded by
    \begin{equation*}
        \begin{split}
            C_{\rho,k} t^{(k+1)/2}\int_{R_t/\sqrt{2t}}^\infty \exp\left(-\frac{u^2}{2}\right)\dif u&\leq C_{\rho,k} t^{(k+1)/2}\exp\left(-\frac{R_t^2}{8t}\right)\\
            &=C_{\rho, k}t^{(k+1)/2}\exp\left(-\frac{1}{2}\log(e+t)\right)\\
            &\leq C_{\rho,k}t^{k/2}.
        \end{split}
    \end{equation*}
    Combining the two estimates gives
    \begin{equation*}
        \E[M_t^k]\leq C_{\rho,k}t^{k/2}\left(\log(e+t)\right)^{k/2}.
    \end{equation*}
This proves the lemma.
\end{proof}
We next estimate the contribution from the old particles. The following bound for moments will also be used later in the proof of the law of large numbers for $M_t$.
\begin{lemma}\label{upper bound of M old}
    There exists $S_0>0$ and $C_k>0$ such that
    \begin{equation*}
        \E\left[\left(\sup_{t\in[0,T]}|M_{s+t}^{\mathrm{old}}-M_s|\right)^k\right]\leq \begin{cases}
            C_k\sqrt{s}T^{k/2}\quad&\text{for }s\leq S_0,\\
            C_k(T\log s)^{k/2}\quad&\text{for }s>S_0,
        \end{cases}
    \end{equation*}
    for all $s,T>0$.
\end{lemma}
\begin{proof}
    Let $(B^{i}_t)_{t\geq 0}$ be a finite collection of independent Brownian motions starting from zero. Let $x_i$ be a finite collection of positive numbers. Define for each $i$ the stopping time $\tau_i=\inf_{u\in[0,T]}\{u:B_u^i+x_i=0\}$ and
    \begin{equation*}
    Y_t^i=
        \begin{cases}
          B_t^i+x_i\quad &{t<\tau_i},\\
          0,\quad &t\geq \tau_i.
        \end{cases}
    \end{equation*}
   Then $Y^i_t$ is a finite collection of Brownian motions starting from $x_i$ absorbing at zero. For any $t\in[0,T]$ we have $|Y_t^i-x_i|=|B_t^i|\ind_{\{t<\tau_i\}}+x_i\ind_{\{t\geq \tau_i\}}$. Note that
    \begin{equation*}
        |B_t^i|\leq\sup_{u\in[0,T]}|B_u^i|,\quad x_i=|B_{\tau_i}^i|\leq\sup_{u\in[0,T]}|B_u^i|.
    \end{equation*}
    Hence we get
    \begin{equation*}
        \sup_{t\in[0,T]}|Y_t^i-x_i|\leq\sup_{u\in[0,T]}|B_u^i|.
    \end{equation*}
    Therefore
    \begin{equation*}
        \begin{split}
            \sup_{t\in[0,T]}|\max_{i}Y_t^i-\max_{i}x_i|&\leq\sup_{t\in[0,T]}\max_i|Y_t^i-x_i|\\
            &\leq\max_{i}\sup_{t\in[0,T]}|Y_t^i-x_i|\\
            &\leq\max_i\sup_{t\in[0,T]}|B_t^i|.
        \end{split}
    \end{equation*}
    Let $\nu_T(\dif y)$ denote the distribution of $\sup_{t\in[0,T]}|B_t^1|$. At time $s$, by Remark \ref{distribution of X_t+}, $X^+_s$ is a Poisson point process on $\R^+$ with intensity $\rho\erfc(x/\sqrt{2s})\dif x$. Assign to each point in $X^+_s$ an independent mark $U^{(i)}_{s,T}$ with distribution $\nu_T(\dif y)$. By the marking theorem, the marked point process is a Poisson point process on $\R^+\times \R^+$ with intensity measure $\rho\erfc(x/\sqrt{2s})\nu_T(\dif y)$. We can interpolate this Poisson point process as follows: at time $s$, the positions of alive particles are distributed as a Poisson point process with intensity $\rho\erfc(x/\sqrt{2s})\dif x$. For each particle alive at time $s$, its path over the time interval $[s,s+T]$ is dominated by its corresponding mark $U_{s,T}^{(i)}$. Hence 
    \begin{equation}\label{upper bound by U_s,T}
        \sup_{t\in[0,T]}|M_{s+t}^\text{old}-M_s|\leq\max_{i} U^{(i)}_{s,T}=:U_{s,T}.
    \end{equation}
    Note that the marks $U^{(i)}_{s,T}$ forms a Poisson point process with intensity 
    \begin{equation*}
        \int_0^\infty \rho\erfc(x/\sqrt{2s})\dif x\nu_T(\dif y)=C\sqrt{s}\nu_T(\dif y).
    \end{equation*}
   Consequently for $r>0$ we have
   \begin{equation*}
        \begin{split}
       \pr[U_{s,T}>r]&=1-\exp\left(-\int_r^\infty C\sqrt{s}\nu_T(\dif y)\right)\\&\leq C\sqrt{s}\pr\left[\sup_{u\in[0,T]}|B_u|>r\right]\\&\leq C\sqrt{s}\exp\left(-\frac{r^2}{2T}\right),
       \end{split}
   \end{equation*}
where we use reflection principle in the last inequality. We use the tail formula for moments given by
\begin{equation}\label{expecation of U_s,T^k}
    \E[U_{s,T}^k]=\int_0^\infty kr^{k-1}\pr[U_{s,T}>r]\dif r.
\end{equation}
Splitting the above integral at $R_{s,T}=\sqrt{2T\log (C\sqrt{s})}$ we obtain
\begin{align*}
    \E[U_{s,T}^k]&\leq \int_0^{R_{s,T}}k r^{k-1}\dif r+\int_{R_{s,T}}^\infty kr^{k-1}C\sqrt{s}\exp\left(-\frac{r^2}{2T}\right)\dif r\\
    &=R_{s,T}^k+Ck\sqrt{s}\sqrt{2T}^k\int_{R_{s,T}/\sqrt{2T}}^\infty u^{k-1}\exp(-u^2)\dif u.
\end{align*}
Note that $R_{s,T}/\sqrt{2T}=\sqrt{\log(C\sqrt{s})}\to\infty$ as $s\to\infty$, so we have for any large enough $s$, the second integral:
\begin{equation*}
    \begin{split}Ck\sqrt{s}\sqrt{2T}^k\int_{R_{s,T}/\sqrt{2T}}^\infty u^{k-1}\exp(-u^2)\dif u&\approx Ck\sqrt{s}\sqrt{2T}^k\left(\frac{R_{s,T}}{\sqrt{2T}}\right)^{k-2}\exp\left(-\log(C\sqrt{s})\right)\\&\leq C_k TR_{s,T}^{k-2}\\&\leq C_kR_{s,T}^k
    \end{split}
\end{equation*}
 Therefore, there exists $S_0$ such that for all $s>S_0$,
 \begin{equation*}
     \E[U_{s,T}^k]\leq C_k\left(T\log(s)\right)^{k/2}
 \end{equation*}
For $s\leq S_0$ by (\ref{expecation of U_s,T^k}) we get
\begin{align*}
    \E[U^k_{s,T}]&\leq \int_0^\infty k r^{k-1}C\sqrt{s}\exp\left(-\frac{r^2}{2T}\right)\dif r\\&=Ck\sqrt{s}\sqrt{2T}^k\int_0^\infty u^{k-1}\exp(-u^2)\dif u\\&=C_kT^{k/2}\sqrt{s}.
\end{align*}
Combining this with (\ref{upper bound by U_s,T}) proves the lemma.
\end{proof}

We now prove Theorem \ref{theorem of M_t} $(i)$, which is the H\"older continuity of $M_t$.
\begin{proof}[proof of Theorem \ref{theorem of M_t} $(i)$]
It is enough to show that for sufficiently large $T>0$ and $t\in[0,T]$, 
\begin{equation}\label{upper bound in proving holder continuous}
    \sup_{s\in[0,T]}\E[|M_{s+t}-M_s|^k]\leq C_k\left(\log (e+T)\right)^{k/2}t^{k/2}.
\end{equation}
Then by Kolmogorov continuity theorem we get $(M_s)_{s\in[0,T]}$ is $\gamma-$H\"older continuous on $[0,T]$ for any $\gamma<1/2$. Note that by Lemma \ref{upper bound of M old} we have
\begin{equation*}
    \begin{split}
        \sup_{s\in[0,T]}\E[|M^{\text{old}}_{s+t}-M_s|^k]&\leq\sup_{s\in[0,T]}\E\left[\left(\sup_{v\in[0,t]}|M^{\text{old}}_{s+v}-M_{s}|\right)^k\right]\\
        &\leq\sup_{s\in[0,T]}C_k\left(\sqrt{s}t^{k/2}\ind_{\{s<S_0\}}+(t\log s)^{k/2}\ind_{\{s\geq S_0\}}\right)\\
        &\leq C_k\left(\sqrt{S_0}t^{k/2}\ind_{\{s<S_0\}}+(t\log T)^{k/2}\ind_{\{s\geq S_0\}}\right).
    \end{split}
\end{equation*}
Since $S_0$ does not depends on $T$ we can choose $T$ being large enough such that $(\log T)^{k/2}\geq\sqrt{S_0}$. Hence
\begin{equation*}
     \sup_{s\in[0,T]}\E[|M^{\text{old}}_{s+t}-M_s|^k]\leq C_k t^{k/2}(\log T)^{k/2}
\end{equation*}
 Combining (\ref{upper bound of M_s+t-M_s}) and Lemma \ref{upper bound of M new} we get
\begin{equation*}
    \begin{split}
        \sup_{s\in[0,T]}\E[|M_{s+t}-M_s|^k]&\leq C_k\left(t^{k/2}(\log T)^{k/2}+t^{k/2}\left(\log(e+t)\right)^{k/2}\right).
    \end{split}
\end{equation*}
Since for any $t\in[0,T]$ we have $\log (e+t)\leq\log (e+T)$. Hence the above upper bound is less equal than $C_k t^{k/2}\left(\log (e+T)\right)^{k/2}$, which is desired upper bound.
\end{proof}

\subsection{Law of large numbers, Theorem \ref{theorem of M_t} (ii)}
Let
\begin{equation*}
    a_n=\exp\left((\log n)^2\right),\quad\mu_t=\sqrt{t\log t}.
\end{equation*}
To prove the law of large numbers for $M_t$ we first establish convergence along the sequence $\left(a_n\right)_{n\geq 1}$.
\begin{lemma}\label{LLN along mu_an}
    Let $a_n$ and $\mu_t$ be defined as before then
    \begin{equation*}
        \frac{M_{a_n}}{\mu_{a_n}}\to 1
    \end{equation*}
    almost surely.
\end{lemma}
\begin{proof}
    Recall that by Remark \ref{distribution of M_t} we have the distribution function of $M_t$. Therefore, on the one hand, we have
    \begin{align*}
        \pr[M_{a_n}>(1+\epsilon)\mu_{a_n}]&=1-\exp\left(-\rho\sqrt{2a_n}\int_{(1+\epsilon)\sqrt{\log (a_n) /2}}^\infty\erfc(x)\dif x\right)\\
        &\approx1-\exp\left(-\frac{\rho C\sqrt{a_n}\exp(-(1+\epsilon)^2\log (a_n)/2)}{(1+\epsilon)^2\log a_n}\right)\\
        &=1-\exp\left(-\frac{\rho C_\epsilon}{\log (a_n) a_n^{\epsilon+\epsilon^2/2}}\right)\\
        &\approx\frac{\rho C_\epsilon}{(\log n)^2 n^{(\epsilon+\epsilon^2/2)\log n}},
    \end{align*}
    which is summable. Hence by Borel-Cantelli lemma we get $\limsup_{n\to\infty}M_{a_n}/\mu_{a_n}\leq 1$ almost surely.
    On the other hand, we have
    \begin{equation*}
        \begin{split}
        \pr[M_{a_n}\leq(1-\epsilon)\mu_{a_n}]&=\exp\left(-\rho\sqrt{2a_n}\int_{(1-\epsilon)\sqrt{\log (a_n)/2}}^\infty \erfc(x)\dif x\right)\\
        &\approx\exp\left(-\frac{\rho C\sqrt{a_n}\exp(-(1-\epsilon)^2\log (a_n)/2)}{(1-\epsilon)^2\log a_n}\right)\\
        &=\exp\left(-\frac{\rho C_\epsilon a_n^{\epsilon-\epsilon^2/2}}{(\log n)^2}\right)\\
        &=\exp\left(-\frac{\rho C_\epsilon n^{(\epsilon-\epsilon^2/2) \log n}}{(\log n)^2}\right),
        \end{split}
    \end{equation*}
    which is summable. Hence $\liminf_{n\to\infty }M_{a_n}\geq 1$. Therefore get $M_{a_n}/\mu_{a_n}\to 1$ almost surely.
\end{proof}
Let $F(t)=\max\{a_n:a_n\leq t\}$, the next step is to show $M_t/\mu_{F(t)}$ converges to one almost surely. To do this we first prove the following maximal moment estimate.
\begin{lemma}\label{upper bound of sup M_t}
    For sufficiently large $T$ and $p$ we have
    \begin{equation*}
        \E\left[\left(\sup_{t\in[0,T]}|M_t|\right)^p\right]\leq C_p\left(\log (e+T)\right)^{p/2}T^{p/2}.
    \end{equation*}
\end{lemma}
\begin{proof}
    Let $\gamma\in(0,1/2)$, then $M_t$ is $\gamma$-H\"older continuous on $[0,T]$ by Theorem \ref{theorem of M_t} $(i)$. Since $M_0=0$ we have $\sup_{t\in[0,T]}|M_t|\leq K_{\gamma,T}T^\gamma$, where $K_{\gamma,T}$ is the H\"older constant. Therefore,
    \begin{equation}\label{upper bound for p moment of sup M_t}
        \E\left[\left(\sup_{t\in[0,T]} M_t\right)^p\right]\leq T^{p\gamma}\E[K_{\gamma,T}^p].
    \end{equation}
    It remains to estimate moments of the H\"older constant. We shall do this by using Lemma \ref{GRR} to get a constant $C_{\gamma,p}$ such that
    \begin{align*}
        \E[K_{\gamma,T}^p]&\leq C_{\gamma,p} \int_0^T\int_0^T \frac{\E[|M_u-M_v|^p]}{|u-v|^{p\gamma+2}}\dif u \dif v\\
        &=C_{\gamma,p}\int_0^T\int_0^{T-t}\frac{\E[|M_{u+t}-M_u|^p]}{t^{p\gamma+2}}\dif u\dif t\\
        &\leq C_{\gamma,p}\int_0^T\int_0^T \frac{\sup_{u\in[0,T]}\E[|M_{u+t}-M_u|^p]}{t^{p\gamma +2}}\dif u\dif t\\
        &\leq C_{\gamma,p}\int_0^T\int_0^T \left(\log (e+T)\right)^{p/2}t^{p(1/2-\gamma)-2}\dif u \dif t,
    \end{align*}
    where in the last inequality we use (\ref{upper bound in proving holder continuous}) in the proof of Theorem \ref{theorem of M_t} $(i)$. Now we can choose $p$ large enough so to get
    \begin{equation*}
        \begin{split}
            \E[K_{\gamma,T}^p]\leq C_{\gamma,p}\left(\log (e+T)\right)^{p/2}T^{p(1/2-\gamma)}.
        \end{split}
    \end{equation*}
    Hence by (\ref{upper bound for p moment of sup M_t}) we get the desired upper bound.
\end{proof}
\begin{lemma}\label{LLN along mu_F(t)}
    Let $F(t)$ and $\mu_t$ be defined as before, then
    \begin{equation*}
        \left|\frac{M_t-M_{F(t)}}{\mu_{F(t)}}\right|\to0
    \end{equation*}
    almost surely.
\end{lemma}
\begin{proof}
    Let $\Delta_n=a_{n+1}-a_n$, it is enough to prove that for any $\epsilon>0$ the probability 
    \begin{equation*}
        \pr\left[\sup_{t\in[0,\Delta_n]}\left|\frac{M_{a_n+t}-M_{a_n}}{\mu_{a_n}}\right|>\epsilon\right]
    \end{equation*}
    is summable in $n$. Then Borel-Cantelli lemma implies the result.
    According to (\ref{upper bound of M_s+t-M_s}) we have
    \begin{equation*}
        \sup_{t\in[0,\Delta_n]}\left|M_{a_n+t}-M_{a_n}\right|\leq\sup_{t\in[0,\Delta_n]}|M_{a_n+t}^\text{old}-M_{a_n}|+\sup_{t\in[0,\Delta_n]}|M^{\text{new}}_{a_n+t}|.
    \end{equation*}
    Hence
    \begin{align*}
        \pr\left[\sup_{t\in[0,\Delta_n]}|M_{a_n+t}-M_{a_n}|>\epsilon \mu_{a_n}\right]&\leq \pr\left[\sup_{t\in[0,\Delta_n]}|M_{a_n+t}^\text{old}-M_{a_n}|>\frac{\epsilon \mu_{a_n}}{2}\right]\\
        &+\pr\left[\sup_{t\in\Delta_n}|M^{\text{new}}_{a_n+t}|>\frac{\epsilon\mu_{a_n}}{2}\right]:=A_{n,\epsilon}+B_{n,\epsilon}.
    \end{align*}
    For $A_{n,\epsilon}$ by Markov's inequality and Lemma \ref{upper bound of M old} we get for large enough $n$:
    \begin{equation*}
        A_{n,\epsilon}\leq \frac{\E[(\sup_{t\in[0,\Delta_n]}|M^{\mathrm{old}}_{a_n+t}-M_{a_n}|)^4]}{(\epsilon/2)^4\mu_{a_n}^4}\leq C_{\epsilon}\left(\frac{\Delta_n}{a_n}\right)^{2}.
    \end{equation*}
    Note that
    \begin{equation*}
    \begin{split}
        \frac{\Delta_n}{a_n}&=\exp\left(\left(\log(1+\frac{1}{n})\right)^2+2\log(1+\frac{1}{n})\log n\right)-1\\&\approx\exp\left(\frac{1}{n^2}+\frac{2\log n}{n}\right)-1\\&\approx\frac{\log n}{n}.
        \end{split}
    \end{equation*}
    Therefore, $A_{n,\epsilon}$ is summable in $n$. As for $B_{n,\epsilon}$, note that $\left(M^{\text{new}}_{a_n+t}\right)_{t\in[0,\Delta_n]}$ has the same distribution as $\left(M_t\right)_{{t\in[0,\Delta_n]}}$, hence by Lemma \ref{upper bound of sup M_t} for a large enough $p$ we get:
    \begin{equation*}
        B_{n,\epsilon}\leq\frac{\E[(\sup_{t\in[0,\Delta_n]}M_t^{\text{new}})^p]}{(\epsilon/2)^p\mu_{a_n}^p}\leq C_{\epsilon}\left(\frac{\Delta_n \log(e+\Delta_n)}{a_n\log a_n}\right)^{p/2}\approx\left(\frac{\Delta_n \log\Delta_n}{a_n\log a_n}\right)^{p/2}.
    \end{equation*}
    Note that $\Delta_n/a_n\approx (\log n)/n$, so we have 
    \begin{equation*}
        \frac{\log\Delta_n}{\log a_n}\approx\frac{\log a_n+\log(\log (n) /n)}{\log a_n}=1+\frac{\log (\log(n)/n)}{(\log n)^2}\to 1
    \end{equation*}
    as $n\to \infty$. therefore $B_{n,\epsilon}$ is summable. We get the desired result.
\end{proof}
\begin{proof}[proof of Theorem \ref{theorem of M_t} $(ii)$]
    Note that
    \begin{equation*}
        \left|\frac{M_t}{\mu_t}-\frac{M_t}{\mu_{F(t)}}\right|=\frac{M_t}{\mu_{F(t)}}\left|1-\frac{\mu_{F(t)}}{\mu_t}\right|.
    \end{equation*}
    By Lemma \ref{LLN along mu_an} and \ref{LLN along mu_F(t)} we know that $M_t/\mu_{F(t)}\to 1$ almost surely. Hence it is enough to show that $|1-\mu_{F(t)}/\mu_t|\to 0$. On the one hand it is easy to see that $\mu_{F(t)}/\mu_t\leq 1$ as $\mu_t$ is increasing in $t$. On the other hand, let $a_n=F(t)$ then according to the definition of $a_n$ we can check that
    \begin{equation*}
        \frac{\mu_{F(t)}}{\mu_t}\geq\frac{\mu_{a_n}}{\mu_{a_{n+1}}}=\sqrt{\frac{a_n(\log n)^2}{a_{n+1}(\log(n+1))^2}}\to 1
    \end{equation*}
    as $n\to \infty$. Therefore almost surely we have $|1-\mu_{F(t)}/\mu_t|\to0$ as desired.
\end{proof}

\subsection{Central limit theorem, Theorem \ref{theorem of M_t} (iii)}
Recall that $A(y)=\int_y^\infty \erfc(x)\dif x$, $u_t$ is the unique solution of $\sqrt{2t}A(y)=1$ and $v_t=\sqrt{t}/(\sqrt{2}u_t)$. Note that
\begin{equation*}
    u_t\approx \sqrt{\frac{\log t}{2}}\quad\text{as }t\to\infty.
\end{equation*}
Therefore $\sqrt{2t} u_t$ is asymptotically equal to $\E[M_t]$ and
\begin{equation*}
    v_t \approx \sqrt{\frac{t}{\log t}}\quad\text{as }t\to\infty,
\end{equation*}
which is asymptotically equal to the standard deviation of $M_t$.
By Remark \ref{distribution of M_t} we get
\begin{equation}\label{distribution of M_t written by A}
    \pr[M_t\leq r]=\exp\left(-\rho \sqrt{2t}A(\frac{r}{\sqrt{2t}})\right).
\end{equation}
We shall check that the distribution function of $(M_t-\sqrt{2t}u_t)/v_t$ converges to the distribution function of a Gumbel distribution.
\begin{proof}[proof of Theorem \ref{theorem of M_t} $(iii)$]
    By (\ref{distribution of M_t written by A}) we have
    \begin{equation*}
        \pr[M_t\leq v_tr+\sqrt{2t}u_t]=\exp\left(-\rho\sqrt{2t}A(\frac{r}{2 u_t}+u_t)\right).
    \end{equation*}
    Using integration by part we can easily see that 
    \begin{equation*}
        A(y)=\frac{\exp(-y^2)}{2\sqrt{\pi}y^2}\left(1+O(y^{-2})\right),\quad{\text{as }y\to\infty.}
    \end{equation*}
    By this asymptotic expansion we get
    \begin{equation*}
    \sqrt{2t}A\left(\frac{r}{2u_t}+u_t\right)=\sqrt{2t}A(u_t)\frac{A\left({r}/{2u_t}+u_t\right)}{A(u_t)}=\frac{A\left({r}/{2u_t}+u_t\right)}{A(u_t)}\to\exp(-r).
    \end{equation*}
    Therefore we get
    \begin{equation*}
        \pr\left[\frac{M_t-\sqrt{2t}u_t}{v_t}\leq r\right]\to\exp(-\rho\exp(-r))=\exp(-\exp(-r+\log \rho)),
    \end{equation*}
which is the distribution function of shifted Gumbel distribution.
\end{proof}

\section{Convergence of $\lambda\n^\lambda$}\label{convergence of n^lambda}
The aim of this section is to prove Theorem \ref{M0 convergence of lambda n lambda}. We start by introducing some notations. For $T\in\R^+$ and $k\in\N$, define the restriction map $R_T: C_0(\R^+)\to C[0,T]$ by
\begin{equation*}
    R_T(e)=e|_{[0,T]}.
\end{equation*}
And define 
\begin{equation*}
    \mathcal{E}_T:=\left\{R_T(e):\|R_T(e)\|>0,\ e\in C_0(\R^+)\right\}\quad\text{and}\quad\mathcal{E}_{k,T}:=\left\{e^{T}\in\mathcal{E}_T:\|e^T\|\geq\frac{1}{k}\right\},
\end{equation*}
where $\|\cdot\|$ is the supremum norm. Recall that $\n^\lambda$ is the distribution of $(B^\lambda_\cdot)$, we define $\n^{\lambda,T}$ to be the distribution of $\left(R_T (B^\lambda_\cdot)\right)$. Hence $\n^{\lambda,T}=\n^\lambda\circ R^{-1}_T$ is the pushforward measure of $\n^\lambda$ under $R_T$. Define
\begin{equation*}
    \mu_\lambda^k:=\lambda \n^{\lambda,T}|_{\mathcal E_{k,T}}
\end{equation*}
be a collection of measures supported on $\mathcal{E}_{k,T}$.
\begin{lemma}
    For any $k\in\N$, there exists a sequence of compact set $(\mathcal{K}_M)_{M\geq1}$ such that
    \begin{equation*}
        \mu_\lambda^k(\mathcal{E}_{k,T}\setminus \mathcal K_M)\to0,\quad\text{as }\lambda\to\infty.
        \end{equation*}
    That is, $\left(\mu^k_{\lambda}\right)_{\lambda>0}$ is a collection of tight measures.
\end{lemma}
\begin{proof}
    For $M\in\N$, we set
    \begin{equation*}
        \mathcal{K}_M=\left\{e^T\in \mathcal{E}_{k,T}: \|e^T\|\leq M,\sup_{{s,t\in[0,T]};|t-s|\leq \delta_N}|e(t)-e(s)|\leq \frac{1}{N},\ \text{for all }N\geq M\right\},
\end{equation*}
where $\delta_N$ is a sequence of positive numbers to be choosen later. By the Arzel\`a–Ascoli theorem, $\mathcal{K}_M$ is compact for every $M\in\N$. Hence for $M\geq 1$,
 \begin{equation}\label{check tightness}
        \begin{split}
            &\mu_\lambda^k\left(\mathcal{E}_{k,T}\setminus \mathcal K_M\right)\leq \mu^k_\lambda\left(\left\{e^T\in\mathcal{E}_{T}:\|e^T\|> M\right\}\right)\\
            &+\mu^k_\lambda\left(\left\{e^T\in\mathcal{E}_{k,T}:\sup_{{s,t\in[0,T]};|t-s|\leq \delta_N}|e^T(t)-e^T(s)|> \frac{1}{N}\ \text{for some }N\geq M\right\}\right)\\
            &\leq \lambda \n^\lambda\left(\left\{e\in\mathcal{E}:\|e\|>M\right\}\right)\\
            &+\lambda \n^\lambda\left(\left\{e\in\mathcal{E}:\sup_{|t-s|\leq \delta_N}|e(t)-e(s)|>\frac{1}{N}\ \text{for some }N\geq M\right\}\right)\\
            &=:I_{1,\lambda}(M)+I_{2,\lambda}(M).
        \end{split}
    \end{equation}
     To bound $I_1^\lambda(M)$, let $B_t^\lambda$ be a Brownian motion killed when its local time exceeds an independent exponential clock $\mathrm{Exp}(\lambda)$. Then we have
     \begin{equation*}
         \n^\lambda\left(\left\{e\in\mathcal{E}:\|e\|>M\right\}\right)=\pr\left[\sup_{t\in[0,\xi(B_t^\lambda)]}|B_t^\lambda|>M\right].
     \end{equation*}
     Recall from It\^o's Brownian excursion theory that Brownian motion can be reconstructed from a Poisson point process $\Gamma$ on $\R^+\times \mathcal{E}$ with intensity measure $\dif t\otimes \n$, where $\R^+$ is the local time axis of the Brownian motion. Consequently, a Brownian motion killed when its local time exceeds $\mathrm{Exp}(\lambda)$ can be reconstructed by  restricting $\Gamma$ to $[0,\mathrm{Exp}(\lambda)]\times \mathcal{E}$.
     In order to have $B_t^\lambda$ reaches level $M$, we need at least one of the excursions of $B^\lambda_t$ reaches level $M$. Hence the event $\{\sup_{t\in[0,\xi(B_t^\lambda)]}|B_t^\lambda|>M\}$ is the same as 
     \begin{equation*}
     \{\text{there exists }(t,e)\in\Gamma\text{ such that }t\leq\mathrm{Exp}(\lambda)\text{ and }\|e\|>M\}.
     \end{equation*}
     Therefore
     \begin{equation*}
         \pr\left[\sup_{t\in[0,\xi(B_t^\lambda)]}|B_t^\lambda|>M\right]=\pr\left[\Poi\left(\mathrm{Exp(\lambda)\n(\{e\in\mathcal{E}:\|e\|>M\})}\right)\geq 1\right].
     \end{equation*}
     By \cite[Proposition 2.1]{Rogers_excursion} we have
    \begin{equation*}
        \n\left(\left\{e\in\mathcal{E}:\|e\|>M\right\}\right)=\frac{1}{M}.
    \end{equation*}
    Hence we obtain
    \begin{align*}
        \pr\left[\sup_{t\in[0,\xi(B_t^\lambda)]}|B_t^\lambda|>M\right]&=\pr\left[\Poi\left(\frac{\mathrm{Exp}(\lambda)}{M}\geq1\right)\right]\\
        &=\int_0^\infty\left(1-\exp\left(-\frac{u}{M}\right)\right)\lambda\exp\left(-\lambda u\right)\dif u\\
        &=\frac{1}{1+M\lambda}.
    \end{align*}
    Therefore,
    \begin{equation*}
        \lim_{M\to\infty}\sup_{\lambda>0}I_{1,\lambda}(M)=\lim_{M\to\infty}\sup_{\lambda>0}\frac{\lambda}{1+M\lambda}\leq\lim_{M\to\infty}\frac{1}{M}=0.
    \end{equation*}
    For $I_{2,\lambda}(M)$,the union bound gives
    \begin{equation}\label{upperbound of I^2_lambda}
        I_{2,\lambda}(M)\leq \sum_{N=M}^\infty\lambda\pr\left[\sup_{|t-s|<\delta_N}|B_t^\lambda-B_s^\lambda|>\frac{1}{N}\right].
    \end{equation}
    For any $N\in\N$ and $\delta>0$, let 
    \begin{equation*}
        A_{N,\delta}=\left\{e\in\mathcal{E}:\sup_{|t-s|<\delta}|e(t)-e(s)|>\frac{1}{2N}\right\}.
    \end{equation*}
    We claim that 
    \begin{equation*}
        \left\{\sup_{|t-s|<\delta_N}|B_t^\lambda-B_s^\lambda|>\frac{1}{N}\right\}\subseteq\left\{\Gamma\left([0,\mathrm{Exp} (\lambda)]\times A_{N,\delta_N}\right)\geq 1\right\}.
    \end{equation*}
    Indeed, suppose that there exists $t_1,t_2$ with $|t_1-t_2|<\delta_N$ such that $|B_{t_1}^\lambda-B_{t_2}^\lambda|>1/N$. We consider two cases. First, suppose that $B_{t_1}^\lambda$ and $B_{t_2}^\lambda$ belong to the same excursion $e$. Then there exist $s_1$ and $s_2$ with $|s_1-s_2|<\delta_N$ such that $|e(s_1)-e(s_2)|>1/N$. Hence the claim follows. Second, suppose $B_{t_1}^\lambda$ and $B_{t_2}^\lambda$ belong to different excursions $e_1$ and $e_2$. Then there exist $s_1,s_2$ such that $B_{t_1}^\lambda=e_1(s_1)$  and $B^\lambda_{t_2}=e_2(s_2)$. Since
    \begin{equation*}
        \frac{1}{N}\leq |B^\lambda_{t_1}-B^\lambda_{t_2}|\leq |B_{t_1}^\lambda|+|B_{t_2}^\lambda|,
    \end{equation*}
it follows that either $|e_1(s_1)|>1/(2N)$ or $|e_2(s_2)|>1/(2N)$. Therefore at least one of excursions $e_1$ or $e_2$ belongs to $A_{N,\delta_N}$, the claims follows. Consequently, 
\begin{equation*}
    \pr\left[\sup_{|t-s|<\delta_N}|B^\lambda_{t}-B_{s}^\lambda|>\frac{1}{N}\right]\leq \pr\left[\left\{\Gamma\left([0,\mathrm{Exp} (\lambda)]\times A_{N,\delta_N}\right)\geq 1\right\}\right]=\frac{\n(A_{N,\delta_N})}{\lambda+\n(A_{N,\delta_N})}
\end{equation*}
Hence by (\ref{upperbound of I^2_lambda})
\begin{equation*}
    \sup_{\lambda>0}I_{2,\lambda}(M)\leq \sum_{N=M}^\infty\n(A_{N,\delta_N}).
\end{equation*}
It remains to choose a sequence of $\delta_N$ such that the right hand side converges. For fixed $N$ and $\delta>0$ observe that $A_{N,\delta}$ is a subset of $\{e\in\mathcal{E}:\|e\|>1/(4N)\}$ and $\n(\{e\in\mathcal{E}:\|e\|>(1/4N)\})=4N<\infty$. Moreover, every excursion is part of a Brownian path and is therefore H\"older continuous. Hence for fixed $N$, we have $\ind_{A_{N,\delta}}\to0$ as $\delta\to0$. By the dominated convergence theorem, $\n(A_{N,\delta})\to0$ as $\delta\to0$. Hence for each $N\in\N$ we can choose $\delta_N$ such that $\n(A_{N,\delta_N})\leq 1/N^2$. It follows that 
\begin{equation*}
    \sup_{\lambda>0}I_{2,\lambda}(M)\leq \sum_{N=M}^\infty \n(A_{N,\delta_N})\leq \sum_{N=M}^\infty\frac{1}{N^2}\to0, \quad \text{as }M\to\infty.
\end{equation*}
Combining this with (\ref{check tightness}) we get $(\mu_\lambda^k)_{\lambda> 0}$ is tight.
\end{proof}
Moreover, 
\begin{equation*}
    \sup_{\lambda>0}\mu^k_\lambda(\mathcal{E}_{k,T})=\sup_{\lambda>0}\lambda\n^{\lambda,T}(\mathcal E_{k,T})\leq \sup_{\lambda>0}I_{1,\lambda}(k^{-1})\leq k.
\end{equation*}

Therefore, by Prokhorov's theorem and and a diagonal argument, for any sequence $(\lambda_j)_{j\geq 1}$ there exists a further subsequence (again denoted by $(\lambda_j)_{j\geq 1}$ for simplicity) and a sequence of measures $(\mu^k)_{k\geq 1}$ such that $\mu_{\lambda_j}^k\to \mu^k$ weakly as $j\to \infty$. 
\begin{lemma}\label{exist a limit mu}
There exists a measure $\mu$ supported on $\mathcal{E}_T$ such that $\lambda_j\n^{\lambda_j,T}\to\mu$ under $M_0$ topology.
\end{lemma}
\begin{proof}
    We first check the consistency of $\mu^k$. Let 
    \begin{equation*}
        \mathcal{O}_{k,T}:=\{e^T\in \mathcal{E}_{k,T}:\|e^T\|>\frac{1}{k}\}
    \end{equation*}
    be the interior of $\mathcal{E}_{k,T}$. For $m\in\N$ define a function $\phi_m:\mathcal E_{k+1,T}\to \R$ by
    \begin{equation*}
        \phi_m(e)=1\wedge m\left(\|e\|-\frac{1}{k}\right)_+
    \end{equation*}
    so that $\phi_m\to 1_{\mathcal{O}_{k,T}}$ as $m\to \infty$. The function $\phi_m$ is continuous because the map $e\mapsto \|e\|$ is continuous. For any $F\in C_b(\mathcal{O}_{k,T})$, let
    \begin{equation*}
        \tilde{F}_m=\begin{cases}
            F(e)\phi_m(e)\ &\text{for }e\in\mathcal{O}_{k,T}\\
            0\ &\text{for }e\in\mathcal{E}_{k+1,T}\setminus\mathcal{O}_{k,T}
        \end{cases},\ 
        \bar{F}_m=\begin{cases}
         F(e)\phi_m(e)\ &\text{for }e\in\mathcal{O}_{k,T}\\
         0\quad&\text{for }e\in\mathcal{E}_{k,T}\setminus\mathcal{O}_{k,T}
        \end{cases}
    \end{equation*}
    be continuous extension of $F\phi_m$ to $\mathcal E_{k+1,T}$ and $\mathcal E_{k,T}$ respectively.  Since $\mu^{k+1}_{\lambda_j}|_{\mathcal{E}_{k,T}}=\mu^k_{\lambda_j}$  we have 
    \begin{equation*}
        \begin{split}
            \int_{\mathcal E_{k+1,T}}\tilde F_m(e)\mu^{k+1}_{\lambda_j}(\dif e)&=\int_{\mathcal{E}_{k,T}}\tilde F_m(e)\mu_{\lambda_j}^{k+1}(\dif e)\\
&=\int_{\mathcal{E}_{k,T}}\tilde F_m(e)\mu^{k+1}_{\lambda_j}|_{\mathcal{E}_{k,T}}(\dif e)\\
&=\int_{\mathcal{E}_{k,T}}\bar F_m(e)\mu^k_{\lambda_j}(\dif e)
        \end{split}
    \end{equation*}
    Since for any $k$ we have $\mu^k_{\lambda_j}\to \mu^k$ as $j\to\infty$, we get
    \begin{equation*}
        \int_{\mathcal{E}_{k+1,T}}\tilde F_m(e)\mu^{k+1}(\dif e)=\int_{\mathcal{E}_{k,T}}\bar F_m(e)\mu^k(\dif e).
    \end{equation*}
    Because $\phi_m\to\ind_{\mathcal{O}_{k,T}}$ as  $m\to\infty$, letting $m\to\infty$ in the preceding equality gives
    \begin{equation*}
        \int_{\mathcal{O}_{k,T}}F(e)\mu^{k+1}(\dif e)=   \int_{\mathcal{O}_{k,T}}F(e)\mu^k(\dif e),
    \end{equation*}
    which implies $\mu^{k+1}|_{\mathcal {O}_{k,T}}=\mu^k|_{\mathcal {O}_{k,T}}$. We now define a measure $\mu$ on $\mathcal E_T$ by
    \begin{equation*}
        \mu(A):=\lim_{k\to\infty}\mu^k(A\cap\mathcal O_{k,T}),\quad A\subseteq \mathcal E_T\text{ measurable}.
    \end{equation*}
    Note that the limit exist because $\left(\mu^k(A\cap \mathcal O_{k,T})\right)_{k\geq 1}$ is an increasing sequence. We next verify the consistency of $\mu$. Let $A\subseteq \mathcal O_{k,T}$. For any $m\geq k$ we have $\mathcal{O}_{k,T}\subseteq \mathcal{O}_{m,T}$ and hence
    \begin{equation*}
        \mu^m(A\cap \mathcal{O}_{m,T})=\mu^m(A\cap\mathcal{O}_{k,T})=\mu^m|_{\mathcal O_{k,T}}(A)=\mu^k(A).
    \end{equation*}
    Therefore, $\mu(A)=\lim_{m\to\infty}\mu^m(A\cap \mathcal O_{m,T})=\mu^k(A)$ so that $\mu|_{\mathcal O_{k,T}}=\mu^k|_{\mathcal{O}_{k,T}}$. Finally,we verify the convergence. For any $F\in\mathcal C_0$, we know that it vanish in a neighborhood of zero. Hence $\mathrm{supp}(F)$ is contained in $\mathcal{O}_{k,T}$ for some $k\in\N$. Thus
    \begin{align*}
        \int_{\mathcal E_T}F(e)(\lambda_j\n^{\lambda_j,T})(\dif e)&=\int_{\mathcal E_T}F(e)\ind_{\mathcal O_{k,T}}\lambda_j \n^{\lambda_j,T}(\dif e)\\
        &=\int_{\mathcal E_{k,T}}F(e)\ind_{\mathcal{O}_{k,T}}\lambda_j \n^{\lambda_j,T}|_{\mathcal{E}_{k,T}}(\dif e)\\
        &\to\int_{\mathcal{E}_{k,T}}F(e)\ind_{\mathcal{O}_{k,T}}\mu^k(\dif e)\\
        &=\int_{\mathcal{E}_T}F(e)\mu^k|_{\mathcal{O}_{k,T}}(\dif e)\\
        &=\int_{\mathcal E_T}F(e)\mu(\dif e),
    \end{align*}
    which proves the convergence under $M_0$ topology.
\end{proof}
According to (\ref{uniequness of n}), for any $t_1<\cdots<t_k<T$, we have
\begin{equation*}
    \lambda \n^\lambda\left(\{e(t_j)\in\dif x_j,\xi(e)\geq t_k\}\right)\to\n\left(\{e(t_j)\in\dif x_j,\xi(e)\geq t_k\}\right),
\end{equation*}
which provides the uniqueness for the limit of $\lambda\n^{\lambda,T}$. Hence by Lemma \ref{exist a limit mu}, for any $T>0$, $\lambda \n^{\lambda,T}\to\n\circ R_T^{-1}$. To finish the proof of Theorem \ref{M0 convergence of lambda n lambda} we define an extension map $L_T: \mathcal E_T\to \mathcal E$ by linearly interpolating between the endpoint of $e^T$ and $(T+1,0)$:
\begin{equation*}
    L_T(e^T)(s)=\begin{cases}
        e^T(s)\quad&\text{for }s\leq T,\\
        e^{T}(T)(T+1-s)\quad&\text{for }T<s\leq T+1,\\
        0\quad&\text{for }s>T+1.
    \end{cases}
\end{equation*}
We can now prove Theorem \ref{M0 convergence of lambda n lambda}.
\begin{proof}[Proof of Theorem \ref{M0 convergence of lambda n lambda}]
    For any $T>0$ and $F\in \mathcal C_0$, we have
    \begin{equation*}
        \int_{C_0(\R^+)}F(e)\n(\dif e)=\int_{C_0(\R^+)}F(e)\ind_{\{\xi(e)\leq T\}}\n(\dif e)+\int_{C_0(\R^+)}F(e)\ind_{\{\xi(e)> T\}}\n(\dif e).
    \end{equation*}
    If $\xi(e)\leq T$, then $L_TR_T(e)=e$. Hence
    \begin{equation*}
        \begin{split}
            \int_{C_0(\R^+)}F(e)&\ind_{\{\xi(e)\leq T\}}\n(\dif e)=\int_{C_0(\R^+)}F(L_TR_Te)\ind_{\{\xi(e)\leq T\}}\n(\dif e)\\
            &=\int_{C_0(\R^+)}F(L_TR_Te)\n(\dif e)-\int_{C_0(\R^+)}F(L_TR_Te)\ind_{\{\xi(e)>T\}}\n(\dif e)\\
            &=\int_{\mathcal{E}_T}F(L_T e^T)\n\circ R_T^{-1}(\dif e^T)-\int_{C_0(\R^+)}F(L_TR_Te)\ind_{\{\xi(e)>T\}}\n(\dif e).
        \end{split}
    \end{equation*}
    Therefore, we get
    \begin{equation*}
        \begin{split}
            \int_{C_0(\R^+)}F(e)\n(\dif e&=\int_{\mathcal{E}_T}F(L_T e^T)\n\circ R_T^{-1}(\dif e^T)\\
            &-\int_{C_0(\R^+)}F(L_TR_Te)\ind_{\{\xi(e)>T\}}\n(\dif e)\\&+\int_{C_0(\R^+)}F(e)\ind_{\{\xi(e)> T\}}\n(\dif e).
        \end{split}
    \end{equation*}
    Similarly, 
   \begin{equation*}
        \begin{split}
            \int_{C_0(\R^+)}F(e)\lambda \n^\lambda(\dif e)&=\int_{\mathcal{E}_T}F(L_T e^T)\lambda\n^\lambda\circ R_T^{-1}(\dif e^T)\\&-\int_{C_0(\R^+)}F(L_TR_Te)\ind_{\{\xi(e)>T\}}\lambda \n^\lambda(\dif e)\\&+\int_{C_0(\R^+)}F(e)\ind_{\{\xi(e)> T\}}\lambda \n^\lambda(\dif e).
        \end{split}
    \end{equation*}
    It follows
    \begin{equation}\label{convergence approximation}
    \begin{split}
        &\left|\int_{C_0(\R^+)}F(e) \n(\dif e)-\int_{C_0(\R^+)}F(e)\lambda \n^\lambda(\dif e)\right|\\
        &\leq \left|\int_{\mathcal{E}_T}F(L_T e^T)\n\circ R_T^{-1}(\dif e^T)-\int_{\mathcal{E}_T}F(L_T e^T)\lambda \n^\lambda\circ R_T^{-1}(\dif e^T)\right|\\
        &+2\|F\|\left(\n(\{\xi(e)> T\})+\lambda\n^\lambda(\{\xi(e)> T\})\right).
        \end{split}
    \end{equation}
   By taking $\lambda \to\infty$, the first term on the right hand side tends to zero because $\lambda \n^{\lambda ,T}\to \n\circ R_T^{-1}$ under $M_0$ topology. For the second term, the entrance law $\n_t$ of $\n$ given in (\ref{entrance law}) yields
    \begin{equation*}
        \n(\{\xi(e)> T\})=\int_\R\n(\{e(T)\in\dif x,\xi(e)>T\})=\int_\R \n_T(\dif x)=\frac{C}{\sqrt{T}}\to 0
    \end{equation*}
    as $T\to \infty$. Similarly using the entrance law of $\n^\lambda$ given in (\ref{first compute of n_t^lambda}) we obtain
    \begin{align*}
        \lambda \n^\lambda(\{\xi(e)>T\})&=\int_\R\lambda\n_T^\lambda(\dif x)\\
        &=\frac{2\lambda}{\sqrt{2\pi T^3}}\int_0^\infty\int_0^\infty (x+l)\exp\left(-\frac{(x+l)^2}{2T}-\lambda l\right)\dif l\dif x\\
        &=\frac{C}{\sqrt{T}}\int_0^\infty \lambda \exp\left(-\frac{l^2}{2T}-\lambda l\right)\dif l\\
        &\leq \frac{C}{\sqrt{T}}.
    \end{align*}
    Thus $\lim_{T\to\infty }\sup_\lambda \lambda\n^\lambda(\{\xi(e)>T\})=0$. Therefore by (\ref{convergence approximation}), we get 
    \begin{equation*}
        \int_{C_0(\R^+)}F(e)\n(\dif e)\to \int_{C_0(\R^+)}F(e)\lambda\n^\lambda(\dif e),
    \end{equation*}
   which completes the proof of $\lambda\n^\lambda$ convergence to $\n$ under the $M_0$ topology. As for the convergence of $\Theta^{\rho\lambda,\lambda}$, let $F$ be a any functional such that 
   \begin{equation*}
       F(t,e)=0\quad\text{if }t>T \text{ or }\|e\|<\epsilon.
   \end{equation*}
   Then we get the Laplace functional for $\Theta^{\rho\lambda,\lambda}$ given by
   \begin{equation*}
       \E\left[\exp(-\langle \Theta^{\rho\lambda,\lambda},F\rangle)\right]=\exp\left(-\int_{\mathcal{E}}\int_0^T(1-\exp(-F(t,e)))\dif t\lambda\n^\lambda(\dif e)\right).
   \end{equation*}
   Note that the functional 
   \begin{equation*}
       G(e):=\int_0^T \left(1-\exp(-F(t,e))\right)\dif t
   \end{equation*}
   vanish near $0$. Hence by the convergence of $\lambda\n^\lambda\to\n$ under $M_0$ topology, we obtain
   \begin{equation*}
   \begin{split}
       \E\left[\exp\left(-\langle \Theta^{\rho\lambda,\lambda},F\rangle\right)\right]&=\exp\left(-\int_{\mathcal E}G(e)\lambda\n^\lambda(\dif e)\right)\\
       &\to\exp\left(-\int_{\mathcal{E}}G(e)\n(\dif e)\right)\\
       &=\E\left[\exp(-\langle \Theta,F\rangle)\right],
       \end{split}
   \end{equation*}
   which implies that $\Theta^{\rho\lambda,\lambda}$ converges to $\Theta$ in distribution under $M_0$ topology.
\end{proof}

\textbf{Acknowledgment}.
The author would like to express his sincere gratitude to Roger Tribe and Oleg Zaboronski for their valuable guidance, insightful discussions, and continuous support throughout this research.

Ruibo Kou is supported by the Warwick Mathematics Institute Centre for Doctoral Training, and gratefully acknowledges funding from the University of Warwick and the UK Engineering and Physical Sciences Research Council (Grant number: EP/W524645/1 )

\textbf{AI Declaration}.
ChatGPT 5.6 plus was the only AI tool used in the preparation of this paper. It helped identify the $M_0$ topology for formulating a convergence result, for which the author had independently developed the proof without recognizing the corresponding topology. It also suggested using Faà di Bruno’s formula, given by (\ref{k-th moement of Z_s+t-Z_s}), to simplify the author’s original approach based on induction. The author verified that ($
\ref{k-th moement of Z_s+t-Z_s}$) can be applied to derive the $n$-th derivative. All mathematical estimates were derived independently by the author. Additionally, the tool was used to polish the language and improve the clarity of the exposition.

\bibliographystyle{plain} 
\bibliography{ref}

  \end{document}